\documentclass[12pt]{amsart}
\usepackage{amscd,amssymb,graphics,color,a4wide,hyperref}

\usepackage{mathrsfs}

\usepackage[all]{xy}
\newtheorem{theorem}{Theorem}[section]
\newtheorem{lemma}[theorem]{Lemma}
\newtheorem{proposition}[theorem]{Proposition}
\newtheorem{corollary}[theorem]{Corollary}

\newtheorem{statement}[theorem]{Statement}

\theoremstyle{definition}
\newtheorem{definition}[theorem]{Definition}
\newtheorem{construction}[theorem]{Construction}

\newtheorem{example}[theorem]{Example}

\theoremstyle{remark}
\newtheorem{remark}[theorem]{Remark}
\newtheorem{remarks}[theorem]{Remarks}

\numberwithin{equation}{section}
\numberwithin{figure}{section}

\newcommand{\bbfamily}{\fontencoding{U}\fontfamily{bbold}\selectfont}
\newcommand{\textbb}[1]{{\bbfamily#1}}

\newcommand {\lfor} {\mbox{\textbb{[}}}
\newcommand {\rfor} {\mbox{\textbb{]}}}

\newcommand{\NN} {\mathbb{N}}
\newcommand{\ZZ} {\mathbb{Z}}
\newcommand{\QQ} {\mathbb{Q}}
\newcommand{\RR} {\mathbb{R}}
\newcommand{\CC} {\mathbb{C}}

\newcommand{\PP} {\mathbb{P}}
\newcommand{\VV} {\mathbb{V}}

\newcommand {\shC}  {\mathcal{C}}

\newcommand {\shK}  {\mathcal{K}}

\newcommand {\shM}  {\mathcal{M}}

\newcommand {\shR}  {\mathcal{R}}

\newcommand {\shX}  {\mathcal{X}}

\newcommand {\foB}  {\mathfrak{B}}

\newcommand {\foD}  {\mathfrak{D}}

\newcommand {\foM}  {\mathfrak{M}}

\newcommand {\foU}  {\mathfrak{U}}

\newcommand {\fob}  {\mathfrak{b}}
\newcommand {\foc}  {\mathfrak{c}}
\newcommand {\fod}  {\mathfrak{d}}

\newcommand {\foj}  {\mathfrak{j}}

\newcommand {\fom}  {\mathfrak{m}}

\newcommand {\fou}  {\mathfrak{u}}
\newcommand {\fov}  {\mathfrak{v}}

\newcommand {\Ass}  {\operatorname{Ass}}
\newcommand {\Asym}  {\operatorname{Asym}}
\newcommand {\Aut}  {\operatorname{Aut}}

\newcommand{\disk}  {\mathrm{disk}}

\newcommand {\dlog} {\operatorname{dlog}}

\newcommand {\dual} {\vee}

\newcommand {\ev}  {\operatorname{ev}}

\newcommand {\Hom}  {\operatorname{Hom}}

\newcommand {\id}  {\operatorname{id}}

\newcommand {\Int}  {\operatorname{Int}}
\newcommand {\Init}  {\operatorname{Init}}

\newcommand {\Interstices} {\mathrm{Interstices}}
\newcommand {\Ideal} {\mathrm{Ideal}}
\newcommand {\Joints} {\operatorname{Joints}}

\newcommand {\kk} {\Bbbk}

\newcommand {\M} {\mathcal{M}}

\renewcommand {\max} {{\operatorname{max}}}

\newcommand {\Mono} {\operatorname{Mono}}

\newcommand {\Mult} {\operatorname{Mult}}

\renewcommand{\O}  {\mathcal{O}}

\newcommand {\out}  {\mathrm{out}}

\newcommand {\Pic}  {\operatorname{Pic}}

\newcommand {\prim}  {\mathrm{prim}}

\renewcommand {\Re} {{\operatorname{Re}}}

\newcommand {\Sing} {\operatorname{Sing}}

\newcommand {\Spec} {\operatorname{Spec}}
\newcommand {\Spf}  {\operatorname{Spf}}

\newcommand {\Supp} {\operatorname{Supp}}

\newcommand {\trop} {\mathrm{trop}}

\newcommand {\Trees}  {\operatorname{Trees}}

\newcommand {\Val}  {\operatorname{Val}}

\newcommand {\X} {\frak X}

\def\dual#1{{#1}^{\scriptscriptstyle \vee}}
\def\mapright#1{\smash{
  \mathop{\longrightarrow}\limits^{#1}}}

\def\mydate{\ifcase\month \or January\or February\or March\or
April\or May\or June\or July\or August\or September\or October\or 
November\or December\fi \space\number\day,\space\number\year}

\begin{document}

\title[Mirror symmetry for $\PP^2$ and tropical geometry]
{Mirror symmetry for $\PP^2$ and tropical geometry}

\author{Mark Gross} \address{UCSD Mathematics,
9500 Gilman Drive, La Jolla, CA 92093-0112, USA}
\email{mgross@math.ucsd.edu}
\thanks{This work was partially supported by NSF grant 0805328.}

\date{\today}

\maketitle
\tableofcontents

\section*{Introduction.}

In \cite{Announce,logmirror, part II,GS}, 
Bernd Siebert and I have been working on
a program designed to understand mirror symmetry via an algebro-geometric
analogue of the Strominger-Yau-Zaslow program \cite{SYZ}. 
The basic idea is
that the controlling objects in mirror symmetry are integral affine manifolds
with singularities. One can view an integral affine manifold as producing
a mirror pair of manifolds, one
a symplectic manifold and one a complex manifold, each 
of twice the real dimension. These correspond to the $A$-
and $B$-models of mirror symmetry. A great deal of the work carried out
by myself and Siebert has been devoted to building up a dictionary
between geometric notions on affine manifolds and 
objects in the $A$- and $B$-models. If mirror symmetry is to become
self-evident from this process, one should be able to find a single geometric
notion on an affine manifold which corresponds to both rational
curves on the $A$-model side and corrections to period calculations on
the $B$-model side. A conceptual proof of mirror symmetry would
identify these objects in the world of integral affine geometry.
(For a survey of this basic approach, the reader may consult \cite{Gr4};
however, while this paper is motivated by this program, it is largely
self-contained.)

In fact, much progress has been made in this direction. In the last five years,
it has become apparent that tropical geometry is the relevant geometry
on integral affine manifolds. In particular, thanks to work of Mikhalkin
\cite{mikhalkin} in dimension two and work of Nishinou and Siebert \cite{nisi} 
in higher dimensions, it is known that there is a correspondence between
tropical rational curves in $\RR^n$ and holomorphic rational curves in
$n$-dimensional toric varieties. On the other hand, work of 
Siebert and myself
\cite{GS} has shown that from an integral affine manifold with
singularities
and some additional data, one can construct a maximally unipotent degeneration
of varieties with effective anti-canonical class
which is controlled by tropical data on the
affine manifold. Thus one begins to see some connection between the
two sides of mirror symmetry.

While this includes the Calabi-Yau case, 
there is still a great deal of work which needs
to be done. In the $B$-model one needs to calculate periods and understand
how these periods are related to the tropical data determining the degeneration.
In the $A$-model, there is as yet no correspondence theorem between tropical
curves and actual curves. So we are still some distance from a proof of mirror
symmetry using these ideas, but the path is clear.

On the other hand, in the case of $\PP^2$, curve-counting via tropical
methods is already understood. So this raises the question: Is it 
possible, using the philosophy of my program with Siebert, to prove
mirror symmetry for $\PP^2$?

The mirror to $\PP^2$ was introduced by Givental in \cite{Giv}, in which
he explicitly described small quantum cohomology of $\PP^2$ via certain
oscillatory integrals. Using work of Dubrovin \cite{Dub} on semi-simple
Frobenius manifolds, one can see that mirror symmetry extends to
big quantum cohomology. Work of Sabbah and Sabbah-Douai \cite{Sabbah,DSI,DSII}
and Barannikov \cite{Bar} showed how to construct Frobenius manifold
structures on the full moduli space of the mirror to $\PP^2$. Here,
we are interested in Barannikov's explicit description of this Frobenius
manifold structure via oscillatory integrals to obtain a description
of \emph{big} quantum cohomology with gravitational descendents. 
The $B$-model description of \emph{small}
quantum cohomology for $\PP^2$ is simpler, and its tropical interpretation
has already been studied by Chan and Leung in \cite{ChanLeung}. 

It is clear that there should be a tropical argument for mirror
symmetry; after all, Mikhalkin
showed that big quantum cohomology for $\PP^2$ can be computed tropically, while
Barannikov showed that his $B$-model for the mirror to $\PP^2$ calculates
big quantum cohomology. The point of this paper is more philosophical.
I claim that mirror symmetry for $\PP^2$ is in fact \emph{equivalent}
to tropical genus zero curve counting in tropical $\PP^2$, in a rather
strong and transparent way. In fact, I believe the most conceptual way to
prove mirror symmetry for $\PP^2$ is via tropical geometry. The
same will be true in any dimension, and the ideas in this paper should also
generalize to other toric varieties, but for simplicity of exposition, 
and for maximum explicitness, for
the most part we will stick to $\PP^2$ in this paper.

The mirror to $X=\PP^2$ is an algebraic
torus $\check X\cong (\CC^{\times})^2$, best written as $V(x_0x_1x_2-1)
\subseteq \CC^3$. However, $\check X$ must come along with a 
\emph{Landau-Ginzburg potential}, a regular function $W:\check X\rightarrow
\CC$. It is usually written as $W=x_0+x_1+x_2$, but the full complex
moduli space of the mirror is the universal unfolding of this function $W$.
The work of Barannikov \cite{Bar} and Douai-Sabbah \cite{DSI}
generalised classic work of Kyoji Saito \cite{Saito},
showing how to construct a Frobenius manifold structure on this universal
unfolding. Barannikov described this structure explicitly in terms of
oscillatory integrals.

This procedure is quite subtle, and is explained in low-tech terms in
\S 1 of this paper. The key point is identifying \emph{flat coordinates}
on the universal unfolding
moduli space. The existence of these flat coordinates tells
us that there is a canonical way of deforming the potential $W$. This
raises the question of whether or not it is possible to write down
such a canonical deformation directly. This paper answers this question
positively.

The main idea is as follows. It has been understood since work of Cho
and Oh \cite{ChoOh} that the terms in the Landau-Ginzburg potential
are in one-to-one correspondence with families of Maslov index two
holomorphic 
disks in $\PP^2$ with boundary on a given $T^2\subseteq \PP^2$,
a fibre of the moment map on $\PP^2$ induced by the toric structure
on $\PP^2$. Therefore, a natural idea for deforming this potential is
to include some points in $\PP^2$ and consider families of disks,
again with boundary in a given $T^2$, which pass through a subset of
the given points. We deform the potential by adding terms corresponding
to families of such disks of the proper dimension.

Such a construction is in fact theoretically possible, and Fukaya,
Oh, Ohta and Ono have shown in independent work \cite{FO3II}
that this idea can be used.
The approach here, in keeping with the general
philosophy I have advocated with Bernd Siebert, is to replace holomorphic
disks with tropical disks. Once this is done, there are no theoretical
difficulties involved in defining the correct perturbation, and its
calculation is a purely combinatorial problem. This construction is
very simple, and is explained in detail in \S 2. These first two sections
can be read as an extended introduction. We also note that work
of Nishinou \cite{Nishinou} shows that holomorphic disks can be counted
tropically, so by working tropically from the beginning, we avoid a
great deal of technical difficulties related to holomorphic disks.

The main theorem of the paper is then stated at the end of \S 3.
Its main content is:
\begin{enumerate}
\item The natural parameters appearing in the deformation of $W$
constructed in \S 2 using tropical disks are in fact flat coordinates.
\item Mirror symmetry for $\PP^2$ is equivalent to specific
tropical formulas for descendent Gromov-Witten invariants 
for $\PP^2$ of the form
\[
\langle [pt],\cdots,[pt], \psi^k \alpha\rangle_d
\]
where $[pt]$ denotes the cohomology class of a point and $\alpha
\in H^*(\PP^2,\CC)$. 
\end{enumerate}

Expanding on this second point, we remark that it is standard
to put the descendent Gromov-Witten invariants together into
a generating function called the $J$-function. Similarly, once
we have tropical definitions for these descendents, we can assemble
these invariants in the same way to form the tropical $J$-function,
which we write as $J^{\trop}$. The philosophy of the paper is then
exhibited by the following diagram:
\bigskip
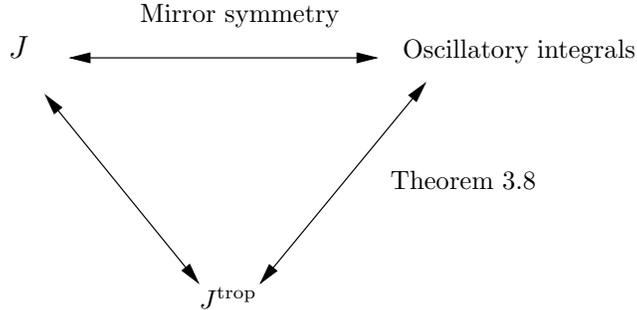
\begin{center}
\input{hannah.pstex_t}
\end{center}
Theorem \ref{maintheorem} identifies the relevant oscillatory integrals
in the B-model with $J^{\trop}$, while mirror symmetry identifies
these same integrals with $J$. Hence the equality $J=J^{\trop}$ suggested
by the third side of the triangle is
equivalent to mirror symmetry.

Note that for $k=0$, we get ordinary 
Gromov-Witten invariants for $\PP^2$, and since mirror symmetry for
$\PP^2$ is known, this gives a new proof of Mikhalkin's tropical curve
counting formula \cite{mikhalkin}.\footnote{It should not be viewed
as a surprise that this gives a new proof of Mikhalkin's formula.
Indeed, the proof of Mikhalkin's formula given by Gathmann and 
Markwig \cite{GathMark} shows that this formula gives a potential satisfying
the WDVV equation; essentially the standard proofs of mirror symmetry
for $\PP^2$ show that the potential produced on the $B$-side
also satisfies the WDVV equations.}
Furthermore, for $\alpha=[pt]$, 
these formulas agree with those recently discovered by Markwig
and Rau in \cite{MarkRau}. This latter work calculates the descendent
invariants via tropical intersection theory. At this point in time,
this method has not been extended to other choices of $\alpha$. The
problem is that tropical intersection theory becomes more difficult
when one has to worry about compactified tropical moduli spaces.
By restricting to the case where $\alpha=[pt]$, this is not an issue.
So there is yet no tropical derivation of the other descendent invariants.

\S\S4 and 5 are devoted to the proof of the main theorem and are
focused on an explicit evaluation of the necessary oscillatory
integrals. However, the crucial
point, which is explored in \S 4,
is that our deformation of the Landau-Ginzburg potential depends
on the choice of a basepoint $Q\in \RR^2$ (which we view as tropical
$\PP^2$), playing the tropical role of
a choice of fibre of the moment map. There turns out to be a chamber
structure in $\RR^2$, so that as $Q$ moves around, the deformed
Landau-Ginzburg potential changes via an explicit wall-crossing formula.
These wall-crossings can be viewed as having to do with Maslov index zero
disks, and this is the same wall-crossing phenomenon 
observed by Auroux in \cite{Auroux}.
This chamber structure is formed by what Siebert and I call a 
\emph{scattering diagram}, which played a vital role in \cite{GS} and
uses ideas originating in \cite{ks}. This is really the one point
where this paper makes contact with the ideas in \cite{GS}.

This point of contact should not be surprising. In \cite{ks} and \cite{GS},
scattering diagrams are used to build, in the former article, non-Archimedean
K3 surfaces, and in the latter article, degenerations of Calabi-Yau manifolds.
It is expected that the scattering diagrams which appear in those
papers should describe Maslov index zero disks on the mirror side. So
again, wall crossing formulas are associated with the presence of
Maslov index zero disks.

In \S 5, we complete the proof by evaluating the necessary period
integrals. This is a rather involved and delicate process, but the calculation
shows that we indeed have constructed flat coordinates and produced
the tropical formulas for the descendent invariants mentioned above.

\smallskip

As mentioned above, there is some overlap between this work and work
of Fukaya, Oh, Ohta and Ono \cite{FO3II}. 
Both Fukaya and I spoke on preliminary versions
of our work at the final conference for the program in Real and Tropical
Algebraic Geometry at EPFL in Lausanne in June 2008.
At that time, I had not yet determined the formula for the descendent
invariants $\langle [pt],\ldots,[pt],\psi^k [\PP^2]\rangle$. While
our approaches are quite different, and Fukaya et al's work does not deal
with descendents, his talk helped lead me to
the correct formulation, for which I thank him.

I would also like to thank Paul Hacking, Claus Hertling,
Ludmil Katzarkov, Sean Keel, Hannah Markwig, and Bernd
Siebert for useful conversations, and thank Kwokwai Chan and
Conan Leung for explaining their work to me during a very pleasant 
visit to CUHK,
which motivated me to work on this problem. I would also like to thank
Alexander Givental, who pointed out to me in 2004 that there was little
hope of using tropical geometry to understand mirror symmetry unless
one could find tropical descriptions of gravitational descendents.
Finally, I would like to thank the referees for useful suggestions.

\section{Barannikov's mirror symmetry for $\PP^n$}

We begin by reviewing Barannikov's description of the $B$-model
for $\PP^n$ in terms of oscillatory integrals, 
giving a precise statement of mirror symmetry.
For the purposes of this paper, we do not need the general
formalism of semi-infinite variation of Hodge structures, but instead
explain the approach as concretely as possible.

Let $X=\PP^n$, $\check X=V(\prod_{i=0}^n x_i-1)\subseteq \Spec \CC[x_0,
\ldots,x_n]$.
Here $\check X$ is isomorphic to $(\CC^{\times})^n$,
but we should consider $\check X$ as a Landau-Ginzburg model, with
potential 
\[
W_0=\sum_{i=0}^n x_i.
\]
The pair $(\check X,W_0)$ is usually viewed as the mirror
to $X$.
We then consider the universal unfolding of $W_0$, with ${\bf t}=(t_0,
\ldots,t_n)$,
\[
W_{{\bf t}}=W_0+\sum_{j=0}^nt_jW_0^j,
\]
parameterized by the moduli space $\shM=\Spf \CC\lfor t_0,\ldots,t_n\rfor$, 
the completion at the origin
of $\CC^{n+1}$ with coordinates $t_0,\ldots,t_n$. Here $\Spf$ denotes
the formal spectrum.
One then considers the local system $\shR$ on $\shM\times \CC^{\times}$
whose fibre at a point $({\bf t},q)$ is the relative homology group
$H_n(\check X,\Re(q W_{{\bf t}})\ll 0;\CC)$.\footnote{While this group is
often written in this way, it should be defined more precisely as
the space of rapid decay homology cycles \cite{Hien},\cite{KaKoPa} as follows.
Choose a variety $Y$ containing $\check X$ with $Y\setminus\check X$
a normal crossings divisor, such that the map $W:\check X\rightarrow
\CC$ extends to a map $W:Y\rightarrow \PP^1$ with $Y\setminus \check X$
mapping to $\infty\in\PP^1$. Let $\tilde Y$ denote a real oriented blow-up
of $Y\setminus\check X$ in $Y$. The exceptional locus contains
a set $Z$ consisting of all points $b$ such that $\Re (qW(z))\rightarrow
-\infty$ as $z\rightarrow b$. We then define the homology group as the
relative homology $H_n(\tilde Y,Z; \CC)$. The main point is that
these are precisely the cycles over which it makes sense to integrate forms
of the sort appearing in (M2).
}

Then Barannikov shows,\footnote{
The discussion in \cite{Bar} considers $\shM$ to be an analytic germ
of $0\in \CC^{n+1}$. This raises certain technical issues, because
for deformations in the directions $t_i$, $i\ge 2$,
we obtain ``non-tame''
behaviour, and this relative homology group jumps, because the critical
locus of the Landau-Ginzburg potential $W_{{\bf t}}$ jumps.
We get around this problem by working formally around a neighbourhood
where the relevant homology group is the correct one. It is not difficult
to check that \cite{Bar} works in this context.}
first of all,
that one can find a unique choice of the following data:
\begin{itemize}
\item[(M1)]
A (multi-valued) basis
of sections of $\shR$, $\Xi_0,\ldots,\Xi_n$, with $\Xi_i$
uniquely determined modulo $\Xi_0,\ldots,\Xi_{i-1}$.
\item[(M2)]
A section $s$ of $\dual{\shR}\otimes_{\CC} \O_{\shM\times \CC^{\times}}$
defined by integration of a family of holomorphic forms on
$\check X\times\shM\times \CC^{\times}$
of the form 
\[
e^{qW_{{\bf t}}}f \dlog x_1\wedge\cdots\wedge\dlog x_n,
\]
where $q$ is the coordinate on $\CC^{\times}$ and 
$f$ is a regular function on $\check X\times\shM\times\CC^{\times}$
with $f|_{\check X\times \{{\bf 0}\}\times \CC^{\times}}=1$ and
which extends to a regular function on $\check X\times\shM\times(\CC^{\times}
\cup \{\infty\})$.
\end{itemize}
This data must satisfy the following conditions:
\begin{itemize}
\item[(M3)] The monodromy associated to $q\mapsto qe^{2\pi i}$ in the local system
$\shR$ is given, in the basis $\Xi_0,\ldots,\Xi_n$, by
$\exp((n+1)2\pi iN)$, where
\[
N=\begin{pmatrix} 0&1&0&0&\ldots&0\\
0&0&1&0&\ldots&0\\
&\vdots&&&&\vdots\\
0&0&0&0&\ldots&1\\
0&0&0&0&\ldots&0
\end{pmatrix}.
\]
Note that this condition determines the basis $\Xi_0,
\ldots,\Xi_n$ up to a change of basis matrix $S=(s_{ij})$ an
upper triangular matrix with $s_{i-1,j-1}=s_{ij}$ for $1\le i,j\le n$.
\item[(M4)] We identify a fibre of the dual local system $\dual{\shR}$
with the ring $\CC[\alpha]/(\alpha^{n+1})$, with $\alpha^i$ dual
to $\Xi_i$. Note that $H^{2*}(\PP^n,\CC)\cong \CC[\alpha]/
(\alpha^{n+1})$, with the primitive positive generator of $H^2(\PP^n,\CC)$
corresponding to $\alpha$. Under this isomorphism, the action of $(n+1)N^{t}$
can be viewed as mirror to cupping with
the anti-canonical class in $H^{2*}(\PP^n,\CC)$.

The section $s$ of $\dual{\shR}\otimes_{\CC}\O_{\shM\times\CC^{\times}}$
yields an element of each fibre of $\dual{\shR}$, 
which we can write as
\[
s({\bf t},q)=\sum_{i=0}^n \alpha^i\int_{\Xi_i} 
e^{qW_{{\bf t}}}f\dlog x_1\wedge\cdots\wedge\dlog x_n.
\]
We then require that we can write
\[
s({\bf t},q)=q^{(n+1)\alpha}\sum_{i=0}^n\varphi_i({\bf t},q)(\alpha/q)^i
\]
for functions $\varphi_i$ satisfying
\[
\varphi_i({\bf t},q)=\delta_{0,i}+\sum_{j=1}^{\infty} \varphi_{i,j}({\bf t})q^j.
\]
for $0\le i\le n$.
(The shape of this formula determines the section $s$
uniquely; Barannikov refers to this as a normalization condition.) 
Here we use the expansion
\[
q^{(n+1)\alpha}=\sum_{i=0}^n {(n+1)^i\over i!}(\log q)^i \alpha^i
\]
to intepret $q^{(n+1)\alpha}$; this takes care of the multi-valuedness
of the integrals. A consequence of these
two conditions is that if we set 
\[
y_i({\bf t})=\varphi_{i,1}({\bf t}),\quad 0\le i \le n,
\]
$y_0,\ldots,y_n$ form a system of coordinates on $\shM$, which are
called \emph{flat coordinates}.
Furthermore,
\[
\lim_{q\rightarrow 0}{\varphi_{i}({\bf 0},q)\over q^i}=\delta_{0,i}.
\]
(This last condition, not explicitly mentioned in \cite{Bar}, 
fixes the basis $\Xi_0,\ldots,\Xi_n$ up to the action of
monodromy.)
\end{itemize}

We will take the existence of such data satisfying these properties
as given, and in some sense
our goal will be to identify
the flat coordinates $y_0,\ldots,y_n$ and the regular
function $f$ tropically, and compare the
results with the $A$-model data on $\PP^n$. Before discussing the
$A$-model data, it is worth noting that it is quite non-trivial to
find $y_0,\ldots,y_n$ and $f$ in terms of ${\bf t}$, and to
the best of my knowledge, unlike in the Calabi-Yau case,
this computation has not been carried out in the literature even to
low order in the $t_i$'s. Rather, Barannikov proves mirror symmetry
for $\PP^n$ by showing that the formalism of semi-infinite variation
of Hodge structures allows one to construct a Frobenius manifold
structure on $\shM$ for which the $y_i$'s are flat coordinates.
The $A$-model for $\PP^n$ also yields a Frobenius manifold, and the fact
that these two Frobenius manifolds are isomorphic is shown first by 
identifying the two algebra structures at one point, and then using
semi-simplicity of this Frobenius algebra structure and results of Dubrovin
\cite{Dub}.

We now consider the $A$-model for $\PP^n$, so we can state mirror
symmetry for $\PP^n$. Let $\overline{\shM}_{0,m}(\PP^n,d)$ denote
the moduli space of stable maps of degree $d$ from $m$-pointed curves
of genus zero into $\PP^n$. Let $\psi_1,\ldots,\psi_m\in
H^2(\overline{\shM}_{0,m}(\PP^n,d),\QQ)$ be the usual $\psi$ classes,
i.e., $\psi_i$ is the first Chern class of the line bundle on
$\overline{\shM}_{0,m}(\PP^n,d)$ whose fibre at a point $[f,C,p_1,\ldots,p_m]$
is the cotangent line to $C$ at the marked point $p_i$. We have evaluation maps
$\ev_i:\overline{\shM}_{0,m}(\PP^n,d)\rightarrow\PP^n$, and define,
for classes $\beta_1,\ldots,\beta_m\in H^*(\PP^n,\QQ)$, the descendent
Gromov-Witten invariant
\[
\langle \psi^{\nu_1}\beta_1,\ldots,\psi^{\nu_m}\beta_m\rangle_{d}
=\int_{[\overline{\shM}_{0,m}(\PP^n,d)]}
\bigwedge_{i=1}^m (\psi_i^{\nu_i}\wedge \ev_i^*\beta_i).
\]
We can then write the precise statement of mirror symmetry in terms
of the Givental $J$-function. This is a function $J_{\PP^n}(y_0,
\ldots,y_n,q)$ with values in $H^*(\PP^n,\CC)$.
Let $T_0,\ldots,T_n$ be generators of $H^*(\PP^n,\ZZ)$, with $T_i\in
H^{2i}(\PP^n,\ZZ)$ positive. Then with $\gamma=\sum_{i=2}^n y_iT_i$,
the $J$-function is defined by
\[
J_{\PP^n}(y_0,\ldots,y_n,q)
=e^{q(y_0T_0+y_1T_1)}\bigg(T_0+
\sum_{i=0}^n\sum_{m\ge 0}\sum_{d\ge 0} {1\over m!} \langle T_0,\gamma^m, T_{n-i}
/(q^{-1}-\psi)
\rangle_d e^{dy_1}T_i\bigg).
\]
(See e.g.,
\cite{iritani}, Definition 2.14 for this description of the $J$-function.)
Here $\gamma^m$ means we take $\gamma$ $m$ times, and $1/(q^{-1}-\psi)
=q/(1-q\psi)$
is expanded formally in $q\psi$.
We define $J_i$ by writing 
\[
J_{\PP^n}(y_0,\ldots,y_n,q)=\sum_{i=0}^n J_{i}(y_0,\ldots,y_n,q)T_i.
\]
We then consider the following statement:

\begin{statement}[Mirror symmetry for $\PP^n$]
\label{mirrorstatement}
In the
$\CC$-vector space $\CC\lfor y_0,\ldots,y_n,q\rfor$,
\[
J_i=\varphi_i.
\]
\end{statement}

This mirror symmetry statement was proved by Barannikov in \cite{Bar} 
for the part of the statement which does not involve gravitational
descendents. In any event, genus zero descendent invariants can
be reconstructed from the non-descendent invariants, but see
\cite{Iritani} for a more direct proof for the statement with
gravitational descendents. However, the philosophy in this paper is
not to prove it, but to prove its equivalence to a tropical statement
which is stated precisely in Statement \ref{tropicalstatement}.

\medskip

It is worthwhile expanding out the expression for the $J$-function.
Recall by the Fundamental Class
Axiom (see e.g.\ \cite{CoxKatz}, page 305) that if $d>0$,
we have
\[
\langle T_0,T_0^{m_0},\ldots,T_n^{m_n},\psi^\nu T_{n-i}\rangle_d
=
\langle T_0^{m_0},\ldots,T_n^{m_n},\psi^{\nu-1} T_{n-i}\rangle_d
\]
where a correlator involving $\psi^{-1}$ is interpreted as zero. 
Also note that $\langle T_0,\psi^{\nu}T_{n-i}\rangle_0=0$ since
$\overline\M_{0,2}(\PP^n,0)$ is empty.
Thus, 
\begin{eqnarray*}
J_{\PP^n}&=&e^{q(y_0T_0+y_1T_1)}\bigg(
T_0+\sum_{i=0}^n\sum_{m_2+\cdots+m_n\ge 0}\sum_{d,\nu\ge 0} \langle 
T_0,T_2^{m_2},\ldots,T_n^{m_n},\psi^\nu T_{n-i}\rangle_d
\,q^{\nu+1}e^{dy_1}{y_2^{m_2}\cdots y_n^{m_n}\over m_2!\cdots m_n!}T_i\bigg)\\
&=&e^{q(y_0T_0+y_1T_1)}\bigg(
T_0+\sum_{i=0}^n\bigg(
\sum_{m_2+\cdots+m_n\ge 1}\sum_{\nu\ge 0}\langle T_0,T_2^{m_2},
\ldots,T_n^{m_n},\psi^{\nu}T_{n-i}\rangle_0
\,q^{\nu+1}{y_2^{m_2}\cdots y_n^{m_n}\over m_2!\cdots m_n!}\\
&&\quad\quad\quad+
\sum_{m_2+\cdots+m_n\ge 0}\sum_{d\ge 1}\sum_{\nu\ge 0}
\langle T_2^{m_2},\ldots,T_n^{m_n},\psi^{\nu} T_{n-i}\rangle_d
\,q^{\nu+2}e^{dy_1}{y_2^{m_2}\cdots y_n^{m_n}\over m_2!\cdots m_n!}\bigg)T_i\bigg).
\end{eqnarray*}

Explicitly for the case of interest in this paper, namely $\PP^2$,
we note that 
\[
\langle T_0,T_2^m,\psi^{\nu}T_{2-i}\rangle_0=0
\]
unless $m+1=2m+2-i+\nu$, for dimension reasons. Thus we need
$0\le\nu=-m-1+i$, so we only get a contribution when $m=1$,
$i=2$, as $i\le 2$. But 
$\langle T_0,T_2,T_{2-i}\rangle_0=\int_{\PP^2} T_0\cup T_2\cup T_{2-i}
=\delta_{2,i}$. Thus we get
\[
J_{\PP^2}
=e^{q(y_0T_0+y_1T_1)}\bigg(T_0+\sum_{i=0}^2\bigg(qy_2\delta_{2,i}
+
\sum_{m_2\ge 0}\sum_{d\ge 1}\sum_{\nu\ge 0}
\langle T_2^{m_2},\psi^{\nu} T_{2-i}\rangle_d
\,q^{\nu+2}
e^{dy_1}{y_2^{m_2}\over m_2!}\bigg)T_i\bigg).
\]

\section{Tropical geometry}

We next review the definition of a tropical curve from
\cite{mikhalkin}, and then
introduce the notion of a tropical disk.

We fix once and for all a
lattice $M=\ZZ^n$, $N=\Hom_{\ZZ}(M,\ZZ)$ the dual lattice,
$M_{\RR}=M\otimes_{\ZZ} \RR$, $N_{\RR}=N\otimes_{\ZZ}\RR$.

Let $\overline{\Gamma}$ be a weighted, connected
finite graph without bivalent vertices. Its set of vertices and
edges are denoted $\overline{\Gamma}^{[0]}$ and $\overline{\Gamma}^{[1]}$,
respectively, and $w_{\overline{\Gamma}}:\overline{\Gamma}^{[1]}
\rightarrow \NN=\{0,1,\ldots\}$ is the weight function. An edge
$E\in \overline{\Gamma}^{[1]}$ has adjacent vertices $\partial E=\{V_1,V_2\}$.
Let $\overline{\Gamma}_{\infty}^{[0]}\subseteq \overline{\Gamma}^{[0]}$
be the set of univalent vertices. We set
\[
\Gamma:=\overline{\Gamma}\setminus\overline{\Gamma}_{\infty}^{[0]}.
\]
We write the set of vertices and edges of $\Gamma$ as $\Gamma^{[0]}$,
$\Gamma^{[1]}$, and we have the weight function $w_{\Gamma}:\Gamma^{[1]}
\rightarrow \NN$.
Some edges of $\Gamma$ are now non-compact, and these are called
\emph{unbounded edges}. We use the convention that the weights of unbounded
edges are always zero or one, and the weights of all bounded edges are
positive. In particular, we do not allow bounded edges to be
contracted in what follows
(as is sometimes the case in the tropical geometry literature).
Write $\Gamma_{\infty}^{[1]}\subseteq
\Gamma^{[1]}$ for the set of unbounded edges.

\begin{definition}
\label{tropcurve}
A \emph{parameterized $d$-pointed marked tropical curve} 
in $M_{\RR}$ with marked points $\{p_1,\ldots,p_d\}$ is a choice of inclusion
$\{p_1,\ldots,p_d\}\hookrightarrow \Gamma_{\infty}^{[1]}$ written as
$p_i\mapsto E_{p_i}$, and a continuous map
$h:\Gamma\rightarrow M_{\RR}$ satisfying the following conditions.
\begin{enumerate}
\item $w_{\Gamma}(E)=0$ if and only if $E=E_{p_i}$ for some $i$.
\item $h|_{E_{p_i}}$ is constant, $1\le i\le d$, while for every other edge 
$E\in\Gamma^{[1]}$, the restriction $h|_E$ is
a proper embedding with image $h(E)$ contained in an affine line with rational
slope.
\item For every vertex $V\in\Gamma^{[0]}$, the following \emph{balancing
condition} holds. Let $E_1,\ldots,E_m\in \Gamma^{[1]}$ be the edges
adjacent to $V$, and let $m_i\in M$ be the primitive integral vector
emanating from $h(V)$ in the direction of $h(E_i)$. Then
\[
\sum_{j=1}^m w_{\Gamma}(E_j)m_j=0.
\]
\end{enumerate}

We write a parameterized $d$-pointed tropical curve as
\[
h:(\Gamma,p_1,\ldots,p_d)\rightarrow M_{\RR}.
\]
We write $h(p_i)$ for the point $h(E_{p_i})$.

An \emph{isomorphism} of tropical curves $h_1:(\Gamma_1,p_1,
\ldots,p_d)\rightarrow M_{\RR}$
and $h_2:(\Gamma_2,p_1,
\ldots,p_d)\rightarrow M_{\RR}$ is a homeomorphism $\Phi:\Gamma_1
\rightarrow\Gamma_2$ respecting the marked edges and 
the weights with
$h_1=h_2\circ\Phi$. A \emph{$d$-pointed tropical curve} is an isomorphism class
of parameterized $d$-pointed tropical curves. We never distinguish
between a $d$-pointed tropical curve and a particular representative.

The \emph{genus} of a tropical curve
$h:\Gamma\rightarrow M_{\RR}$ is the first Betti number of $\Gamma$.
A \emph{rational} tropical curve is a tropical curve of genus zero.

The \emph{combinatorial type} of a marked tropical curve
$h:(\Gamma,p_1,\ldots,p_d)\rightarrow M_{\RR}$ 
is defined to be the homeomorphism
class of $\overline{\Gamma}$ with the marked points and weights, together
with, for every vertex $V$ and edge $E$ containing $V$, the
primitive tangent vector to $h(E)$ in $M$ pointing away from $V$.
\end{definition}

We modify this definition slightly to define a \emph{tropical disk}:
(see \cite{Nishinou}, where these are called tropical curves with
stops).

\begin{definition} Let $\overline{\Gamma}$ be a weighted, connected finite
graph without bivalent vertices as above, with the additional data
of a choice of univalent vertex $V_{\out}$, adjacent to a unique edge
$E_{\out}$. Let 
\[
\Gamma':=(\overline{\Gamma}\setminus\overline{\Gamma}_{\infty}^{[0]})
\cup\{V_{\out}\}\subseteq\overline{\Gamma}.
\]
Suppose furthermore that $\Gamma'$ has first Betti number zero (i.e., 
$\Gamma'$ is a tree with one compact external edge and a number of
non-compact external edges). Then a \emph{parameterized $d$-pointed 
tropical disk}
in $M_{\RR}$ is a choice of inclusion $\{p_1,\ldots,p_d\}\hookrightarrow 
\Gamma^{[1]}_{\infty}\setminus\{E_{\out}\}$ written as $p_i
\mapsto E_{p_i}$ and a
map $h:\Gamma'\rightarrow M_{\RR}$ satisfying
the same conditions as Definition \ref{tropcurve},
except there is \emph{no} balancing condition at $V_{\out}$.

An \emph{isomorphism} of tropical disks $h_1:(\Gamma_1',
p_1,\ldots,p_d)\rightarrow M_{\RR}$
and $h_2:(\Gamma'_2,
p_1,\ldots,p_d)\rightarrow M_{\RR}$ is a homeomorphism $\Phi:\Gamma_1'
\rightarrow\Gamma_2'$ respecting the marked edges and the 
weights with
$h_1=h_2\circ\Phi$. A \emph{tropical disk} is an isomorphism class
of parameterized tropical disks.

The \emph{combinatorial type} of a tropical disk
$h:(\Gamma',p_1,\ldots,p_d)\rightarrow M_{\RR}$ 
is defined to be the homeomorphism
class of $\overline{\Gamma}$ with the marked points, weights, and $V_{\out}$, 
together
with, for every vertex $V$ and edge $E$ containing $V$, the
primitive tangent vector to $h(E)$ in $M$ pointing away from $V$.
\end{definition}

We also recall Mikhalkin's notion of multiplicity. 
For the remainder of the paper, we
restrict to the case that $M$ is rank two, i.e., $M=\ZZ^2$. Much of
what we say can be generalized to higher dimension, but for ease of
exposition, we restrict to dimension two.

\begin{definition}
\label{multiplicitydef}
Let $h:\Gamma\rightarrow M_{\RR}$ be a marked tropical
curve or $h:\Gamma'\rightarrow M_{\RR}$ be a marked tropical disk
such that $\bar\Gamma$ only has vertices of valency one and
three. The \emph{multiplicity} of a vertex $V\in\Gamma^{[0]}$ in $h$ is
$\Mult_V(h)=1$ if one of the edges adjacent to $V$ has weight zero (i.e,
is a marked unbounded edge), and otherwise
\[
\Mult_V(h)=w_1w_2|m_1\wedge m_2|=w_1w_3|m_1\wedge m_3|=w_2w_3|m_2\wedge
m_3|,
\]
where $E_1,E_2,E_3\in\Gamma^{[1]}$ are the edges containing $V$,
$w_i=w_{\Gamma}(E_i)$, and $m_i\in M$ is a primitive integral vector
emanating from $h(V)$ in the direction of $h(E_i)$. The equality of the
three expressions follows from the balancing condition.

The \emph{multiplicity} of the curve or disk $h$ is then
\[
\Mult(h):=\prod_{V\in\Gamma^{[0]}}\Mult_V(h).
\]
Note that in the case of the tropical disk, there is no contribution
from $V_{\out}$.
\end{definition}

We now fix a complete rational polyhedral fan $\Sigma$ in $M_{\RR}$,
with $\Sigma^{[1]}$ denoting the set of one-dimensional cones in $\Sigma$.
We denote by $T_{\Sigma}$ the free abelian group generated by $\Sigma^{[1]}$,
and for $\rho\in\Sigma^{[1]}$ we denote by $t_{\rho}$ the corresponding
generator for $T_{\Sigma}$ and by $m_{\rho}$ the primitive generator of
the ray $\rho$. Let $X_{\Sigma}$ denote the toric surface defined by
$\Sigma$. 

\begin{definition} 
A tropical curve (or disk) $h$ is a \emph{tropical curve (or disk) in
$X_{\Sigma}$} if every $E\in\Gamma^{[1]}_{\infty}$ (or
$E\in \Gamma^{[1]}_{\infty}\setminus \{E_{\out}\}$) has $h(E)$ 
either a point or
a translate of some $\rho\in\Sigma^{[1]}$.

If the tropical curve or disk has $d_{\rho}$ unbounded edges which
are translates of $\rho\in\Sigma^{[1]}$, (remember these unbounded
edges always have weight one), then
the \emph{degree} of $h$ is 
\[
\Delta(h):=\sum_{\rho\in\Sigma^{[1]}}d_{\rho}t_{\rho}\in T_{\Sigma}.
\]
We define
\[
|\Delta(h)|:=\sum_{\rho\in\Sigma^{[1]}}d_{\rho}.
\]
\end{definition}

\bigskip

Fix points $P_1,\ldots,P_k\in M_{\RR}$ general, and fix a general 
\emph{base-point}
$Q\in M_{\RR}$. When we talk about general points in the sequel,
we mean that there is an open dense subset (typically the complement
of a finite union of polyhedra of codimension at least one) of $M_{\RR}^{k+1}$
such that $(P_1,\ldots,P_k,Q)\in M_{\RR}^{k+1}$ lies in this open subset.
This choice of open subset will depend on particular needs.

Associate to the points $P_1,\ldots,P_k$ the variables
$u_1,\ldots,u_k$ in the ring
\[
R_k:={\CC[u_1,\ldots,u_k]\over (u_1^2,\ldots,u_k^2)}.
\]

\begin{definition} 
Let $h:(\Gamma',p_1,\ldots,p_d)\rightarrow M_{\RR}$ be a
tropical disk in $X_{\Sigma}$ with $h(V_{\out})=Q$,
$h(p_j)=P_{i_j}$, $1\le i_1<\cdots<i_d\le k$. (This ordering
removes a $d!$ ambiguity about the labelling of the marked points.)
We say $h$ is a \emph{tropical disk in $(X_{\Sigma},P_1,
\ldots,P_k)$ with boundary $Q$}.

The \emph{Maslov index} of the disk $h$ is
\[
MI(h):=2(|\Delta(h)|-d).
\]
\end{definition}

\begin{lemma}
\label{diskmoduli}
If $P_1,\ldots,P_k,Q$ are chosen in general position,
then the set of Maslov index $2n$ tropical disks in $(X_{\Sigma},P_1,\ldots,
P_k)$ with boundary $Q$
is an $(n-1)$-dimensional polyhedral complex. The set of Maslov
index $2n$ tropical disks with arbitrary boundary is an
$(n+1)$-dimensional polyhedral complex.
\end{lemma}

\proof This is a standard tropical general position argument. We sketch
it here. Fix a combinatorial type of tropical disk with $d$ marked points,
with degree $\Delta$.
If the combinatorial type is general, then the domain
$\Gamma'$ only has trivalent vertices apart from $V_{\out}$. Such a tree
has $|\Delta|+d-1$ bounded edges (including $E_{\out}$). A tropical disk 
$h:\Gamma'\rightarrow M_{\RR}$ of this given combinatorial type
is then completely determined by the position of $h(V_{\out})\in M_{\RR}$
and the affine lengths of the bounded edges. This produces a cell
in the moduli space $\shM^{\disk}_{\Delta,d}(X_{\Sigma})$
of all $d$-pointed tropical disks of degree 
$\Delta$.
The closure of this cell is $(\RR_{\ge 0})^{|\Delta(h)|+d-1}
\times M_{\RR}$. Also, there are only a finite number
of combinatorial types of disks of a given degree.
Thus $\shM^{\disk}_{\Delta,d}
(X_{\Sigma})$ is a finite 
$(|\Delta|+d+1)$-dimensional polyhedral complex. 
Furthermore, we have a piecewise
linear map $\ev:\shM^{\disk}_{\Delta,d}(X_{\Sigma})\rightarrow M_{\RR}^d$,
taking a disk $h$ to the tuple $(h(p_1),\ldots,h(p_d))$.
Let $E\subseteq\shM^{\disk}_{\Delta,d}(X_{\Sigma})$ be the union of
cells mapping under $\ev$ to cells of codimension $\ge 1$ in $M_{\RR}^d$;
then $h(E)$ is a closed subset of $M_{\RR}^d$. Thus, 
if $(P_{i_1},\ldots,P_{i_d})
\in M_{\RR}^d$ is not in this closed subset, for $1\le i_1<\ldots<i_d\le k$
distinct indices, then $\ev^{-1}
(P_{i_1},\ldots,P_{i_d})$ is a codimension $2d$ subset of 
$\shM^{\disk}_{\Delta,d}(X_{\Sigma})$. Thus the dimension of the
moduli space of
tropical disks of a given degree $\Delta$ with arbitrary boundary in
$(X_{\Sigma},P_1,\ldots,P_k)$ is $|\Delta|+1-d=MI(h)/2+1$. Similarly, if we fix
a general boundary point $Q$, the dimension is $MI(h)/2-1$, as claimed.
\qed

\bigskip

\begin{definition}
\label{LGdef}
Given the data $P_1,\ldots,P_k,Q\in M_{\RR}$ general,
let $h:(\Gamma',p_1,\ldots,p_d)\rightarrow M_{\RR}$ be a Maslov
index two marked
tropical disk with boundary $Q$ in $(X_{\Sigma},
P_1,\ldots,P_k)$.
Then we can associate to $h$
a monomial in $\CC[T_{\Sigma}]\otimes_{\CC} R_k\lfor y_0\rfor$,
\[
\Mono(h):=\Mult(h)z^{\Delta(h)} 
u_{I(h)}, 
\]
where $z^{\Delta(h)}\in\CC[T_{\Sigma}]$ is the monomial corresponding to
$\Delta(h)\in T_{\Sigma}$, the subset
$I(h)\subseteq \{1,\ldots,k\}$ is defined by
\[
I(h):=\{i\,|\,\hbox{$h(p_j)=P_i$ for some $j$}\},
\]
and 
\[
u_{I(h)}=\prod_{i\in I(h)} u_i.
\]
Define the \emph{$k$-pointed Landau-Ginzburg potential}
\[
W_k(Q):=y_0+\sum_h \Mono(h)\in \CC[T_{\Sigma}]\otimes_{\CC} R_k\lfor y_0\rfor
\]
where the sum is over all Maslov index two disks $h$ in 
$(X_{\Sigma},P_1,\ldots,P_k)$ with boundary $Q$. By Lemma \ref{diskmoduli},
this is a finite sum for $P_1,\ldots,P_k,Q$ general.
\end{definition}

Now restrict further to the case that $X_{\Sigma}$ is non-singular.
There is an obvious map $r:T_{\Sigma}\rightarrow M$
given by $r(t_{\rho})=m_{\rho}$, the primitive generator of $\rho$,
and the assumption of non-singularity gives $r$ surjective.
So there is a natural exact sequence
\[
0\mapright{} K_{\Sigma}\mapright{} T_{\Sigma}\mapright{r} M\mapright{} 0
\]
defining $K_{\Sigma}$. Dualizing this sequence gives
\[
0\rightarrow N\rightarrow \Hom_{\ZZ}(T_{\Sigma},\ZZ)\rightarrow \Pic X_{\Sigma}
\rightarrow 0.
\]
See, e.g., \cite{Fulton}, \S 3.4, for this description of the
Picard group of $X_{\Sigma}$.
After tensoring with $\CC^{\times}$, we get an exact sequence
\[
0\mapright{} N\otimes \CC^{\times} \mapright{}
\Hom(T_{\Sigma},\CC^{\times})\mapright{\kappa}\Pic X_{\Sigma}\otimes\CC^{\times}
\mapright{} 0.
\]
We set 
\[
\check\shX_{\Sigma}:=\Hom(T_{\Sigma},\CC^{\times})
=\Spec \CC[T_{\Sigma}]
\]
and we define the \emph{K\"ahler moduli space}\footnote{This is almost,
but not quite, the usual K\"ahler moduli space in the context of mirror
symmetry (see, e.g., \cite{CoxKatz},\S6.2). The K\"ahler moduli space of a
K\"ahler manifold $X$
is generally the tube domain given by $(H^2(X,\RR)+
i\shK)/H^2(X,\ZZ)$ for $\shK\subseteq H^2(X,\RR)$ the K\"ahler cone of
$X$. In the case $X=X_{\Sigma}$, this is naturally
an analytic open subset of $\shM_{\Sigma}$.}
of $X_{\Sigma}$ to be
\[
\shM_{\Sigma}:=\Pic X_{\Sigma}\otimes\CC^{\times}=\Spec \CC[K_{\Sigma}].
\]
We also have a morphism
\[
\kappa:\check\shX_{\Sigma}\rightarrow\shM_{\Sigma}.
\]
We thicken $\shM_{\Sigma}$ by setting 
\[
\shM_{\Sigma,k}=\shM_{\Sigma}\times\Spf R_k\lfor y_0\rfor,
\]
and 
\[
\check\shX_{\Sigma,k}=\check\shX_{\Sigma}\times \Spf R_k\lfor y_0\rfor.
\]
Then we have the family 
\[
\kappa:\check\shX_{\Sigma,k}\rightarrow\shM_{\Sigma,k}.
\]
$W_k(Q)$ is a regular function on $\check\shX_{\Sigma,k}$,
so we can think of this as providing a family of Landau-Ginzburg
potentials. Note that a fibre of $\kappa$ over a closed point
of $\shM_{\Sigma,k}$ is isomorphic to $N\otimes\CC^{\times}$.

The sheaf of relative differentials 
$\Omega^1_{\check\shX_{\Sigma,k}/\M_{\Sigma,k}}$ is canonically
isomorphic to $M\otimes_{\ZZ}\O_{\check \shX_{\Sigma,k}}$, with
$m\otimes 1$ corresponding to the differential
\begin{equation}
\label{dlogdef}
\dlog m:={d(z^{\overline{m}})\over z^{\overline{m}}};
\end{equation}
here $\overline{m}\in T_{\Sigma}$ is any lift of $m\in M$, and $\dlog m$
is well-defined as a \emph{relative} differential independently
of the lift. Thus a choice of generator of $\bigwedge^2 M\cong
\ZZ$ determines a nowhere-vanishing relative holomorphic two-form $\Omega$,
canonical up to sign. Explicitly, if $e_1,e_2\in M$ is a
positively oriented basis, then
\begin{equation}
\label{Omegadef} \Omega=\dlog e_1\wedge \dlog e_2.
\end{equation}

\begin{remark}
\label{correctLGremark}
The function $W_k(Q)$ is intended to be the
``correct'' Landau-Ginzburg potential to describe the mirror
to $X_{\Sigma}$, in the sense that $W_0(Q)$ is the expression
usually taken to be the Landau-Ginzburg potential, and
$W_k(Q)$ should be a canonical perturbation, in the sense
that the parameters appearing in $W_k(Q)$ are closely related
to flat coordinates. However, for general choice of $\Sigma$
this is not true, the chief problem being that there will be
copies of $\PP^1$ in the toric boundary of $X_{\Sigma}$
which do not deform to curves intersecting the big torus orbit
of $X_{\Sigma}$. This is a standard problem in tropical geometry:
so far, tropical geometry cannot ``see'' these curves. This
is not a problem as long as $X_{\Sigma}$ is a product of projective
spaces, so in particular, we will now restrict to the case of
$X_{\Sigma}=\PP^2$.
\end{remark}

\begin{example}
\label{firstP2example}
Let $\Sigma$ be
the fan depicted in Figure \ref{P2fanfigure}, so that
$X_{\Sigma}=\PP^2$.
Here $T_{\Sigma}=\ZZ^3$ with basis $t_0,t_1,t_2$
corresponding to $\rho_0,\rho_1,\rho_2$,
and we write $x_i$ for the monomial $z^{t_i}\in\kk[T_{\Sigma}]$. The map
\[
\kappa:\Spec \kk[T_{\Sigma}]\rightarrow \shM_{\Sigma}=\Spec \kk[K_{\Sigma}]
\]
is then a map
\[
\kappa:(\kk^{\times})^3\rightarrow\kk^{\times}
\]
given by $\kappa(x_0,x_1,x_2)=x_0x_1x_2$.
\begin{figure}
\input{P2fan.pstex_t}
\caption{The fan for $\PP^2$.}
\label{P2fanfigure}
\end{figure}

If we take $k=0$, then there are precisely three Maslov index two
tropical disks, as depicted in Figure \ref{NoPointsDisks}. Thus
we take 
\[
W_0(Q)=y_0+x_0+x_1+x_2.
\]
This is the standard Landau-Ginzburg potential for the mirror
to $\PP^2$ (except for the additional variable $y_0$). The
formula in Definition \ref{LGdef} gives a deformation of this
potential over the thickened moduli space $\shM_{\Sigma,k}$.

If we take $k=1$, marking one point in $\PP^2$, we obtain one
additional disk, as depicted in Figure \ref{OnePointDisks},
and if we take $k=2$ with $P_1$ and $P_2$ chosen as in Figure
\ref{TwoPointDisks}, we have three additional disks.
Note the potential
depends on the particular choices of the points $P_1,\ldots,P_k$
as well as $Q$. In the given examples, we have respectively
\begin{eqnarray*}
W_1(Q)&=&y_0+x_0+x_1+x_2+u_1x_1x_2\\
W_2(Q)&=&y_0+x_0+x_1+x_2+u_1x_0x_1+u_2x_0x_1+u_1u_2x_0x_1^2.
\end{eqnarray*}
\qed

\begin{figure}
\input{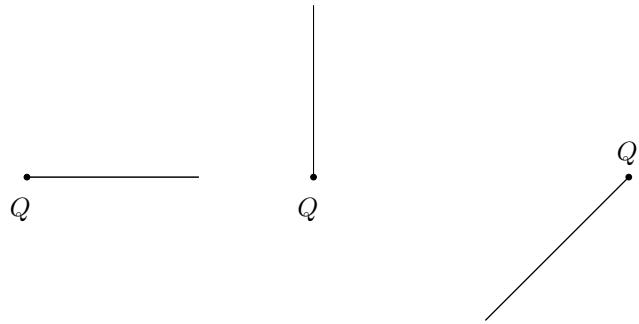}
\caption{Maslov index two tropical disks with no marked points.}
\label{NoPointsDisks}
\end{figure}
\begin{figure}
\input{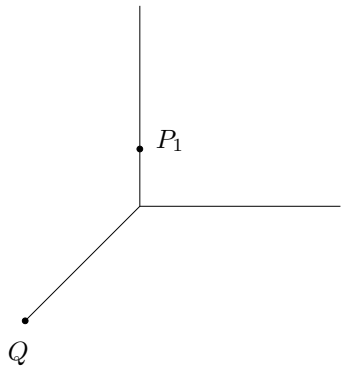}
\caption{The one additional 
Maslov index two tropical disk with $k=1$.}
\label{OnePointDisks}
\end{figure}
\begin{figure}
\input{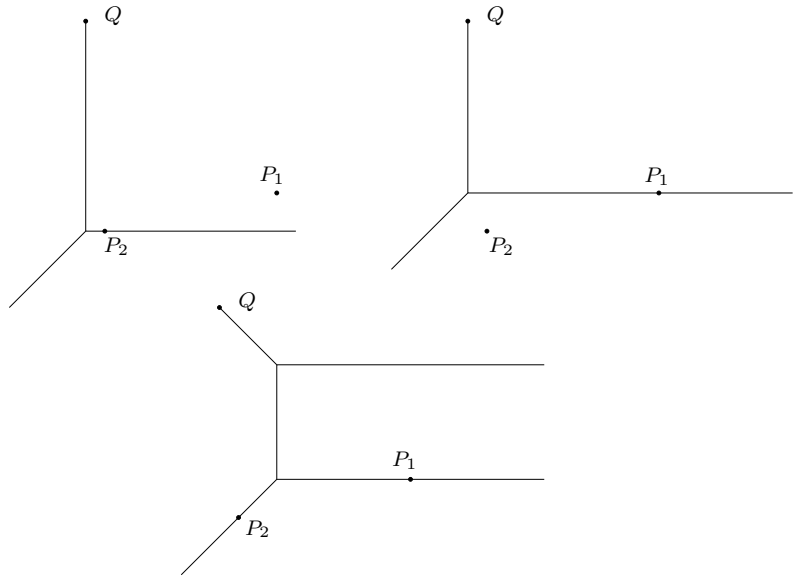}
\caption{The additional Maslov index two tropical disks with $k=2$.}
\label{TwoPointDisks}
\end{figure}
\end{example}

\bigskip

\section{Tropical descendent invariants and the main theorem}

We will now define the tropical version of the descendent
invariants and state the main theorem of the paper. 
These invariants will be defined as a count of
certain tropical curves with vertices of high valency.
These curves need to be counted with certain multiplicities. 
We begin with

\begin{definition}
\label{psiclassdef}
Let $P_1,\ldots,P_k\in M_{\RR}$ be general.
Let $\Sigma$ be a complete fan in $M_{\RR}$ defining a toric surface.
Let $S\subseteq M_{\RR}$ be a subset.
Define
\[
\shM^{\trop}_{\Delta,n}(X_{\Sigma},P_1,\ldots,P_k,\psi^{\nu}S)
\]
to be the moduli space of rational $(n+1)$-pointed tropical
curves in $X_{\Sigma}$
\[
h:(\Gamma,p_1,\ldots,p_n,x)\rightarrow M_{\RR}
\]
of degree $\Delta$ such that
\begin{enumerate}
\item $h(p_j)=P_{i_j}$, $1\le i_1<\cdots<i_n\le k$.
\item The edge $E_x$ is attached to a vertex $V_x$ of $\Gamma$; the
valency of this vertex is denoted $\Val(x)$. Then 
\[
\Val(x)=\nu+3.
\]
\item $h(x)\in S$.
\end{enumerate}
\end{definition}

\begin{lemma}
For $P_1,\ldots,P_k\in M_{\RR}$ general, 
\begin{enumerate}
\item
$\shM^{\trop}_{\Delta,n}(X_{\Sigma},P_1,\ldots,P_k,\psi^{\nu}M_{\RR})$
is a polyhedral complex of dimension $|\Delta|-n-\nu$.
\item 
$\shM^{\trop}_{\Delta,n}(X_{\Sigma},P_1,\ldots,P_k,\psi^{\nu}C)$
is a polyhedral complex of dimension $|\Delta|-n-\nu-1$ for $C$
a general translate of a tropical curve in $M_{\RR}$.
\item
$\shM^{\trop}_{\Delta,n}(X_{\Sigma},P_1,\ldots,P_k,\psi^{\nu}Q)$
is a polyhedral complex of dimension $|\Delta|-n-\nu-2$
for $Q\in M_{\RR}$ a general point.
\end{enumerate}
\end{lemma}

\proof This is again standard, as in Lemma \ref{diskmoduli}. The
dimension count is as follows. Fix the combinatorial type of the
curve to be generic, so that all vertices of $\Gamma$ are trivalent except
for the vertex adjacent to $E_x$, which is $(\nu+3)$-valent. Such
a tree has $|\Delta|+n+1$ unbounded edges, and then has
$|\Delta|+n+1-(\nu+3)$ bounded edges. The curves of this combinatorial
type are then determined by the location of $h(x)\in M_{\RR}$ and
the lengths of the bounded edges, giving a cell of the form
$(\RR_{\ge 0})^{|\Delta|+n-\nu-2}\times M_{\RR}$. Fixing $h(p_1),\ldots,
h(p_n)$ then yields the desired dimension of the moduli space in (1)
being
\[
|\Delta|+n-\nu-2n. 
\]
This gives (1). For (2) and (3), we consider the map
\[
\ev_x:\shM^{\trop}_{\Delta,n}(X_{\Sigma},P_1,\ldots,P_k,\psi^{\nu}M_{\RR})
\rightarrow M_{\RR}
\]
given by $\ev_x(h):=h(x)$. Let $E_1$ be the union of
cells in $\shM^{\trop}_{\Delta,n}(X_{\Sigma},P_1,\ldots,P_k,\psi^{\nu}
M_{\RR})$ which map to codimension $\ge 1$ sets in $M_{\RR}$, and
let $E_2$ be the union of cells which map to points
in $M_{\RR}$. Then we need to choose the translate $C$ so that
$C\cap h(E_1)$ is zero dimensional and $C\cap h(E_2)=\emptyset$. 
Similarly, we need to choose $Q\not\in h(E_1)$. Then $\ev_x^{-1}(C)$
or $\ev_x^{-1}(Q)$ are the desired moduli spaces in cases (2) and
(3) and are of the desired dimension.
\qed

\begin{lemma}
\label{curvechop}
Let $P_1,\ldots,P_k\in M_{\RR}$ be general and $S\subseteq M_{\RR}$
a subset.
Let
\[
h\in\shM^{\trop}_{\Delta,n}(X_{\Sigma},P_1,\ldots,P_k,\psi^{\nu}S).
\]
Let $\Gamma_1',\ldots,\Gamma_{\nu+2}'$
denote the closures of the connected components of $\Gamma\setminus
\{E_x\}$, with $h_i:\Gamma_i'\rightarrow M_{\RR}$ the restrictions of 
$h$. Each disk $h_i$ is viewed as being marked by those points $p\in
\{p_1,
\ldots,p_n\}$ with $E_p\subseteq \Gamma_i'$. 
There is one special case to consider here: if $E_x$ and $E_{p_i}$
share a common vertex $V$, then we discard the edge $E_{p_i}$ from
consideration as well, so we have disks $h_1,\ldots,h_{\nu+1}$. 
(Note that since $h(p_i)\not=h(p_j)$ for $i\not=j$, 
we never have $E_{p_i}$ and $E_{p_j}$ sharing a common vertex.)
\begin{enumerate}
\item If $S=M_{\RR}$ and $n=|\Delta|-\nu$, then either
\begin{enumerate}
\item $E_x$ does not share a vertex with any of the edges $E_{p_i}$, and
then $MI(h_i)=2$ for all but two
choices of $i$, for which $MI(h_i)=0$.
\item $E_x$ does share a vertex with one of the edges $E_{p_i}$,
and then $MI(h_i)=2$ for all $i$.
\end{enumerate}
\item If $S=C$ is a general translate of a tropical curve in $M_{\RR}$ 
and $n=|\Delta|-\nu-1$, 
then $MI(h_i)=2$ for all but one
$i$, and for this $i$, $MI(h_i)=0$.
\item If  $S=\{Q\}$ for a general point $Q$ and $n=|\Delta|-\nu-2$,
then $MI(h_i)=2$ for all $i$.
\end{enumerate}
\end{lemma}

\proof First note that the condition on $n$ and the generality of
$P_1,\ldots,P_k$, $C$, and $Q$ guarantee by the previous
lemma that the moduli space under consideration is zero dimensional.
If any of the disks $h_i$ can be deformed
while keeping its boundary $h_i(x)$ fixed, then this yields a non-trivial
deformation of $h$, which does not exist. Thus by
Lemma \ref{diskmoduli} we must have $MI(h_i)\le 2$ in all cases.
Let $n_i$ be the number of marked points on $h_i$.
We note that 
\begin{eqnarray*}
\sum_i {MI(h_i)\over 2}&=&\sum_i (|\Delta(h_i)|-n_i)\\
&=&\begin{cases}
|\Delta(h)|-(n-1)&\hbox{Case (1) (b)}\\
|\Delta(h)|-n&\hbox{otherwise}\end{cases}\\
&=&\begin{cases}
\nu&\hbox{Case (1) (a)}\\
\nu+1&\hbox{Case (1) (b)}\\
\nu+1&\hbox{Case (2)}\\
\nu+2&\hbox{Case (3)}\\
\end{cases}
\end{eqnarray*}
Since there are $\nu+2$ disks except in Case (1) (b),
when there are $\nu+1$ disks, the result follows.
\qed

\medskip

We can now define tropical analogues of the descendent Gromov-Witten invariants
which appear in the Givental $J$-function. From now on in this section, 
we are only
considering the case of $X_{\Sigma}=\PP^2$, with $\Sigma$ the fan
with rays generated by $m_0=(-1,-1)$, $m_1=(1,0)$ and $m_2=(0,1)$,
and $t_0,t_1,t_2$ the generators of $T_{\Sigma}$, with $r(t_i)=m_i$.
Let 
\[
\Delta_d:=d(t_0+t_1+t_2)\in T_{\Sigma};
\]
curves of degree $\Delta_d$ should be viewed as degree $d$ curves
in $\PP^2$.

\begin{definition} 
\label{tropicaldescinv}
Fix general points $Q,P_1,P_2,\cdots\in M_{\RR}$. Let
$L$ be the tropical line (a translate of the fan of $\PP^2$)
with vertex $Q$.

For a tropical curve $h$ in $\PP^2$ with a marked point $x$, 
let $n_0(x),n_1(x)$ and $n_2(x)$ be the number of unbounded rays
sharing a common vertex with $E_x$ 
in the directions $m_0,m_1$ and $m_2$ respectively.
As in Lemma
\ref{curvechop}, we denote by $h_1,\ldots$ the tropical
disks obtained by removing $E_x$ from $\Gamma$, with the outgoing edge of
$h_i$ being $E_{i,\out}$. Let $m(h_i)=w(E_{i,\out})m^{\prim}(h_i)$,
where $m^{\prim}(h_i)\in M$ is a primitive vector tangent to
$h_i(E_{i,\out})$ pointing away from $h(x)$.

Define
\begin{eqnarray*}
\Mult^0_x(h)&=&{1\over n_0(x)!n_1(x)!n_2(x)!}\\
\Mult^1_x(h)&=&
-{\sum_{k=1}^{n_0(x)} {1\over k}+\sum_{k=1}^{n_1(x)} {1\over k}
+\sum_{k=1}^{n_2(x)} {1\over k}\over n_0(x)! n_1(x)! n_2(x)!}\\
\Mult^2_x(h)&=&
{\bigg(\sum_{l=0}^2\sum_{k=1}^{n_l(x)} {1\over k}\bigg)^2+\sum_{l=0}^2
\sum_{k=1}^{n_l(x)} {1\over k^2}
\over 2n_0(x)! n_1(x)! n_2(x)!}
\end{eqnarray*}
\begin{enumerate}
\item
We define 
\[
\langle P_1,\ldots,P_{3d-2-\nu},\psi^{\nu} Q\rangle_d^{\trop}
\]
to be 
\[
\sum_{h} \Mult(h)
\]
where the sum is over all marked tropical rational curves 
\[
h\in \shM^{\trop}_{\Delta_d,3d-2-\nu}(P_1,\ldots,P_{3d-2-\nu},\psi^{\nu}Q).
\]
We define
\[
\Mult(h):=\Mult^0_x(h)\prod_{V\in \Gamma^{[0]}\atop V\not\in E_x} \Mult_V(h).
\] 
\item
We define 
\[
\langle P_1,\ldots,P_{3d-1-\nu},\psi^{\nu} L\rangle_d^{\trop}
\]
as a sum
\[
\sum_{h} \Mult(h)
\]
where the sum is again over all marked tropical rational
curves  
\[
h:(\Gamma,p_1,\ldots,p_{3d-1-\nu},x)\rightarrow M_{\RR}
\]
with $h(p_i)=P_i$ and
satisfying one of the following two conditions.
\begin{enumerate}
\item 
\[
h\in \shM^{\trop}_{\Delta_d,3d-1-\nu}(P_1,\ldots,P_{3d-1-\nu},\psi^{\nu}L).
\]
Furthermore, no unbounded edge of $\Gamma$ having a common vertex
with $E_x$ other than $E_x$
maps into the connected component of $L\setminus \{Q\}$ containing
$h(x)$. By Lemma \ref{curvechop}, there is precisely one $j$, $1\le j\le
\nu+2$, with $MI(h_j)=0$. Suppose also 
the connected component of $L\setminus\{Q\}$ containing $h(x)$
is $Q+\RR_{>0} m_i$.
Then we define
\[
\Mult(h)=|m(h_{j})\wedge m_i|
\Mult^0_x(h)\prod_{V\in\Gamma^{[0]}\atop V\not\in E_x}\Mult_V(h).
\]
\item $\nu\ge 1$ and
\[
h\in \shM^{\trop}_{\Delta_d,3d-1-\nu}(P_1,\ldots,P_{3d-1-\nu},\psi^{\nu-1}Q).
\]
Then
\[
\Mult(h)=\Mult_x^1(h)\prod_{V\in\Gamma^{[0]}\atop V\not\in E_x}\Mult_V(h).
\]
\end{enumerate}
\item
We define 
\[
\langle P_1,\ldots,P_{3d-\nu},\psi^{\nu} M_{\RR}\rangle_d^{\trop}
\]
as a sum
\[
\sum_h \Mult(h)
\]
where the sum is over all marked tropical rational curves
\[
h:(\Gamma,p_1,\ldots,p_{3d-\nu},x)\rightarrow M_{\RR}
\]
such that $h(p_i)=P_i$ and either
\begin{enumerate}
\item 
\[
h\in \shM^{\trop}_{\Delta_d,3d-\nu}(P_1,\ldots,P_{3d-\nu},\psi^{\nu}M_{\RR})
\]
and $E_x$ does not share a vertex with any of the $E_{p_i}$'s.
Furthermore, no unbounded
edge of $\Gamma$ having a common vertex with $E_x$ other than
$E_x$ maps into the connected
component of $M_{\RR}\setminus L$ containing $h(x)$.
By Lemma \ref{curvechop}, there are precisely two distinct
$j_1,j_2$ with $1\le j_1,j_2\le \nu+2$ such that $MI(h_{j_i})=0$.
Then we define
\[
\Mult(h)=
|m(h_{j_1})\wedge m(h_{j_2})|
\Mult^0_x(h)\prod_{V\in\Gamma^{[0]}\atop V\not\in E_{x}} \Mult_V(h).
\]
\item 
\[
h\in \shM^{\trop}_{\Delta_d,3d-\nu}(P_1,\ldots,P_{3d-\nu},\psi^{\nu}M_{\RR})
\]
and $E_x$ shares a vertex with $E_{p_i}$.
Furthermore, no unbounded edge of $\Gamma$ having a common vertex with
$E_x$ other than $E_x$ and $E_{p_i}$
maps into the connected component of $M_{\RR}\setminus L$ containing $h(x)$.
Then we define
\[
\Mult(h)=\Mult^0_x(h)\prod_{V\in\Gamma^{[0]}\atop V\not\in E_x} \Mult_V(h).
\]
\item $\nu\ge 1$ and
\[
h\in \shM^{\trop}_{\Delta_d,3d-\nu}(P_1,\ldots,P_{3d-\nu},\psi^{\nu-1}L).
\]
Furthermore, no unbounded edge of $\Gamma$ having a common vertex with
$E_x$ other than $E_x$ maps into the connected component of $L\setminus\{Q\}$
containing $h(x)$. By Lemma \ref{curvechop}, there is precisely one
$j$, $1\le j\le \nu+1$, with $MI(h_j)=0$. Suppose the connected
component of $L\setminus \{Q\}$ containing $h(x)$
is $Q+\RR_{>0} m_i$.
Then we define
\[
\Mult(h)=|m(h_{j})\wedge m_i|\Mult^1_x(h)\prod_{V\in\Gamma^{[0]}\atop V\not
\in E_x}\Mult_V(h).
\]
\item $\nu\ge 2$ and
\[
h\in \shM^{\trop}_{\Delta_d,3d-\nu}(P_1,\ldots,P_{3d-\nu},\psi^{\nu-2}Q).
\]
Then
\[
\Mult(h)=\Mult^2_x(h)\prod_{V\in \Gamma^{[0]}\atop V\not\in E_x} \Mult_V(h).
\]
\end{enumerate}
\end{enumerate}

In all cases $S=\{Q\}, L$ or $M_{\RR}$, we define for $\sigma\in
\Sigma$,
\begin{equation}
\label{tropdesczonedef}
\langle P_1,\ldots,P_{3d-\nu-(2-\dim S)},\psi^{\nu} S\rangle_{d,\sigma}^{\trop}
\end{equation}
to be the contribution to 
$\langle P_1,\ldots,P_{3d-\nu-(2-\dim S)},\psi^{\nu} S\rangle_d^{\trop}$
coming from curves $h$ with $h(x)$ in the interior of $\sigma+Q$.
In (1), the only contribution comes from $\sigma=\{0\}$, in (2),
the contributions come from the zero and one-dimensional cones of
$\Sigma$, and in (3), the contributions come from all cones of $\Sigma$.
\end{definition}

\begin{remarks}
(1) Note that all moduli spaces involved are zero-dimensional
for general choices of $Q,P_1,\ldots$, so the sums make sense.

(2) The formula in Definition \ref{tropicaldescinv}, (1), for $\nu=0$,
gives the standard tropical curve counting formula for the number
of rational curves of degree $d$ passing through $3d-1$ points.
For $\nu>0$, this coincides with the formula given by Markwig and Rau
in \cite{MarkRau}. In particular, by the results of that paper,
\[
\langle P_1,\ldots,P_{3d-2-\nu},\psi^{\nu}Q\rangle^{\trop}_d
=\langle T_2^{3d-2-\nu},\psi^{\nu}T_2\rangle_d.
\]

(3) Clearly the formulas for the descendent invariants involving
$\psi^{\nu}L$ or $\psi^{\nu}M_{\RR}$ are rather more complicated
and mysterious. Hannah Markwig informs me that if the tropical
intersection procedure
of \cite{MarkRau} is carried out for the
classes $\psi^{\nu}L$ or $\psi^{\nu}M_{\RR}$, one obtains the contributions
in Definition \ref{tropicaldescinv} only of the type (2) (a) and (3)
(a) and (b). 
This does not give the correct formula for the descendent invariants.
The likely explanation for this phenomenon is that \cite{MarkRau} works
with the non-compact moduli space of tropical rational 
curves in tropical $\PP^2$. There should be a suitable compactification
of this moduli space, and it is possible that tropical intersection theory
applied to this compactified moduli space would yield the correct 
descendent invariants. If this is the case, then the contributions of the
form (2) (b) or (3) (c)-(d) in Definition \ref{tropicaldescinv}
could be viewed as boundary contributions.

It would be very interesting to learn if this is indeed the case. However,
exploring this question would be considerably outside the scope of
this paper. As a consequence, it is difficult to motivate the rather
mysterious multiplicity formulas given in Definition \ref{tropicaldescinv}.
These are simply the formulas which emerge naturally from period
integrals in \S 5.

(4) It is easy to see that
$\langle P_1,\ldots,P_{3d-1},\psi^0L\rangle_d^{\trop}$
is $d$ times the number of rational curves through $3d-1$ points.
Indeed, the only contribution to this number comes from 
Definition \ref{tropicaldescinv}, (2) (a). For each tropical rational
curve $h:\Gamma\rightarrow M_{\RR}$ with $3d-1$ marked points
passing through $P_1,\ldots,P_{3d-1}$
we obtain a contribution
for every point of $h^{-1}(L)$ by marking that point with $x$. 
The factor $|m(h_j)\wedge m_i|
\Mult^0_x(h)$ for the multiplicity
in this case gives the intersection multiplicity of $h(\Gamma)$ with $L$ at each
point of $h^{-1}(L)$,
as defined in \cite{Sturmfels}, \S 4. By the tropical B\'ezout
theorem (\cite{Sturmfels}, Theorem 4.2), the total contribution 
from $h$ is then $(h(\Gamma).L) \Mult(h)=d\Mult(h)$.

Thus
\[
\langle P_1,\ldots,P_{3d-1},\psi^{0}L\rangle^{\trop}_d
=d\langle P_1,\ldots,P_{3d-1}\rangle_d^{\trop}
=d\langle T_2^{3d-1}\rangle_d=\langle T_2^{3d-1},T_1\rangle_d,
\]
by Mikhalkin's formula and the Divisor Axiom.

(5) $\langle P_1,\ldots,P_{3d},\psi^0 M_{\RR}\rangle^{\trop}_d=0$.
Indeed, the only possible contributions come from Definition
\ref{tropicaldescinv} (3) (a), but there are no rational curves
of degree $d$ through $3d$ general points. Thus
\[
\langle P_1,\ldots,P_{3d},\psi^{0}M_{\RR}\rangle^{\trop}_d
=\langle T_2^{3d},T_0\rangle_d,
\]
as both are zero.
\end{remarks}

One of our results will be the following:

\begin{theorem}
\label{PQinvariant}
The invariants defined in Definition \ref{tropicaldescinv}
are independent of the choice of the $P_i$'s and $Q$.
\end{theorem}

This will be proved in \S 5. However, this allows us to make the following
definition.

\begin{definition}
We define
\begin{eqnarray*}
\langle T_2^{3d-2-\nu},\psi^{\nu}T_2\rangle^{\trop}_d
&:=&\langle P_1,\ldots,P_{3d-2-\nu},\psi^{\nu}Q\rangle^{\trop}_d\\
\langle T_2^{3d-1-\nu},\psi^{\nu}T_1\rangle^{\trop}_d
&:=&\langle P_1,\ldots,P_{3d-1-\nu},\psi^{\nu}L\rangle^{\trop}_d\\
\langle T_2^{3d-\nu},\psi^{\nu}T_0\rangle^{\trop}_d
&:=&\langle P_1,\ldots,P_{3d-\nu},\psi^{\nu}M_{\RR}\rangle^{\trop}_d
\end{eqnarray*}
where the $P_i$'s and $Q$ have been chosen generally.

We define
\[
\langle T_2^m,\psi^{\nu} T_i\rangle_d^{\trop}:=0
\]
if $m+i+\nu\not= 3d$.

We define the \emph{tropical $J$-function} for $\PP^2$ by
\begin{eqnarray*}
J^{\trop}_{\PP^2}&:=&e^{q(y_0T_0+y_1T_1)}\bigg(T_0+\sum_{i=0}^2\bigg(qy_2\delta_{i,2}\\
&&\quad\quad\quad+
\sum_{d\ge 1}\sum_{\nu\ge 0}
\langle T_2^{3d+i-2-\nu},\psi^{\nu} T_{2-i}\rangle^{\trop}_d
\,q^{\nu+2}e^{dy_1}{y_2^{3d+i-2-\nu}\over (3d+i-2-\nu)!}\bigg)T_i\bigg)\\
&=:&\sum_{i=0}^2 J_i^{\trop} T_i.
\end{eqnarray*}
\end{definition}

The main theorem of the paper, to be proved in \S 5, is now

\begin{theorem}
\label{maintheorem}
A choice of general points $P_1,\ldots,P_k$ and $Q$ gives rise to
a function $W_k(Q)\in \CC[T_{\Sigma}]\otimes_{\CC} R_k\lfor y_0\rfor$
by Definition \ref{LGdef}, and hence a family of Landau-Ginzburg
potentials on the family $\kappa:\check\shX_{\Sigma,k}
\rightarrow\shM_{\Sigma,k}$
with a relative nowhere-vanishing two-form $\Omega$ as defined in 
\eqref{Omegadef}.
This data gives rise to a local system $\shR$
on $\M_{\Sigma,k}\otimes \Spec\CC[q,q^{-1}]$ whose fibre over $(\kappa,q)$
is $H_2((\check\shX_{\Sigma,k})_{\kappa},\Re(qW_0(Q))\ll 0)$. 
There exists a multi-valued basis $\Xi_0,\Xi_1,\Xi_2$ 
of sections of $\shR$ satisfying conditions (M1) and (M3) of \S 1
such that 
\begin{equation}
\label{phidefeq}
\sum_{i=0}^2 \alpha^i \int_{\Xi_i} e^{qW_k(Q)}\Omega
=q^{3\alpha}\sum_{i=0}^2\varphi_i\cdot(\alpha/q)^i
\end{equation}
where
\[
\varphi_i(y_0,\kappa,u_1,\ldots,u_k)
=\delta_{0,i}+\sum_{j=1}^{\infty}\varphi_{i,j}(y_0,\kappa,u_1,\ldots,u_k)q^j,
\]
for $0\le i\le 2$,
with
\begin{eqnarray*}
\varphi_{0,1}&=&y_0\\
\varphi_{1,1}&=&y_1:=\log(\kappa)\\
\varphi_{2,1}&=&y_2:=\sum_{i=1}^k u_i.
\end{eqnarray*}
Furthermore, 
\[
\varphi_i=J^{\trop}_i(y_0,y_1,y_2).
\]
\end{theorem}

Note one immediate corollary of this theorem:

\begin{corollary} Let $\foM_{\Sigma,k}$ be the formal spectrum of the
completion of $\CC[K_{\Sigma}]\otimes_{\CC} R_k\lfor y_0\rfor$ at the
maximal ideal $(y_0,\kappa-1,u_1,\ldots,u_k)$. This completion
is isomorphic to $\CC\lfor y_0,y_1\rfor\otimes_{\CC} R_k$, with
$y_1:=\log\kappa$, the latter expanded in a power series at $\kappa=1$.
Let 
\[
\check\X_{\Sigma,k}=\check\shX_{\Sigma,k}\times_{\M_{\Sigma,k}}
\foM_{\Sigma,k}.
\]
The regular function $W_k(Q)$ on $\check\X_{\Sigma,k}$ restricts to 
$W_0(Q)=x_0+x_1+x_2$ on the closed fibre 
of $\check\X_{\Sigma,k}\rightarrow\foM_{\Sigma,k}$, and hence $W_k(Q)$ 
provides a deformation of this function over $\foM_{\Sigma,k}$.
Thus we have a morphism from $\foM_{\Sigma,k}$
to the universal unfolding moduli space $\Spec \CC\lfor y_0,y_1,y_2\rfor$
with $y_0,y_1,y_2$ the flat coordinates of \S 1, (M4). This map is
given by
\begin{eqnarray*}
y_0&\mapsto&y_0\\
y_1&\mapsto&\log\kappa\\
y_2&\mapsto&\sum_{i} u_i.
\end{eqnarray*}
\end{corollary}

\proof The universal unfolding of $x_0+x_1+x_2$ on $(\CC^\times)^2
=V(x_0x_1x_2-1)\subseteq\CC^3$ is a regular function $W$ on $(\CC^\times)^2
\times \Spec\CC\lfor y_0,y_1,y_2\rfor$, and there is a map
$\xi:\foM_{\Sigma,k}\rightarrow\Spf\CC\lfor y_0,y_1,y_2\rfor$ and
an isomorphism 
\[
\eta:\check\X_{\Sigma,k}\rightarrow(\CC^{\times})^2\times\foM_{\Sigma,k}
\]
which is the identity on the closed fibre and such that, if $W'=W\circ(\id
\times \xi)$ is the pull-back of $W$ to $(\CC^{\times})^2\times
\foM_{\Sigma,k}$, then $W_k(Q)=W'\circ\eta$.
Then
\[
(\eta^{-1})^*(e^{qW_k(Q)}\Omega)=e^{qW'}(\eta^{-1})^*(\Omega).
\]
Of course $(\eta^{-1})^*\Omega=f\Omega$ for a regular function $f$
which is identically $1$ on the closed fibre. Thus $(\eta^{-1})^*(e^{qW_k(Q)}
\Omega)$ is a form of the sort allowed in (M2), and Theorem \ref{maintheorem}
shows this form satisfies the normalization conditions of (M4). Thus 
$\xi$ is given by $y_0\mapsto y_0$, $y_1\mapsto y_1$, $y_2\mapsto
\sum_{i=1}^k u_i$ because of the calculation of $\varphi_{i,1}$, $i=0,1,2$
given in Theorem \ref{maintheorem}.
\qed

\bigskip

Another corollary is the strong equivalence between mirror symmetry
for $\PP^2$ and tropical curve counting with descendents. 
We have the following statement:

\begin{statement}[Tropical curve counting with descendents]
\label{tropicalstatement}
\[
J_{\PP^2}=J^{\trop}_{\PP^2}.
\]
\end{statement}

Then

\begin{corollary}
Statement \ref{mirrorstatement} and Statement \ref{tropicalstatement}
are equivalent. In particular, since Statement \ref{mirrorstatement}
is known to be true, Statement \ref{tropicalstatement} is true.
\end{corollary}

This is the promised strong equivalence between mirror symmetry for
$\PP^2$ and tropical curve counting.

The proofs of Theorems \ref{PQinvariant} and 
\ref{maintheorem} require explicit evaluation of the integrals
$\int_{\Xi} e^{qW_k(Q)}\Omega$ for various $\Xi$. While it is not,
in general, difficult to get explicit answers for these integrals,
it is actually quite difficult to extract a useful combinatorial
result from these answers. There is one particular case where this
is not difficult, however, namely the case when $\Xi=\Xi_0$. 
Because of the simplicity of this case, in which one can see exactly
how the oscillatory integral extracts Gromov-Witten invariants, I believe
it is worth presenting this calculation before embarking on the main part
of the proof of the main theorems.

We will need the following lemma, whose proof is given at the very
end of \S 4. 

\begin{lemma} 
\label{integralindependentlemma}
For $\Xi\in 
H_2((\check\shX_{\Sigma,k})_{\kappa},\Re(qW_0(Q))\ll 0)$,
the integral
\[
\int_{\Xi} e^{qW_k(Q)}\Omega
\]
is independent of the choice of $Q$ and $P_1,\ldots,P_k$.
\end{lemma}

\begin{proposition}
\label{J0case}
$\langle P_1,\ldots,P_{3d-2-\nu},\psi^{\nu}Q\rangle_d^{\trop}$ is
independent of the choice of the points $P_1,\ldots,P_{3d-2-\nu},Q$ general.
Furthermore, there exists a well-defined cycle $\Xi_0$ such that
\[
\int_{\Xi_0}e^{qW_k(Q)}\Omega=J_0^{\trop}.
\]
\end{proposition}

\proof
Each fibre of $\kappa$ is an algebraic torus
isomorphic to $N\otimes\CC^\times$. We will take 
$\Xi_0$ to be the element of $H_2(N\otimes \CC^{\times},\CC)$
which is $1/(2\pi i)^2$ times the homology class of the compact
torus $N\otimes U(1)\subseteq N\otimes\CC^{\times}$. Since the choice
of $\Omega$ depended on a choice of orientation on $M_{\RR}$, hence
an orientation on $N_{\RR}$, this defines an orientation on
$N\otimes U(1)$.
Then it is very easy to compute the
integral
\begin{equation}
\label{Delta0int}
\int_{\Xi_0} e^{qW_k(Q)}\Omega
=e^{qy_0}\int_{\Xi_0} e^{q(W_k(Q)-y_0)}\Omega.
\end{equation}
Indeed, this is just a residue calculation. We expand
$e^{q(W_k(Q)-y_0)}\Omega$ in a power series, obtaining a sum
of terms of the form 
$Cz^m\Omega$ for various $m\in T_{\Sigma}$, $C\in R_k[q]$. Such
a term contributes to the
integral if and only if the image $r(m)$ of $m$ in $M$ is zero, i.e., 
$m=d(t_0+t_1+t_2)=\Delta_d$ for some $d\ge 0$,
in which case the contribution is just $C\kappa^d$.

These terms are easily interpreted. We write
\[
\exp(q(W_k(Q)-y_0))=\sum_{n=0}^{\infty} {q^n(W_k(Q)-y_0)^n\over n!}.
\]
The coefficient of $q^0$ in this expansion is $1$, contributing $1$
to the integral. The coefficient of $q$ is $W_k(Q)-y_0$, 
and no term of this expression
is of the form $Cz^m$ with $r(m)=0$. So there is no contribution
from these terms. Thus any other possible contribution comes from
terms in $q^{\nu+2}(W_k(Q)-y_0)^{\nu+2}/(\nu+2)!$ for $\nu\ge 0$, 
and such a term can be written as
\[
q^{\nu+2}\prod_i {1\over \nu_i!}\Mono(h_i)^{\nu_i}
\]
with $\sum \nu_i=\nu+2$, with the $h_i:\Gamma_i'\rightarrow M_{\RR}$ 
being Maslov index two tropical disks.
Note that $\Mono(h_i)^2=0$ unless $h_i$ has no marked points, so $\nu_i=1$
unless $h_i$ consists of just one unbounded edge. 

The condition that this term contributes to the integral is saying
precisely that these tropical disks can be glued to obtain a tropical
curve. Indeed, we can take $\nu_1$ copies of $\Gamma_1'$, $\nu_2$ copies
of $\Gamma_2'$, and so on, and glue these graphs together by identifying
the outgoing vertices $V_{\out}$ on all of them to get a graph
$\Gamma$. 
We can add an additional marked unbounded edge to $\Gamma$
at $V_{\out}$, labelled $x$, so that now
$V_{\out}$ is a $\nu+3$-valent vertex of $\Gamma$. 
The maps $h_i$ then glue to give a map $h:\Gamma\rightarrow M_{\RR}$ 
with $h(x)=Q$.
This map satisfies the balancing condition at all vertices
except perhaps at $V_{\out}$.
To check balancing at $V_{\out}$, let $E_{\out,i}$ 
denote the edge of the tropical disk $h_i$ adjacent
to the outgoing vertex $V_{\out,i}$ of $\Gamma'_i$. Suppose $h_i$ has degree
$\Delta(h_i)$. 
Let $v_i\in M$ be a primitive tangent vector based at $Q$ pointing away
from $Q$ and tangent to $h_i(E_{\out,i})$. Then by summing the balancing
conditions at all vertices of $\Gamma_i'$ other than $V_{\out,i}$, one
obtains
\begin{equation}
\label{balanceeq1}
w_{\Gamma_i}(E_{\out,i})v_i=r(\Delta(h_i)).
\end{equation}
So the condition that $\prod_i \Mono(h_i)^{\nu_i}$ gives a non-zero
contribution to the integral \eqref{Delta0int} is equivalent to
\begin{equation}
\label{balanceeq2}
0=\sum_{i}\nu_i r(\Delta(h_i))=\sum_i \nu_iw_{\Gamma_i}(E_{\out,i})v_i,
\end{equation}
which is precisely the balancing condition for $\Gamma$ at $V_{\out}$.
Furthermore, the exponent of $q$ determines the valency of $\Gamma$
at $V_{\out}$. 
Thus we see that a term of
$\int_{\Xi_0} e^{q(W_k(q)-y_0)}\Omega$ of the form
\begin{equation}
\label{intexpand}
q^{\nu+2}\kappa^d \prod_i {1\over \nu_i!}\Mult(h_i)^{\nu_i} u_{I(h_i)}
\end{equation}
gives a marked tropical curve $h:(\Gamma,p_1,\ldots,p_{3d-2-\nu},x)
\rightarrow M_{\RR}$ of degree $d$ in $(\PP^2,P_1,\ldots,P_k)$
with $V_{\out}\in\Gamma^{[0]}$ the vertex of $E_x$ and
$\Val(V_{\out})=\nu+3$, $h(x)=Q$. 

Conversely, by Lemma \ref{curvechop}, (3), if we have such a tropical
curve, we can split this curve up into a collection of Maslov index
two disks with endpoint $Q$. This in turn gives a term of 
$\int_{\Xi_0} e^{q(W_k(q)-y_0)}\Omega$ of the form \eqref{intexpand}.
This gives a one-to-one correspondence between such terms and curves.

Now given such a tropical curve $h$, let $n_0(x)$,
$n_1(x)$ and $n_2(x)$ be the number of unbounded edges of $h$
with vertex $V_{\out}$
in the directions $m_0$, $m_1$ and $m_2$ respectively.
Then the term corresponding to such an $h$ 
is
\[
{1\over n_0(x)! n_1(x)! n_2(x)!}q^{\nu+2} \kappa^d u_{I(h)}
\prod_{i}\Mult(h_i)
=
\Mult_x^0(h)q^{\nu+2} \kappa^d u_{I(h)}
\prod_{V\in\Gamma\atop V\not\in E_x}\Mult_V(h).
\]
Thus we see that for an index set $I=\{i_1,\ldots,i_{3d-2-\nu}\}
\subseteq \{1,\ldots,k\}$, $i_1<\cdots<i_{3d-2-\nu}$,
the coefficient of $q^{\nu+2}\kappa^d u_{I}$ in
$\int_{\Xi_0} e^{q(W_k(q)-y_0)}\Omega$ is
\[
\langle P_{i_1},\ldots,P_{i_{3d-2-\nu}},\psi^{\nu} Q\rangle_d^{\trop}.
\]
By Lemma \ref{integralindependentlemma}, 
the integral is independent of the position
of the $P_i$'s and $Q$, and hence this number is independent of $I$,
as can be seen simply by permuting the $P_i$'s. Thus this
number can be interpreted as 
\[
\langle T_2^{3d-2-\nu},\psi^{\nu} T_2\rangle_d^{\trop},
\]
now shown to be completely independent of the choice of
$P_1,\ldots,P_k$ and $Q$. In addition,  we see that
\begin{eqnarray*}
\int_{\Xi_0} e^{qW_k(Q)}\Omega&=&
e^{qy_0}\bigg(1+\sum_{d\ge 1}\sum_{\nu\ge 0}
\langle T_2^{3d-2-\nu},\psi^{\nu} T_2\rangle^{\trop}_d\kappa^dq^{\nu+2}
\sum_{I\subseteq \{1,\ldots,k\}
\atop \# I=3d-2-\nu} u_I\bigg)\\
&=&e^{qy_0}\bigg(1+\sum_{d\ge 1}\sum_{\nu\ge 0}
\langle T_2^{3d-2-\nu},\psi^{\nu} T_2\rangle_d^{\trop}
q^{\nu+2} e^{dy_1} {y_2^{3d-2-\nu}
\over (3d-2-\nu)!}\bigg),
\end{eqnarray*}
where we take $y_2=\sum_{i=1}^k u_i$ and formally take
$y_1=\log\kappa$. The latter is the expression for $J_0^{\trop}$,
hence the result.
\qed

\begin{example}
$\langle \psi^{3d-2}T_2\rangle_d^{\trop}$ is easily computed: there is only one
tropical curve of degree $d$ with a vertex of valency $3d$ at $Q$,
namely the curve which has $d$ legs of weight one in each of the three
directions $(-1,-1)$, $(1,0)$, and $(0,1)$, and hence contributes a multiplicity
of $1/(d!)^3$, so
\[
\langle\psi^{3d-2}T_2\rangle_d^{\trop}={1\over (d!)^3}.
\]

Using the formula for $W_1(Q)=y_0+x_0+x_1+x_2+u_1x_1x_2$ given in 
Example~\ref{firstP2example}, we see that contributions to
$\langle T_2,\psi^{3d-3}T_2\rangle_d^{\trop}$ comes from
terms of the form $q^{3d-1}x_0^dx_1^{d-1}x_2^{d-1}(u_1x_1x_2)$,
which shows that 
\[
\langle T_2,\psi^{3d-3}T_2\rangle_d^{\trop}={1\over (d!)(d-1)!(d-1)!}.
\]

Both these give the correct non-tropical descendent invariants.
\end{example}

\section{Scattering diagrams}

It will be important for our task of computing the period
integrals $\int_{\Xi} e^{qW_k(Q)}\Omega$ to understand
how the functions $W_k(Q)$ depend on the choice of the
base point $Q$ as well as the $P_i$'s. It turns out that
there is a chamber structure to $M_{\RR}$. For $Q$ varying within
a chamber, $W_k(Q)$ is constant, and there are wall-crossing
formulas for when $Q$ moves between chambers. 
In fact, these wall-crossing formulas have
already appeared in the context of \emph{scattering diagrams},
which appeared in \cite{GS}, following 
Kontsevich and Soibelman \cite{ks}, to construct toric degenerations
of Calabi-Yau manifolds. In addition, scattering diagrams 
have also been seen to play a role in enumerative and tropical
geometry in \cite{GPS} as well as in wall-crossing
formulas for Donaldson-Thomas invariants in \cite{KSDT}. 
These types of wall-crossing formulas also appeared in \cite{Auroux},
arising there, as here, from Maslov index zero disks.
Our discussion here will make the enumerative relevance of
these ideas clearer.

We shall repeat the definition of scattering diagram here, with slightly
different definitions and conventions than was used in \cite{GS} or
\cite{GPS}. As in \S 2, we fix $M$ a rank two lattice, and a complete fan
$\Sigma$ in $M_{\RR}$ defining a non-singular toric
surface, giving an exact sequence
\[
0\mapright{} K_{\Sigma}\mapright{} T_{\Sigma}\mapright{r} M
\mapright{}0.
\]

\begin{definition} Fix $k\ge 0$.
\begin{enumerate}
\item A \emph{ray} or \emph{line} is a pair $(\fod,f_{\fod})$ such that
\begin{itemize}
\item $\fod\subseteq M_{\RR}$
is given by
\[
\fod=m_0'-\RR_{\ge 0}r(m_0)
\]
if $\fod$ is a ray and
\[
\fod=m_0'-\RR r(m_0)
\]
if $\fod$ is a line,
for some $m_0'\in M_{\RR}$ and $m_0\in T_{\Sigma}$
with $r(m_0)\not=0$. The set
$\fod$ is called the \emph{support} of the line or ray. If $\fod$ is
a ray, $m_0'$ is called the \emph{initial point} of the ray, written
as $\Init(\fod)$.
\item $f_{\fod}\in \CC[z^{m_0}]\otimes_{\CC} R_k\subseteq 
\CC[T_{\Sigma}]\otimes_{\CC} R_k\lfor y_0\rfor.$
\item $f_{\fod}\equiv 1 \mod (u_1,\ldots,u_k)z^{m_0}$.
\end{itemize}
\item A \emph{scattering diagram} $\foD$ 
is a finite collection of lines and rays.
\end{enumerate}
\end{definition}

If $\foD$ is a scattering diagram, we write
\[
\Supp(\foD):=\bigcup_{\fod\in\foD} \fod\subseteq M_{\RR}
\]
and
\[
\Sing(\foD):=\bigcup_{\fod\in\foD} \partial\fod \cup \bigcup_{\scriptstyle
\fod_1,\fod_2
\atop\scriptstyle
\dim \fod_1\cap\fod_2=0} \fod_1\cap \fod_2.
\]
Here $\partial\fod=\{\Init(\fod)\}$ if $\fod$ is a ray, and is empty
if $\fod$ is a line.

\begin{construction} Given a smooth immersion $\gamma:[0,1]\rightarrow
M_{\RR}\setminus\Sing(\foD)$
with endpoints not in any element of a scattering diagram $\foD$,
such that $\gamma$ intersects elements of $\foD$ transversally,
we can define a ring automorphism $\theta_{\gamma,\foD}
\in \Aut(\CC[T_{\Sigma}]\otimes_{\CC} R_k\lfor y_0\rfor)$, 
the \emph{$\gamma$-ordered product of $\foD$}. Explicitly, we can find
numbers
\[
0<t_1\le t_2\le \cdots \le t_s<1
\]
and elements $\fod_i\in\foD$ such that $\gamma(t_i)\in\fod_i$ and
$\fod_i\not=\fod_j$ if $t_i=t_j$, $i\not=j$, and $s$ taken as large as
possible. Then for each $i$, define $\theta_{\gamma,\fod_i}\in
\Aut(\CC[T_{\Sigma}]\otimes_{\CC} R_k\lfor y_0\rfor)$ to be
\begin{eqnarray*}
\theta_{\gamma,\fod_i}(z^m)&=&z^mf_{\fod_i}^{\langle n_0,r(m)\rangle}\\
\theta_{\gamma,\fod_i}(a)&=&a
\end{eqnarray*}
for $m\in T_{\Sigma}$, $a\in R_k\lfor y_0\rfor$,
where $n_0\in N$ is chosen to be primitive, annihilate the tangent
space to $\fod_i$, and is finally completely determined by the sign
convention that
\[
\langle n_0,\gamma'(t_i)\rangle <0.
\]
We then define
\[
\theta_{\gamma,\foD}=\theta_{\gamma,\fod_s}\circ\cdots\circ
\theta_{\gamma,\fod_1}.
\]

There is still some ambiguity to the ordering if $\gamma$
crosses several overlapping rays. However, an easy check
shows that automorphisms associated to parallel rays commute, so
the order is irrelevant.
Automorphisms associated with non-parallel rays do not necessarily
commute, hence the need for $\gamma$ to avoid points of $\Sing(\foD)$.

We will also allow the possibility that $\gamma$ is piecewise linear
so that $\gamma'$ may not be defined at $t_i$, but still have
$\gamma$ pass from one side of $\fod$ to the other, in which
case we take $n_0$ so that $\gamma$ passes from the side of $\fod$
where $n_0$ is larger to the side it is smaller.

It is easy to check that $\theta_{\gamma,\foD}$ only depends on
the homotopy class of the path $\gamma$ inside $M_{\RR}\setminus
\Sing(\foD)$.
\end{construction}

\bigskip

\begin{example}
\label{basicscatterexample}
Let $\foD=\{(\fod_1,f_{\fod_1}),(\fod_2,f_{\fod_2}),(\fod_3,f_{\fod_3})\}$
with
\begin{eqnarray*}
\fod_1=\RR r(m_1),&\quad& f_{\fod_1}=1+c_1w_1z^{m_1},\\
\fod_2=\RR r(m_2),&\quad& f_{\fod_2}=1+c_2w_2z^{m_2},\\
\fod_3=-\RR_{\ge 0}(r(m_1+m_2)),&\quad&f_{\fod_3}
=1+c_1c_2w_{\out}|r(m_1)\wedge r(m_2)|z^{m_1+m_2},
\end{eqnarray*}
where $m_1,m_2\in T_{\Sigma}\setminus K_{\Sigma}$, and $w_1,w_2$
and $w_{\out}$ are the indices\footnote{If $m\in M\setminus \{0\}$,
the \emph{index} of $m$ is the largest positive integer
$w$ such that $m=wm'$ with $m'\in M$.} of $r(m_1)$, $r(m_2)$ and $r(m_1+m_2)$
respectively.
The expression $|r(m_1)\wedge r(m_2)|$ denotes the absolute
value of $r(m_1)\wedge r(m_2)\in \bigwedge^2 M\cong\ZZ$.
Suppose $c_1,c_2\in R_k$ satisfy $c_1^2=c_2^2=0$.
Then one can calculate that if $\gamma$ is a loop around the origin,
then $\theta_{\gamma,\foD}$ is the identity. See Figure \ref{sampleloop},
where $\theta_{\gamma,\foD}=\theta_2^{-1}\theta_3\theta_1^{-1}\theta_2\theta_1$,
with $\theta_i$ the automorphism coming from the first crossing of $\fod_i$.
\end{example}

\begin{figure}
\input{sampleloop.pstex_t}
\caption{}
\label{sampleloop}
\end{figure}

We now relate this concept to the behaviour of the Landau-Ginzburg
potential under change of the base-point $Q$. We fix general
$P_1,\ldots,P_k\in M_{\RR}$, and will study the behaviour
of $W_k(Q)$ as a function of the
base-point $Q$. The discussion in this section
will be for arbitrary complete fans $\Sigma$ in $M_{\RR}$
defining a non-singular toric surface.

First, we define a variant of tropical disk.

\begin{definition} A \emph{tropical tree in $(X_{\Sigma},P_1,
\ldots,P_k)$} is
a $d$-pointed tropical curve $h:(\Gamma,p_1,\ldots,p_d)
\rightarrow M_{\RR}$ with $h(p_j)=P_{i_j}$ 
$1\le i_1<\cdots<i_d\le k$, along with the additional
data of a choice of unmarked unbounded edge $E_{\out}\in \Gamma^{[1]}_{\infty}$
such that for any $E\in \Gamma^{[1]}_{\infty}\setminus \{E_{\out}\}$,
$h(E)$ is a point or is a translate of some $\rho\in\Sigma^{[1]}$.
The degree of $h$, $\Delta(h)$, is defined \emph{without} counting the edge
$E_{\out}$, which need not be a translate of any $\rho\in\Sigma^{[1]}$.

The \emph{Maslov index} of $h$ is
\[
MI(h):=2(|\Delta(h)|-d).
\]
\end{definition}

Given $h$ and a point $V_{\out}$ in the interior of $E_{\out}$, 
we can remove the
unbounded component of $E_{\out}\setminus \{V_{\out}\}$ from $\Gamma$
to obtain $\Gamma'$. Note $V_{\out}$ is a univalent vertex
of $\Gamma'$. Take $h':\Gamma'\rightarrow M_{\RR}$ with $h'=h|_{\Gamma'}$.
Then $h'$ is a tropical disk with boundary $h(V_{\out})$ and Maslov
index $MI(h')=MI(h)$, since $|\Delta(h')|=|\Delta(h)|$.

As in Lemma \ref{diskmoduli}, 
standard tropical dimension counting arguments show that,
for general choice of $P_1,\ldots,P_k$, a tropical tree $h$
moves in a family of dimension $MI(h)/2$. In particular, the
set of Maslov index zero trees is a finite set, which we denote
by $\Trees(\Sigma,P_1,\ldots,P_k)$. As usual, with general
choice of $P_1,\ldots,P_k$, we can assume all these trees are
trivalent.

\begin{definition} We define $\foD(\Sigma,P_1,\ldots,P_k)$ to be the
scattering diagram which contains one ray
for each element $h$ of $\Trees(\Sigma,P_1,\ldots,P_k)$. The ray
corresponding to $h$ is of the form $(\fod,f_{\fod})$, where 
\begin{itemize}
\item $\fod=h(E_{\out})$.
\item $f_{\fod}=1+w_{\Gamma}(E_{\out})\Mult(h)z^{\Delta(h)}u_{I(h)}$,
where $u_{I(h)}=\prod_{i\in I(h)} u_i$
and $I(h)\subseteq \{1,\ldots,k\}$ is defined by
\[
I(h):=\{i\,|\, \hbox{$h(p_j)=P_i$ for some $j$}\}.
\] 
\end{itemize}
\end{definition}

\begin{figure}
\input{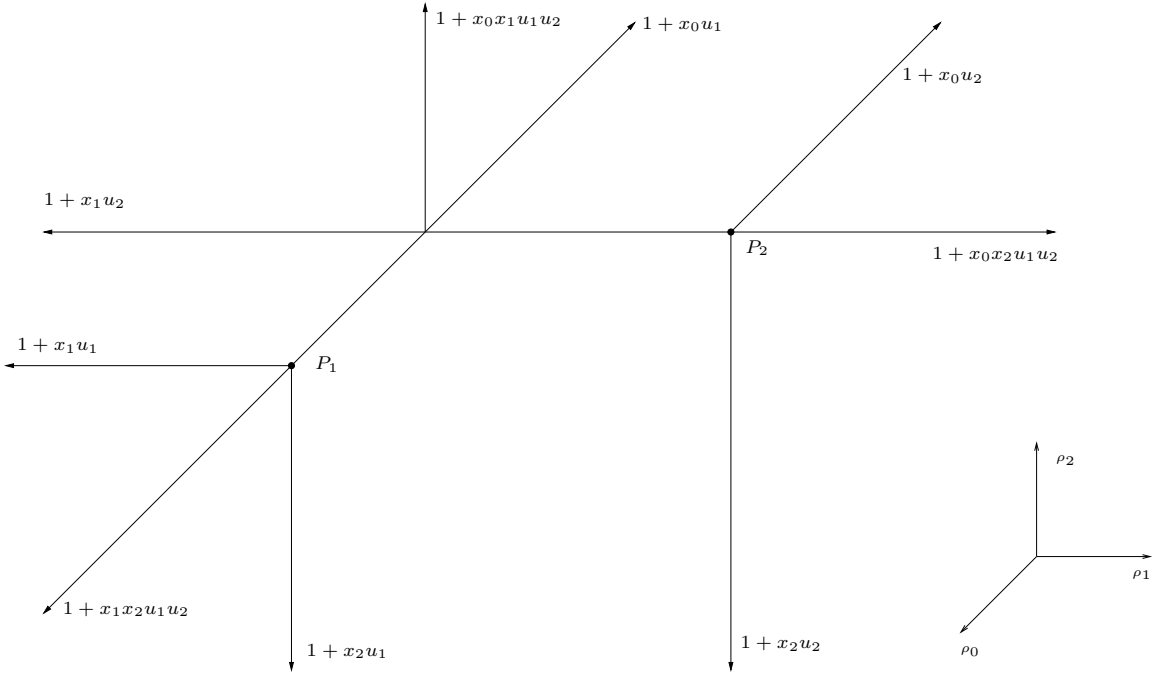}
\caption{The scattering diagram $\foD(\Sigma,P_1,P_2)$. For reference,
the fan for $\PP^2$ is given in the lower right-hand corner.}
\label{twopointscatter}
\end{figure}

\begin{example}
$\foD(\Sigma,P_1,P_2)$ is illustrated in Figure \ref{twopointscatter},
where $\Sigma$ is the fan for $\PP^2$, as given in Example \ref{firstP2example}.
\end{example}

\begin{proposition}
\label{gammaid}
If $P\in\Sing(\foD(\Sigma,P_1,\ldots,P_k))$ is a singular point
with $P\not\in \{P_1,\ldots,P_k\}$, and $\gamma_P$ is a small loop
around $P$, then $\theta_{\gamma_P,\foD(\Sigma,P_1,\ldots,P_k)}=
\id$.
\end{proposition}

\proof Let $P$ be such a singular point. Suppose that 
$\fod\in\foD(\Sigma,P_1,\ldots,P_k)$ has $\Init(\fod)=P$,
and let $h$ be the corresponding tree. Then
the unique vertex $V$ of $\Gamma$ on $E_{\out}$
has $h(V)=P$. By the generality assumption, $V$ is trivalent, so
if we cut $\Gamma$ at $V$, we obtain two tropical disks $h'_1:\Gamma_1'
\rightarrow M_{\RR}$ and $h'_2:\Gamma_2\rightarrow M_{\RR}$ with boundary
$P$ and $V_{\out}=V$ in both cases. Now $MI(h)=MI(h'_1)+MI(h'_2)$, 
so $MI(h_1')=MI(h_2')=0$ is the only possibility. Thus $h'_1$, $h'_2$
extend to tropical trees $h_i:\Gamma_i\rightarrow M_{\RR}$, with
corresponding rays $\fod_1,\fod_2$. Note that $P\not=\Init(\fod_1),
\Init(\fod_2)$ and $I(h_1)\cap I(h_2)=\emptyset$. 
So every ray $\fod$ with $P=\Init(\fod)$ arises from
the collision of two rays $\fod_1$, $\fod_2$ with $P\not=\Init(\fod_i)$.

Conversely, if we are given two such rays $\fod_1$, $\fod_2$ passing
through $P$ corresponding
to trees $h_1$ and $h_2$ with $I(h_1)\cap I(h_2)=\emptyset$, 
we obtain a new tree by cutting $h_1$
and $h_2$ at $P$ to get Maslov index zero disks 
$h_i':\Gamma_i'\rightarrow M_{\RR}$ with boundary at
$P$. Next glue $\Gamma_1'$ and $\Gamma_2'$ at the outgoing
vertex $V$,
and add an additional unbounded edge $E_{\out}$ with endpoint
$V$ to get a graph $\Gamma$. If 
$E_{\out,1},E_{\out,2}$ are the two outgoing edges of $\Gamma_1'$
and $\Gamma_2'$ respectively, with primitive tangent vector
to $h_i'(E_{\out,i})$ pointing away from $P$ being $m_i'$, then
we define $h:\Gamma\rightarrow M_{\RR}$ to restrict to $h_i'$ on $\Gamma_i'$
and to take $E_{\out}$ to the ray $P-\RR_{\ge 0}(w_{\Gamma}(E_{\out,1}) m_1'
+w_{\Gamma}(E_{\out,2})m_2')$. By taking $w_{\Gamma}(E_{\out})$ to
be the index of $w_{\Gamma}(E_{\out,1}) m_1'
+w_{\Gamma}(E_{\out,2})m_2'$, we find $h$ is balanced at $V$. Thus
$h$ is a tropical tree, whose Maslov index is zero.

To prove the proposition, define a new scattering diagram $\foD_P$,
whose elements are in one-to-one correspondence with elements of
$\foD(\Sigma,P_1,\ldots,P_k)$ 
containing $P$. If $(\fod,f_{\fod})\in\foD(\Sigma,P_1,\ldots,P_k)$ is a ray 
containing $P$
then the corresponding element of $\foD_P$ will be $(\fod',f_{\fod})$,
where $\fod'$ is the tangent line (through the origin) of $\fod$
if $P\not=\Init(\fod)$ and is the ray $\fod-P$ with endpoint the origin 
otherwise. If $\gamma_0$ is a loop around the origin with the same orientation
as $\gamma_P$, then 
$\theta_{\gamma_0,\foD_P}=\theta_{\gamma_P,\foD(\Sigma,P_1,\ldots,P_k)}$.

First consider the simplest case, when $\foD_P$ contains two lines and
at most one ray. If the two lines correspond to trees $h_1$ and $h_2$,
and $I(h_1)\cap I(h_2)\not=\emptyset$, then $h_1$ and
$h_2$ \emph{cannot} be glued as above since they pass through some common
marked point $P_i$. Thus $\foD_P$ contains no rays.
In this case, the automorphisms associated to
$\fod_1$ and $\fod_2$ commute by Example \ref{basicscatterexample}
since $u_{I(h_1)}u_{I(h_2)}=0$, and
so $\theta_{\gamma_0,\foD_P}$ is the identity.

If, on the other hand, $I(h_1)\cap I(h_2)=\emptyset$, then $h_1$ and $h_2$
can be glued to obtain a new tree $h$, and $\foD_P$ consists of 
three elements $\fod_1$, $\fod_2$ and $\fod$, corresponding to $h_1, h_2$
and $h$ respectively. Now
\[
f_{\fod_i}=1+w_{\Gamma}(E_{\out,i})\Mult(h_i) z^{\Delta(h_i)}u_{I(h_i)}
\]
for $i=1,2$ and
\begin{align*}
f_{\fod}\,=\,1+w_{\Gamma}(E_{\out})&\Mult(h) z^{\Delta(h)}u_{I(h)}\\ 
=\,1+w_{\Gamma}(E_{\out})&\Mult(h_1)\Mult(h_2)\Mult_V(h)
z^{\Delta(h_1)+\Delta(h_2)}u_{I(h_1)}u_{I(h_2)}\\
=\,1+w_{\Gamma}(E_{\out})\Mult(h_1)\Mult&(h_2)
w_{\Gamma_1}(E_{\out,1})w_{\Gamma_2}(E_{\out,2})|m_1'\wedge m_2'|
z^{\Delta(h_1)+\Delta(h_2)}u_{I(h_1)}u_{I(h_2)}\\
=\,1+w_{\Gamma}(E_{\out})\Mult(h_1)\Mult&(h_2)
|r(\Delta(h_1))\wedge r(\Delta(h_2))|
z^{\Delta(h_1)+\Delta(h_2)}u_{I(h_1)}u_{I(h_2)}.
\end{align*}
Thus from Example \ref{basicscatterexample}, $\theta_{\gamma_0,\foD_P}$
is the identity.

For the general case, we have some finite set of lines in $\foD_P$,
along with some rays. Suppose that there are three lines in $\foD_P$
corresponding
to trees $h_1$, $h_2$ and $h_3$ with $I(h_1),I(h_2)$ and $I(h_3)$
mutually disjoint. Then as in the case of two lines above, these trees
can be glued at $P$, obtaining a Maslov index zero tree with a quadrivalent
vertex. However, since $P_1,\ldots,P_k$ are in general position, no
Maslov index zero tree has a vertex with valence $>3$. Thus this case
does not occur. On the other hand,
given two lines corresponding to trees $h_1$, $h_2$ with $I(h_1)\cap I(h_2)
=\emptyset$, these two trees can be glued as above at $P$ to obtain
a new Maslov index zero tree. Thus the
rays in $\foD_P$ are in one-to-one correspondence with
pairs of lines $\fod_1,\fod_2\in\foD_P$ corresponding to trees $h_1$ and
$h_2$ with $I(h_1)\cap I(h_2)=\emptyset$. 
So we can write 
\[
\foD_P=\{\fod_1,\ldots,\fod_n\}\cup \bigcup_{j=1}^m \foD_i
\]
where $\fod_1,\ldots,\fod_n$ are lines corresponding to trees $h$
such that $I(h)\cap I(h')\not=\emptyset$ for any Maslov index zero
tree $h'$ with outgoing edge
passing through $P$, and $\foD_1,\ldots,\foD_m$ are
scattering diagrams each consisting of two lines and one ray,
with the lines corresponding to trees $h_1$ and $h_2$ with 
$I(h_1)\cap I(h_2)=\emptyset$ and the ray corresponding to the
tree obtained by gluing $h_1$ and $h_2$ at $P$.

Now computing $\theta_{\gamma_0,\foD_P}$ is an exercise in commutators.
Note that if $\fod_1,\fod_2\in\foD_P$ correspond to two trees $h_1,h_2$
with $I(h_1)\cap I(h_2)\not=\emptyset$, then as already observed,
$\theta_{\gamma_0,\fod_1}$ and $\theta_{\gamma_0,\fod_2}$ commute. Thus
after using this commutation, one can write
\[
\theta_{\gamma_0,\foD_P}=
\bigg(\prod_{i=1}^n \theta_{\gamma_0,\fod_i}\circ\theta_{\gamma_0,\fod_i}^{-1}
\bigg)\circ \prod_{j=1}^m \theta_{\gamma_0,\foD_j}.
\]
Of course $\theta_{\gamma_0,\fod_i}\circ\theta_{\gamma_0,\fod_i}^{-1}=\id$
and $\theta_{\gamma_0,\foD_j}=\id$ by the special case already carried out.
Thus $\theta_{\gamma_0,\foD_P}=\id$ in this general case.
\qed

\medskip

\begin{remark}
\label{Piremark}
Note that the rays in $\foD=\foD(\Sigma,P_1,\ldots,P_k)$ with endpoint
$P_i$ are in one-to-one correspondence with Maslov index two
disks in $(X_{\Sigma},P_1,\ldots,P_{i-1},P_{i+1},\ldots,P_k)$
with boundary $P_i$. Indeed, taking any such Maslov index two
disk, extending the
outgoing edge to get a tropical tree, we can mark the point on this
outgoing edge which maps to $P_i$, thus getting a tropical tree
with Maslov index zero in $(X_{\Sigma},P_1,\ldots,P_k)$. The corresponding
ray in $\foD$ has endpoint $P_i$.
Conversely, given a ray in $\foD$ with endpoint $P_i$, this corresponds
to a Maslov index zero tree such that the vertex $V$ adjacent to $V_{\out}$
is the vertex of a marked edge $E_x$ mapping to $P_i$. By cutting 
this tree at $V$, removing the marked edge mapping to $P_i$, we get
a Maslov index two disk with boundary $P_i$.

Furthermore, by the general position of the $P_j$, there are 
no rays in $\foD$ containing $P_i$ but
which don't have $P_i$ as an endpoint.
\end{remark}

\medskip

One benefit of this scattering diagram approach is that it is
easy to describe the Maslov index two disks with boundary a
general point $Q$, using what we call \emph{broken lines}:

\begin{definition}
\label{brokenlinedef}
A broken line is a continuous proper piecewise linear map $\beta:(-\infty,0]
\rightarrow M_{\RR}$ with endpoint $Q=\beta(0)$, along with some
additional data described as follows. Let 
\[
-\infty=t_0<t_1<\cdots< t_n=0
\]
be the smallest set of real numbers such that $\beta|_{(t_{i-1},t_i)}$
is linear. Then for each $1\le i\le n$, we are given the additional
data of a monomial $c_iz^{m^{\beta}_i}\in
\CC[T_{\Sigma}]\otimes_{\CC} R_k\lfor y_0\rfor$ with $m^{\beta}_i\in T_{\Sigma}$ and
$0\not=
c_i\in R_k\lfor y_0\rfor$. Furthermore, this data should satisfy the following
properties:
\begin{enumerate}
\item
For each $i$, $r(m^{\beta}_i)$ points in the same direction as $-\beta'(t)$
for $t\in (t_{i-1},t_i)$.
\item $m_1^{\beta}=t_{\rho}$ for some $\rho\in\Sigma^{[1]}$ and $c_1=1$.
\item $\beta(t_i)\in\Supp(\foD(\Sigma,P_1,\ldots,P_k))$ for $1\le i\le n$.
\item If the image of $\beta$ is disjoint from $\Sing(\foD(\Sigma,P_1,
\ldots,P_k))$, and $\beta(t_i)\in \fod_1\cap\cdots\cap\fod_s$ (necessarily
this intersection is one-dimensional), then $c_{i+1}z^{m^{\beta}_{i+1}}$
is a term in 
\[
(\theta_{\beta,\fod_1}\circ\cdots\circ\theta_{\beta,\fod_s})(c_iz^{m^{\beta}_i}).
\]
By this, we mean the following. Suppose $f_{\fod_j}=1+c_{\fod_j}z^{m_{\fod_j}}$,
$1\le j\le s$, with $c^2_{\fod_j}=0$, 
and $n\in N$ is primitive, orthogonal to all the $\fod_j$'s,
chosen so that
\begin{eqnarray*}
(\theta_{\beta,\fod_1}\circ\cdots\circ\theta_{\beta,\fod_s})(c_iz^{m^{\beta}_i})
&=&c_iz^{m^{\beta}_i}\prod_{j=1}^s(1+c_{\fod_j}z^{m_{\fod_j}})^{\langle n,
r(m^{\beta}_i)\rangle}\cr
&=&c_iz^{m^{\beta}_i}\prod_{j=1}^s(1+\langle n,r(m^{\beta}_i)\rangle
c_{\fod_j}z^{m_{\fod_j}}).
\end{eqnarray*}
Then we must have
\[
c_{i+1}z^{m^{\beta}_{i+1}}=c_iz^{m^{\beta}_i}\prod_{j\in J}
\big(\langle n,r(m^{\beta}_j)\rangle c_{\fod_j}z^{m_{\fod_j}}\big)
\]
for some index set $J\subseteq \{1,\ldots,s\}$.
We think of this as $\beta$ being bent by time $t_i$ by the rays
$\{\fod_j\,|\,j\in  J\}$.

\item If the image of $\beta$ is not disjoint from $\Sing(\foD(\Sigma,P_1,
\ldots,P_k))$, then $\beta$ is the limit of a family of broken lines
which are disjoint from $\Sing(\foD(\Sigma,P_1,\ldots,P_k))$.
More precisely, there is: 
\begin{itemize}
\item A continuous map $B:(-\infty,0]\times
[0,1]\rightarrow M_{\RR}$. 
\item Continuous functions $t_0,
\ldots,t_n:[0,1]\rightarrow [-\infty,0]$ such that
\[
-\infty=t_0(s)\le t_1(s)\le\cdots\le t_n(s)=0
\]
for $s\in[0,1]$, with strict inequality for $s<1$.
\item Monomials $c_iz^{m^i_B}$ for $1\le i\le n$.
\end{itemize}
This data satisfies $B_s:=B|_{(-\infty,0]\times \{s\}}$
(with the data $c_iz^{m^i_B}$) 
is a broken line not passing through a point of $\Sing(\foD(\Sigma,
P_1,\ldots,P_k))$ for $s<1$, and $\beta:=B_1$.

Note that in taking such a limit, we might have $t_{i-1}$ and $t_i$
coming together for various $i$, so the limit might have fewer
linear segments.
\end{enumerate}
\end{definition}

\begin{figure}
\input{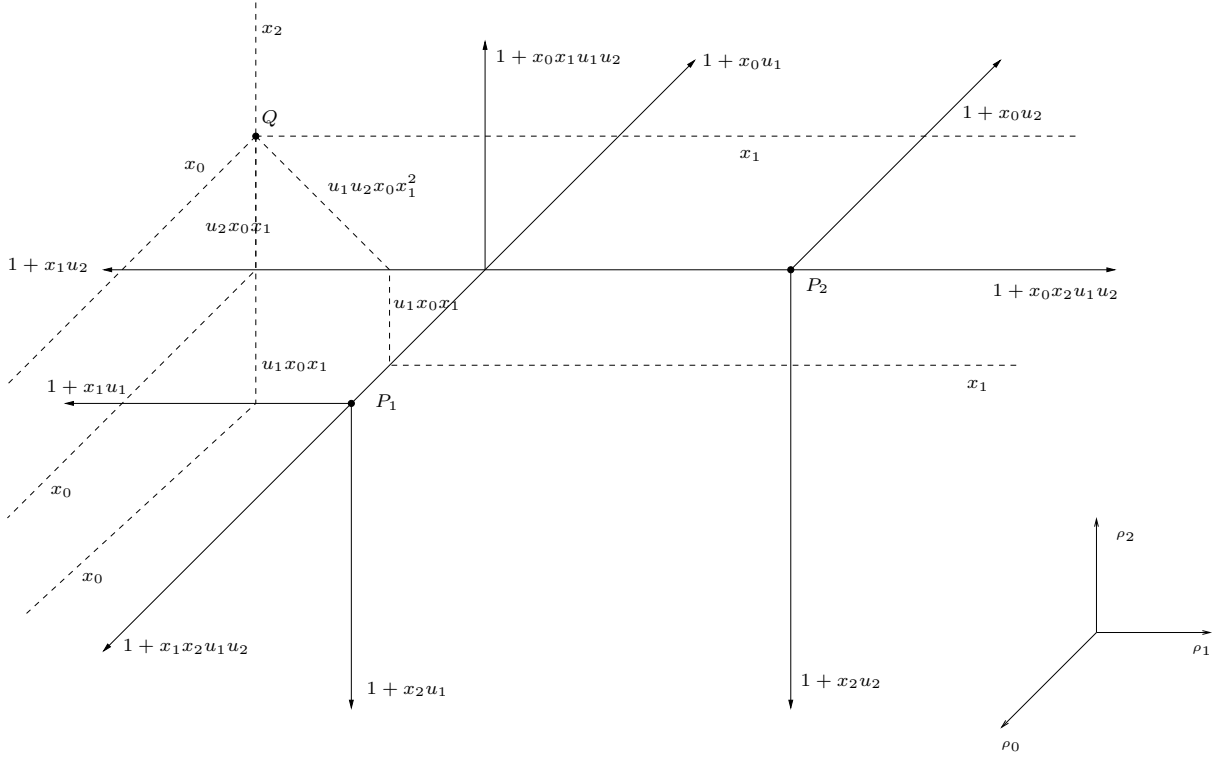}
\caption{The broken lines ending at $Q$.}
\label{twopointscatter2}
\end{figure}

\begin{example}
Again, in the case of $\PP^2$, $k=2$, Figure \ref{twopointscatter2}
shows the broken lines with $\beta(0)$ the given point $Q$. The
segments of each broken line are labelled with their corresponding
monomial.
\end{example}

\begin{proposition}
\label{brokenlinemi2}
If $Q\not\in\Supp(\foD(\Sigma,P_1,\ldots,P_k))$ is general, then 
there is a one-to-one correspondence between broken lines
with endpoint $Q$ and Maslov index two disks with boundary
$Q$. In addition, if $\beta$ is a broken line corresponding
to a disk $h$, and $cz^m$ is the monomial associated to the
last segment of $\beta$ (the one whose endpoint is $t_n=0$), then
\[
cz^m=\Mono(h).
\]
\end{proposition}

\proof 
We first prove the following claim:

\medskip

\emph{Claim:} Let $h:\Gamma'\rightarrow M_{\RR}$ be a Maslov index
two disk in $(X_{\Sigma},P_1,\ldots,P_k)$ 
with boundary $Q'\in M_{\RR}$. Suppose furthermore that all vertices of
$\Gamma'$ except $V_{\out}$ are trivalent and $h$ cannot be deformed
continuously in a family of Maslov index two disks with boundary $Q'$.
Then there is a uniquely determined
subset $\Xi=\Xi(h)\subseteq\overline{\Gamma}$ which is a union of edges
of $\overline{\Gamma}$ and is homeomorphic to $[-\infty,0]$,
connecting some point in $\overline{\Gamma}^{[0]}_{\infty}\setminus\{V_{\out}\}$
to $V_{\out}$, satisfying:
\begin{enumerate}
\item $\Xi$ is disjoint from $\partial E_{p_i}$ for all $i$. 
\item The restriction of $h$ to the closure of any connected component
of $\Gamma'\setminus\Xi$ is a Maslov index zero disk.
\end{enumerate}

\proof
We proceed inductively on the number of vertices of $\Gamma'$. If
$\Gamma'$ has only one vertex, $V_{\out}$, then $\Gamma'$ has only
one edge and no marked edges. We simply take $\Xi(h)$ to consist of this
edge.

For the induction step, let $\Gamma'$ have
outgoing edge $E_{\out}$ with vertices $V_{\out}$ and $V$,
and $h(V_{\out})=Q'$. 

First, we will show that $V$ cannot be a vertex of
some marked edge $E_{p_i}$. If it were, so that
$h(V)=h(p_i)=P_j$ for some $j$, then we can cut $h$ at $V$ and remove
the marked edge $E_{p_i}$. This gives a disk $h'$ with boundary $P_j$ 
but with one less marked point than $h$. Hence $h'$ is Maslov index four, and
thus by Lemma \ref{diskmoduli}, $h'$ can be deformed in a one-parameter
family while keeping the endpoint $P_j$ fixed. Note that for small
deformations of $h'$, the edge $h'(E'_{\out})$ does not change its tangent
direction. Thus a deformation of $h'$ can be extended to a deformation
of $h$. This contradicts the assumption that $h$ cannot be deformed.
Thus $V$ cannot be a vertex of some $E_{p_i}$. 

Now split $h$ at $V$, letting $\Gamma_1'$ and $\Gamma_2'$ be
the closures of the two connected components of $\Gamma'\setminus\{V\}$
not containing $V_{\out}$. Let $h_i=h|_{\Gamma_i'}$.
This gives two disks $h_1$, $h_2$ with boundary
$h(V)$. We have $MI(h)=MI(h_1)+MI(h_2)$. Suppose $MI(h_1)\ge 4$. Then $h_1$
can be deformed leaving the endpoint $h(V)$ fixed, and by gluing such a 
deformation to $h_2$, we obtain a deformation of $h$, again a contradiction.
Thus $MI(h_1),MI(h_2)\le 2$, so we must have $MI(h_1)=0$ and $MI(h_2)=2$
or vice versa. 

Without loss of generality, assume $MI(h_2)=2$. Note
that $h_2$ is now a Maslov index two disk with boundary $Q''=h_2(V)$.
If $h_2$ could be deformed in a family of disks with boundary $Q''$,
then by gluing these deformations to $h_1$, we obtain a deformation of $h$,
a contradiction.
Thus $h_2$ satisfies the hypotheses of the Claim, but $\Gamma_2'$ has
fewer vertices than $\Gamma'$. So $\Xi(h_2)$ exists by the induction
hypothesis, and we can take $\Xi(h)$ to be
\[
\Xi(h)=E_{\out}\cup \Xi(h_2).
\]
$\Xi(h)$ satisfies the two desired properties because
$\Xi(h_2)$ does, $E_{\out}$ is disjoint from $\partial E_{p_i}$ for
all $i$, and $h_1$ is a Maslov index zero disk.
\qed
\medskip

Now fix a Maslov index two disk $h:\Gamma'\rightarrow M_{\RR}$ with
boundary $Q$. By the generality of $Q,P_1,\ldots,P_k$, $h$ satisfies
the hypotheses of the Claim.
Taking $\beta=h|_{\Xi(h)}$,
we see that $\beta$ is piecewise linear. Let $-\infty=t_0<\cdots<t_n=0$
be chosen as in the definition of broken line. 
Each $t_i$
corresponds to a vertex $V_i$ of $\Gamma'$. Of course 
$\Gamma'\setminus \{V_i\}$ for $i\not=n$
has two connected components not containing $V_{\out}$, and the proof of
the claim shows that restricting $h$ to the closure of one of these two
connected components yields a Maslov index two disk with boundary
$h(V_i)$ which we now call $h_i$.
The other component similarly yields a Maslov index zero disk. Hence
$\beta(t_i)
\in\Supp(\foD(\Sigma,P_1,\ldots,P_k))$ for $1\le i\le n$. 
We take the monomial $c_iz^{m_i}$ to be $\Mono(h_i)$, and need to check
that with this data $\beta$ is now a broken line.

We have just shown condition (3), and condition (2) is obvious.
Condition (1) is immediate from the balancing condition for $h_i$.
So for $Q$ sufficiently general, we only need to verify condition 
(4). We need to show that the monomial $\Mono(h_{i+1})$ attached to the
edge joining $V_i$ and $V_{i+1}$ arises from the monomial $\Mono(h_i)$
attached to the edge joining $V_{i-1}$ and $V_i$, as in Condition (4).
Suppose that the two subtrees in $\Gamma'$ rooted at $V_i$ are $g$ and 
$h_i$,
with $MI(g)=0$. Now
\[
\Mono(h_i)=\Mult(h_i)z^{\Delta(h_i)}u_{I(h_i)},
\]
and if $\fod$ is the ray corresponding to the tropical
tree obtained from $g$, then 
\[
f_{\fod}=1+w_{\out}(g)\Mult(g)z^{\Delta(g)}u_{I(g)}.
\]
Here $w_{\out}(g)$ denotes the weight of the outgoing edge of
$g$.
By the balancing condition, we can write 
$-w_{\out}(g)m_1'=r(\Delta(g))$, and $-w_{\out}(h_i)m_2'=r(\Delta(h_i))$
with $m_i'\in M$ primitive, with $m_1'$ tangent to the outgoing edge
of $g$ and $m_2'$ tangent to the outgoing edge of $h_i$.
Choosing a basis for $M$ and writing
$m_i'=(m'_{i1},m'_{i2})$, we see that to define $\theta_{\beta,\fod}$,
we take $n_0=\pm (m'_{12},-m'_{11})$, so
\begin{eqnarray*}
\theta_{\beta,\fod}(\Mono(h_i))&=&\Mult(h_i)z^{\Delta(h_i)}u_{I(h_i)}
\cdot\\
&&
\cdot(1+w_{\out}(g)\Mult(g)z^{\Delta(g)}u_{I(g)})^
{\langle\pm(m'_{12},-m'_{11}),-w_{\out}(h_i)(m'_{21},m'_{22})\rangle}\\
&=&\Mult(h_i)z^{\Delta(h_i)}u_{I(h_i)}\\
&&+\Mult(h_i)\Mult(g)|m_1'\wedge m_2'| w_{\out}(g)w_{\out}(h_i)
z^{\Delta(h_i)+\Delta(g)}u_{I(h_i)}u_{I(g)}.
\end{eqnarray*}
Here the second term occurs with a plus sign since the exponent
is always positive--- the convention on $n_0$ says that $n_0$
should be negative on vectors pointing in the direction we cross
$\fod$; but $m_2'$ is such a vector so $n_0$ is positive on $-m_2'$.

Now it is the second term we are interested in, and this is
\begin{eqnarray*}
\Mult(h_i)\Mult(g)\Mult_{V_i}(h_{i+1})z^{\Delta(h_{i+1})}u_{I(h_{i+1})}
&=&\Mult(h_{i+1})z^{\Delta(h_{i+1})}u_{I(h_{i+1})}\\
&=&\Mono(h_{i+1})
\end{eqnarray*}
as desired.

Conversely, given a broken line $\beta$, it is easy to construct
the corresponding Maslov index two disk, by attaching Maslov
index zero disks to the domain $(-\infty,0]$ of $\beta$ at each
bending point. In particular, if $\beta(t_i)$
lies in rays $\fod_1,\ldots,\fod_s\in\foD(\Sigma,P_1,\ldots,P_k)$,
and $\beta$ is bent at time $t_i$ by a subset $\{\fod_j\,|\,j\in J\}$
of these rays, then for each $j\in J$,
we attach the Maslov index zero disk with endpoint $\beta(t_i)$ 
corresponding to $\fod_j$ to $t_i\in (-\infty,0]$.
(Note that by general position of the $P_i$'s and $Q$, in fact we can
assume that $\#J=1$.) It is clear that this reverses the above process
of passing from a Maslov index two disk to a broken line.
\qed

\bigskip

The first main theorem of this section explores how $W_k(Q)$
depends on $Q$.

\begin{theorem}
\label{wallcrossingtheorem}
If $Q,Q'\in M_{\RR}\setminus \Supp(\foD(\Sigma,P_1,\ldots,P_k))$ are
general, 
and $\gamma$ is a path connecting $Q$ and $Q'$ for which
$\theta_{\gamma,\foD(\Sigma,P_1,\ldots,P_k)}$ is defined, then
\[
W_k(Q')=\theta_{\gamma,\foD(\Sigma,P_1,\ldots,P_k)}(W_k(Q)).
\]
\end{theorem}

\proof 
Let $\foD=\foD(\Sigma,P_1,\ldots,P_k)$. 
Let $\foU$ be the union of $\Supp(\foD)$ and the union
of images of all broken lines,
with arbitrary endpoint, which pass through points
of $\Sing(\foD)$.
Recall from Definition \ref{brokenlinedef}, (5), that such broken
lines are limits of broken lines that don't pass through singular
points of $\foD$. It is clear that $\dim\foU\le 1$ (of course equal
to $1$ provided $k\ge 1$). 

We will now define a continuous deformation of a broken line, much
as we did in Definition \ref{brokenlinedef}, (5).
This is a continuous
map $B:(-\infty,0]\times I\rightarrow M_{\RR}$ with $I\subseteq\RR$
an interval,
continuous functions $t_0,\ldots,t_n:I\rightarrow [-\infty,0]$
such that $-\infty=t_0(s)<t_1(s)<\cdots<t_n(s)=0$ for $s\in I$,
and monomials $c_iz^{m_i^{B}}$, $1\le i \le n$. This data satisfies
the condition that
$B_s:=B|_{(-\infty,0]\times\{s\}}$ is
a broken line in the usual sense for all $s\in I$, 
with the data $t_0(s)<\cdots<t_n(s)$
and monomials $c_iz^{m^B_i}$, $1\le i\le n$. 

We say $B_{s'}$ is a \emph{deformation} of $B_s$ for $s,s'\in I$.

Note a broken line $\beta$ which does not pass through a point
of $\Sing(\foD)$ can always be deformed continuously. This can
be done as follows. We translate the initial
ray $\beta((-\infty,t_1])$ of $\beta$. Inductively, this deforms
all the remaining segments of $\beta$. As long as one of the bending
points
does not reach a singular point of $\foD$, each bending point remains inside
exactly the same set of rays in $\foD$, and therefore the deformed
broken line can bend in exactly the same way as $\beta$.
Thus we run into trouble building this deformation only when
this deformation of $\beta$ converges to a broken line which passes
through a point of $\Sing(\foD)$, as then the set of rays
containing a bending point may jump.

From this it is clear that as long as the
endpoint of $\beta$ stays within one connected component of
$M_{\RR}\setminus\foU$, $\beta$ can be deformed continuously.
More precisely, if we consider a path $\gamma:[0,1]\rightarrow\fou$, for
$\fou$ a connected component of $M_{\RR}\setminus\foU$, and
$\beta$ is a broken line with endpoint $\gamma(0)$, then there is a
continuous deformation $B$ with $\beta=B_0$ and with $B_s(0)=\gamma(s)$,
$0\le s\le 1$. 

By Proposition \ref{brokenlinemi2},
the Maslov index
two disks with boundary $Q$ are in one-to-one correspondence with the 
broken lines with endpoint $Q$ for $Q$ general. Thus by the above discussion,
$W_k(Q)$ is constant for $Q$ varying inside a connected component
of $M_{\RR}\setminus\foU$.

We will now analyze carefully how broken lines change if their endpoint
passes in between different connected components of $M_{\RR}\setminus
\foU$.
So now consider two connected components $\fou_1$ and $\fou_2$
of $M_{\RR}\setminus\foU$. Let $L=\overline{\fou_1}\cap\overline{\fou_2}$, 
and assume $\dim L=1$.
Let $Q_1$ and $Q_2$ be general 
points in $\fou_1$ and $\fou_2$, near $L$, positioned on opposite
sides of $L$. Let
$\gamma:[0,1]\rightarrow M_{\RR}$ be
a short general path connecting $Q_1$ and $Q_2$ crossing 
$L$ precisely once. Let $s_0$ be the only
time at which $\gamma(s_0)\in L$. By choosing $\gamma$
sufficiently generally, we can assume that $\gamma(s_0)$ is 
a point in a neighbourhood of which $\foU$ is a manifold.

Let $\foB(Q_i)$ be the set of broken lines with
endpoint $Q_i$. Let $n_0\in N$ be a primitive
vector annihilating the tangent space to $L$ and taking a smaller
value on $Q_1$ than $Q_2$. We can decompose $\foB(Q_i)$ into
three sets $\foB^+(Q_i)$, $\foB^-(Q_i)$, and $\foB^0(Q_i)$ as follows.
For $\beta\in\foB(Q_i)$, let $m_{\beta}=\beta_*(-\partial/\partial t|_{t=0})$.
Then $\beta\in\foB^+(Q_i), \foB^-(Q_i)$, or $\foB^0(Q_i)$ depending
on whether $\langle n_0,m_{\beta}\rangle >0$, $\langle n_0,m_{\beta}\rangle
<0$, or $\langle n_0,m_{\beta}\rangle=0$. This gives
decompositions
\begin{eqnarray*}
W_k(Q_1)&=&W_k^-(Q_1)+W^0_k(Q_1)+W_k^+(Q_1)\\
W_k(Q_2)&=&W_k^-(Q_2)+W^0_k(Q_2)+W_k^+(Q_2).
\end{eqnarray*}
We will show 
\begin{align}
\label{eqminus}
\theta_{\gamma,\foD}(W^-_k(Q_1)){} = & W_k^-(Q_2),\\
\label{eqplus}
\theta_{\gamma,\foD}^{-1}(W_k^+(Q_2)){} = &W_k^+(Q_1),\\
\label{eqzero}
W^0_k(Q_2){} =&W^0_k(Q_1).
\end{align}
From this follows the desired identity
\[
\theta_{\gamma,\foD}(W_k(Q_1))=W_k(Q_2),
\]
as $\theta_{\gamma,\foD}$ is necessarily the identity on
$W^0_k(Q_1)$.
One then uses this inductively to see that this holds
for any path $\gamma$ with endpoints in $M_{\RR}\setminus
\foU$ for which $\theta_{\gamma,\foD}$ is defined.

\medskip

\emph{Proof of \eqref{eqminus} and \eqref{eqplus}}.
If $\beta$ is a broken line with endpoint $Q_1$, then 
$\beta([t_{n-1},0])\cap L=\emptyset$ if 
$\langle n_0,m_{\beta}\rangle\le 0$, while $\beta([t_{n-1},0])\cap
L\not=\emptyset$ if $\langle n_0,m_{\beta}\rangle>0$. (Here we
are using $Q$ very close to $L$.) On the other hand, if $\beta$
has endpoint $Q_2$, then $\beta([t_{n-1},0])\cap L=\emptyset$ if
$\langle n_0,m_{\beta}\rangle\ge 0$ and $\beta([t_{n-1},0])\cap L 
\not=\emptyset$ if $\langle n_0,m_{\beta}\rangle <0$.

To see, say, \eqref{eqminus}, we proceed as follows.
Let $\beta\in\foB^-(Q_1)$. By the previous paragraph,
$\beta([t_{n-1},0])\cap L=\emptyset$. Let $c_nz^{m^{\beta}_{n}}$
be the monomial associated to the last segment of $\beta$, and
write
$\theta_{\gamma,\foD}(c_{n}z^{m^{\beta}_{n}})$ as a sum of
monomials $\sum_{i=1}^s d_iz^{m_i}$ as in Definition \ref{brokenlinedef},
(4). We can then deform $\beta$ continuously along $\gamma$ to
time $s_0$. Indeed, by the definition of $\foU$, if 
$\beta$ converged to a broken line through $\Sing(\foD)$, the image
of this broken line would be contained in $\foU$, and then $\foU$,
already containing $L$, would not be a manifold in a neighbourhood
of $\gamma(s_0)$.

\begin{figure}
\input{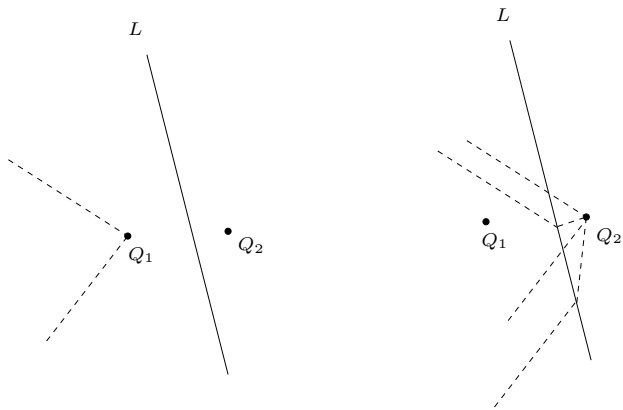}
\caption{Broken lines with endpoints $Q_1$ and $Q_2$.}
\label{bendycrossing}
\end{figure}

Let $\beta'$ be the deformation of $\beta$ with endpoint $\gamma(s_0)$.
For $1\le i\le s$, we then get a broken line $\beta_i'$ by adding
a short line segment to $\beta'$ in the direction $-r(m_i)$, with
attached monomial $d_iz^{m_i}$. This new broken line has endpoint
in $\fou_2$, and hence can be deformed to a broken line $\beta_i''
\in\foB^-(Q_2)$. We note that the line may not actually bend at $L$
if $d_iz^{m_i}$ is the term $c_{n}z^{m^{\beta}_n}$ appearing in
$\theta_{\gamma,\foD}(c_{n}z^{m^{\beta}_{n}})$. See Figure
\ref{bendycrossing}.

Conversely, any broken line $\beta\in\foB^-(Q_2)$ clearly
arises in this way. 

From this, \eqref{eqminus} becomes clear. \eqref{eqplus}
is identical. \qed

\medskip

\emph{Proof of \eqref{eqzero}}.
We will show that there are partitions
$\foB^0(Q_1)=\coprod_{i=1}^s \foB_i^1$ and
$\foB^0(Q_2)=\coprod_{i=1}^s \foB_i^2$ such that for each
$i$,
the contributions to $W_k(Q_1)$ and $W_k(Q_2)$ from
$\foB_i^1$ and $\foB_i^2$ are the same.

For simplicity of exposition, we will describe this very
explicitly in the case that any broken line with endpoint
$\gamma(s_0)$ passes through at most one point of $\Sing(\foD)$;
we leave it to the reader to deal with the general case: this
is notationally, but not conceptually, more complicated.

Let $\beta_1\in\foB^0(Q_1)$. If $\beta_1$ deforms continuously
to a broken line $\beta_2$ in $\foB^0(Q_2)$, then $\beta_1$ and
$\beta_2$ will each appear in one-element sets in the partition,
say $\beta_1\in\foB_i^1$, $\beta_2\in\foB_i^2$, and clearly
both these sets contribute the same term to $W_k(Q_1)$ and
$W_k(Q_2)$.

Now suppose $\beta_1\in\foB^0(Q_1)$ cannot be deformed continuously
to any $\beta_2\in\foB^0(Q_2)$. This means that there is a 
$B:(-\infty,0]\times [0,s_0]\rightarrow M_{\RR}$ as in Definition
\ref{brokenlinedef}, (5) such that $B|_{(-\infty,0]\times [0,s_0)}$
is a continuous deformation and $B_{s_0}$ is a broken line passing
through a point of $\Sing(\foD)$.
Furthermore, there must be some $j$ such that 
$B(t_j(s_0),s_0)=P\in\Sing(\foD)$.
In other words, for $s$ near $s_0$, $B_s$ must bend near $P$, and
this bending point approaches $P$. Otherwise,
we could continue to deform our broken line continuously through $P$ as the
line does not bend near $P$.

There are two cases
we need to analyze: either $P\in \{P_1,\ldots,P_k\}$ or
$P\not\in \{P_1,\ldots,P_k\}$.

\medskip

\begin{figure}
\input{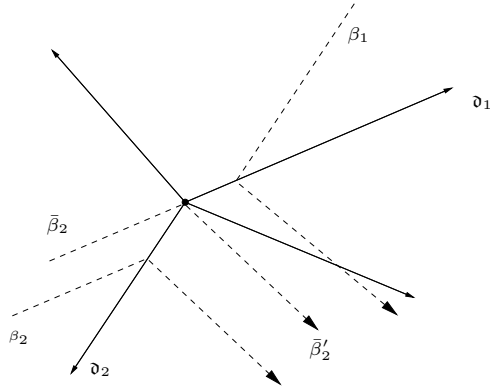}
\caption{The change in a broken line as it passes through a singular
point $P\in\{P_1,\ldots,P_k\}$.}
\label{secondcase}
\end{figure}

\emph{Case 1. $P=P_i$ for some $i$.} Because all rays with endpoint $P_i$
involve the monomial $u_i$, a broken line can only bend along at most
one ray with endpoint $P_i$, and as observed above,
$\beta_1$ must bend along at
least one such ray.
So call the ray with endpoint $P_i$ along which
$\beta_1$ bends $\fod_1$,
corresponding to a Maslov index zero
tree $\tilde h_1:\tilde
\Gamma_1\rightarrow M_{\RR}$. This tree passes through $P_i$,
and by cutting this tree at $P_i$ and removing the marked edge
mapping to $P_i$, we obtain a Maslov index two disk
$h_1:\Gamma'_1\rightarrow M_{\RR}$ with boundary $P_i$. Let
$\bar\beta_2$ be the broken line with endpoint $P_i$
corresponding to this Maslov index 2 disk. See Figure \ref{secondcase}.

Next, recalling that $B(t_j(s_0),s_0)=P_i$,
let $\bar\beta_2':[t_j(s_0),0]\rightarrow M_{\RR}$ be the restriction
of $B$ to $[t_j(s_0),0]\times \{s_0\}$: this is a piece of a broken
line starting at $P_i$. We can then concatenate $\bar\beta_2$
with $\bar\beta_2'$ by identifying $0$ in the domain of $\bar\beta_2$
with $t_j(s_0)$ in the domain of $\bar\beta_2'$, 
obtaining what we hope will be a broken line $\beta_2'$ passing through $P_i$.

Note that the broken line $B|_{(-\infty,t_j(s_0)]\times \{s_0\}}$
is a broken line with endpoint $P_i$, and hence corresponds to
a Maslov index two disk $h_2:\Gamma_2'\rightarrow M_{\RR}$ with
endpoint $P_i$. By extending the edge $E'_{\out}$ of
$\Gamma_2'$ to an unbounded edge, we get a tropical tree
$\tilde h_2:\tilde\Gamma_2\rightarrow M_{\RR}$, and once we mark
the point on $\tilde\Gamma_2$ which maps to $P_i$, it becomes
a Maslov index zero tree and hence corresponds to a ray $\fod_2\in
\foD$ with endpoint $P_i$.

Note that the function attached to $\fod_i$ is $1+w_{\Gamma_i}(E_{\out,i})u_i
\Mono(h_i)$. On the other hand, the monomial attached to
the last segment of $B|_{(-\infty,t_j(s_0)]\times \{s_0\}}$,
i.e. $c_jz^{m^{\beta_1}_j}$, is $\Mono(h_2)$, while the
monomial attached to the last segment of $\bar\beta_2$ is
$\Mono(h_1)$. Thus, in particular, the monomial $c_{j+1}z^{m^{\beta_1}_{j+1}}$
is obtained from the bend of $\beta_1$ at $\fod_1$, and hence is
\begin{equation}
\label{betamono}
w_{\Gamma_1}(E_{\out,1})\langle n_1,r(\Delta(h_2))\rangle u_i
\Mono(h_1)\Mono(h_2).
\end{equation}
Here $n_1\in N$ is primitive, orthogonal to $\fod_1$, and
positive on $r(\Delta(h_2))$.

We can now deform $\beta_2'$ by moving the endpoint of
$\bar\beta_2$ along $\fod_2$ away from $P_i$, moving $\bar\beta_2'$
along with it. However, we also need to keep track of monomials:
we have to make sure that the monomial on the first segment
of $\bar\beta_2'$ is the one which would arise when $\beta_2'$
bends along $\fod_2$.
However, this latter monomial is a term obtained by applying
the automorphism associated to crossing $\fod_2$ to $\Mono(h_1)$,
and is thus precisely
\begin{equation}
\label{beta2mono}
w_{\Gamma_2}(E_{\out,2})\langle n_2,r(\Delta(h_1))\rangle u_i
\Mono(h_1)\Mono(h_2).
\end{equation}
Again, $n_2\in N$ is primitive, orthogonal to $\fod_2$, and positive
on $r(\Delta(h_1))$. However, one sees easily that 
\[
w_{\Gamma_1}(E_{\out,1})\langle n_1,r(\Delta(h_2))\rangle=w_{\Gamma_2}
(E_{\out,2})\langle n_2,r(\Delta(h_1))\rangle,
\]
so \eqref{beta2mono} coincides with
\eqref{betamono}. As a result, $\beta_2'$ can now be deformed
away from the singular point $P_i$, giving a broken line $\beta_2$
with endpoint $Q_2$. Note that in no way does this represent
a continuous deformation: the broken line really jumps as it
passes through $P_i$.

Note this process is reversible. If we start with $\beta_2$
and try to deform it through $P_i$ as above, we obtain $\beta_1$.

To conclude, in this case, we can take one-element sets in the
partition of the form $\beta_1\in\foB^1_i$ and $\beta_2\in\foB^2_i$
for some $i$. They both give the same contribution to $W_k(Q)$.

\medskip

{\emph Case 2. $P\not\in \{P_1,\ldots,P_k\}$}. 
Let $\foD_P\subseteq \foD(\Sigma,P_1,\ldots,P_k)$ be the subset
of rays passing through $P$. Recall from the proof of Proposition 
\ref{gammaid} that we can write
\[
\foD_P=\{\fod_1,\ldots,\fod_n\}\cup \bigcup_{i=1}^m \foD_i
\]
where $\fod_1,\ldots,\fod_n$ are rays with $\Init(\fod_i)\not=P$
corresponding to trees $h$
such that $I(h)\cap I(h')\not=\emptyset$ for any Maslov index zero tree
$h'$ with outgoing edge
passing through $P$, and $\foD_1,\ldots,\foD_m$ are scattering
diagrams consisting of three rays, with two rays corresponding
to trees $h_1$ and $h_2$ with $I(h_1)\cap I(h_2)=\emptyset$ and
the third ray having initial point $P$, corresponding to the gluing
of the two trees $h_1$ and $h_2$ at $P$. 

Now $\beta_1$ cannot bend at two different rays near $P$ corresponding to 
trees $h_1$ and $h_2$ with $I(h_1)\cap I(h_2)\not=\emptyset$, as
$u_{I(h_1)}u_{I(h_2)}=0$. 
This gives us two immediate possibilities. 
The first possibility is that
$\beta_1$ only bends once near $P$ along a ray $\fod$ for which $P\not=
\Init(\fod)$. The second possibility is that there is some $i$
such that
$\beta$ either bends once along a ray $\fod\in\foD_i$
with $P=\Init(\fod)$ or $\beta$ bends at exactly two rays $\fod_1,\fod_2
\in\foD_i$ with $P\not=\Init(\fod_1),\Init(\fod_2)$.

The first possibility is easy to analyze: the singular point in fact
plays no role, and the broken line $\beta$ can be deformed continuously
through $P$. We were assuming this was not possible.

The second possibility is quite delicate. Since the broken line only
interacts with rays in some $\foD_i$, we can assume $\foD_i=\foD_P$,
so we are in the situation
depicted in Figure \ref{scatteringcones}, where we have precisely three
rays through $P$. We will write $\foD_P=\{\fod_1,\fod_2,\fod_3\}$
as depicted in that figure. 
We will find in this case that the broken lines we
are interested in cannot be deformed through $P$. Rather, one broken
line will split up into two broken lines, or two broken lines will
turn into one. 

We note $\foD_P$ splits $M_{\RR}$ up into five cones, labelled
$\sigma_1,\ldots,\sigma_5$ in Figure \ref{scatteringcones}. 
Now a broken line $\beta$
which bends twice near $P$ or bends along the ray with initial point $P$
will have an attached monomial $cz^m$ as $\beta$ approaches a neighbourhood
of $P$ and an attached monomial $c'z^{m+m_{\fod_1}+m_{\fod_2}}$ once $\beta$ leaves
the neighbourhood of $P$. Without loss of generality, assume $P=0$.
The behaviour of $\beta$ as it moves through $P$ will depend
on which of the cones $\sigma_1,\ldots,\sigma_5$ contain $r(m)$ and
$-r(m+m_{\fod_1}+m_{\fod_2})$.

\begin{figure}
\input{scatteringcones.pstex_t}
\caption{}
\label{scatteringcones}
\end{figure}

\begin{figure}
\input{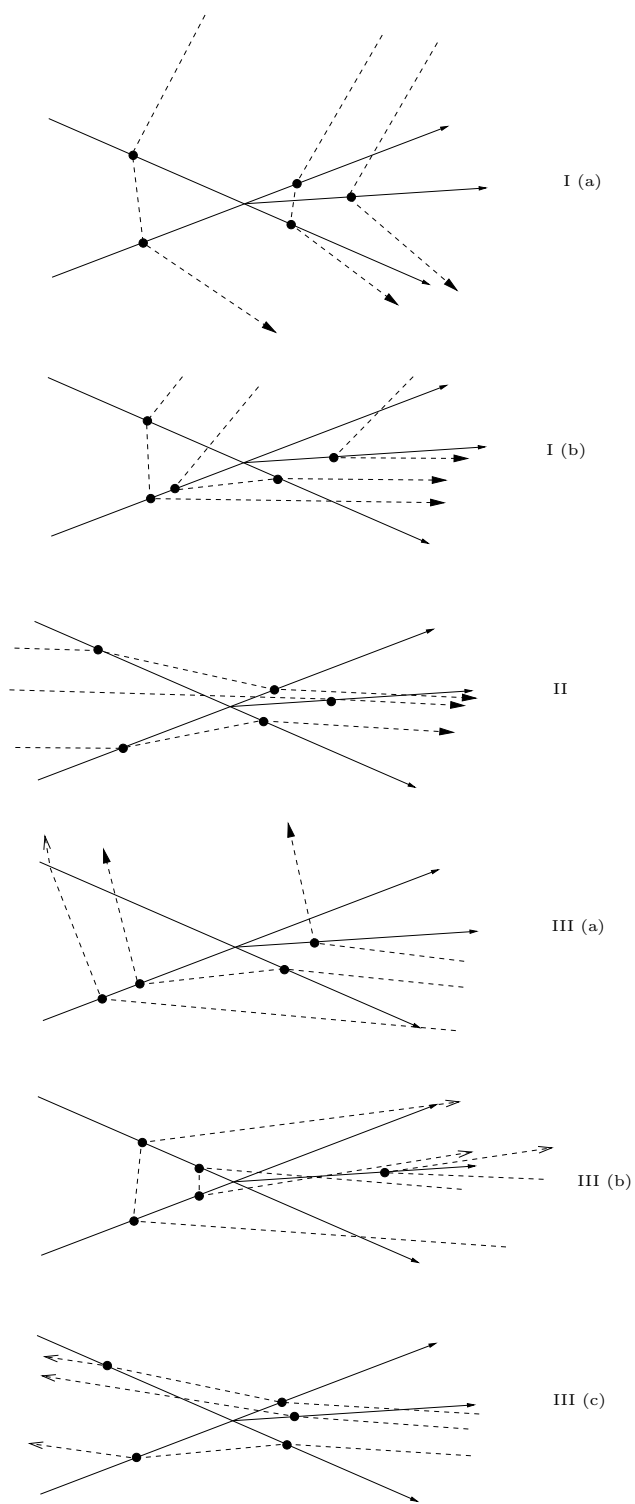}
\caption{Various possibilities for broken lines passing through singular
points. The dots indicate bending points.}
\label{casebycase}
\end{figure}

Figure \ref{casebycase} now gives a case by case description in the
case that $r(m)$ and $r(m+m_{\fod_1}+m_{\fod_2})$
are not tangent to any of the rays in $\foD_P$,
so $r(m)$ and $-r(m+m_{\fod_1}+m_{\fod_2})$ each lie in the interior of
one of $\sigma_1,\ldots,\sigma_5$.
In Case I, $r(m)\in\sigma_1$; in Case II, $r(m)\in\sigma_2$; and in
Case III, $r(m)\in\sigma_4$.
We can ignore the possibility of $r(m)$ being in $\sigma_3$ or
$\sigma_5$ as by symmetry, these behave in the same way as $\sigma_1$
or $\sigma_4$ respectively.

The finer classification comes from considering in which cone
$-r(m+m_{\fod_1}+m_{\fod_2})$ lies. We will go through this classification for
each of the three cases.

Consider first Case I. As $-r(m)$ lies in the cone $\sigma_3$ spanned
by $r(m_{\fod_2})$ and $-r(m_{\fod_1})$, $-r(m+m_{\fod_1}+m_{\fod_2})$ 
lies in the cone
spanned by $r(m_{\fod_2})$ and $-r(m_{\fod 1}+m_{\fod_2})$, i.e. the union of $\sigma_3$
and $\sigma_4$. In Case I (a), $-r(m+m_{\fod_1}+m_{\fod_2})\in \sigma_3$,
and in Case I (b), $-r(m+m_{\fod_1}+m_{\fod_2})\in\sigma_4$.

Consider next Case II. As $r(m)$ is in the interior of the cone generated by
$r(m_{\fod_1})$ and $r(m_{\fod_2})$, $-r(m+m_{\fod_1}+m_{\fod_2})$ lies in the interior of the
cone generated by $-r(m_{\fod_1})$ and $-r(m_{\fod_2})$. Hence $-r(m+m_{\fod_1}+m_{\fod_2})$
lies in either $\sigma_4$ or $\sigma_5$. 
The situation is symmetric with
respect to $\sigma_4$ and $\sigma_5$, so we only consider 
the case of $\sigma_4$. This gives Case II.

Finally, consider Case III. As $r(m)$ lies in the interior of the cone
generated by $-r(m_{\fod_1}+m_{\fod_2})$ and $-r(m_{\fod_1})$, $-r(m+m_{\fod_1}+m_{\fod_2})$
lies in the cone generated by $\pm r(m_{\fod_1}+m_{\fod_2})$ and $r(m_{\fod_1})$.
Thus $-r(m+m_{\fod_1}+m_{\fod_2})$ lies in $\sigma_1$, $\sigma_5$ or $\sigma_2$,
giving cases III (a), III (b), and III (c) respectively. 

Figure \ref{casebycase} illustrates all six of these possibilities.
In each case we draw all possible types, up to continuous
deformation, of broken lines near $P$
of the sort we are considering. There is always one that bends
precisely once near $P$, along $\fod_3$.
There are then always two other types of lines which bend twice, once
on each of $\fod_1$ and $\fod_2$. For example, consider Case I (a).
A broken line coming in from cone $\sigma_1$ hits $\fod_1$ or $\fod_2$
first, depending on its position. It would then bend as depicted,
and next hits $\fod_2$ or $\fod_1$, bending again. Note that
if a line hits $\fod_1$ first but does not bend there, but then
bends at $\fod_2$, the next line segment remains in the cone 
$\sigma_3$ and thus cannot hit $\fod_1$ again. Similarly
a broken line which hits $\fod_2$ first but fails to bend there
never hits $\fod_2$ again. Thus the possibilities depicted for
Case I (a) are the only possible broken lines which bend along both
$\fod_1$ and $\fod_2$ or along $\fod_3$.

The remaining cases are similar; we leave it to the reader to check
in each case that there are precisely two types of broken lines
which bend both at $\fod_1$ and $\fod_2$, as depicted.

Consider case I (a). The broken line $\beta_1$ is either the
line on the left or one of the two on the right. In the first
case, take $\foB^1_i=\{\beta_1\}$ for some choice of $i$, and
take $\foB^2_i$ to consist of the two broken lines on the right.
In the second case, we interchange this choice of $\foB^1_i$ and
$\foB^2_i$. We need only check that these two sets contribute
the same amount to $W_k$.

This in fact follows from Proposition \ref{gammaid}.
Indeed, consider two paths $\gamma_1$ and $\gamma_2$ with the same
endpoints, starting in $\sigma_1$ and ending in $\sigma_3$,
but with $\gamma_1$ going around the left of $P$ and
$\gamma_2$ going around the right. Then by Proposition \ref{gammaid},
$\theta_{\gamma_1,\foD}=\theta_{\gamma_2,\foD}$. 
Call the three broken lines in the figure for Case I (a), 
from left to right, $\beta_1$, $\beta_2$ and $\beta_3$, so
that $\foB^1_i=\{\beta_1\}$ and $\foB^2_i=\{\beta_2,\beta_3\}$
or vice versa.

Recalling that the monomial attached to the $j$-th segment
of each of these broken lines
is $cz^{m}$,
then the monomial attached to the $j+2$-nd
segment of $\beta_1$ is, by definition, the term of 
$\theta_{\gamma_1,\foD}(cz^{m})$ of the form
$c'z^{m+m_{\fod_1}+m_{\fod_2}}$. Similarly, the sum of the two
monomials attached to the
$j+2$-nd segment of $\beta_2$ and the $j+1$-st segment of $\beta_3$
is the term of
$\theta_{\gamma_2,\foD}(cz^{m})$ of the form
$c'z^{m+m_{\fod_1}+m_{\fod_2}}$, precisely because such terms appear
either when we cross $\fod_1$ and $\fod_2$, or when we cross $\fod_3$.
But since 
\[
\theta_{\gamma_1,\foD}(cz^{m})
=
\theta_{\gamma_2,\foD}(cz^{m})
\]
by Proposition \ref{gammaid}, we see that the contribution
to $W_k$ from $\foB^1_i$ and $\foB^2_i$ are the same.

The other cases are essentially the same: In Case I (b), we take
$\foB^1_i$ to be the two left-most lines and $\foB^2_i$ the right-most
line, or vice-versa. In Case II, we take $\foB^1_i$ to be the two
upper lines and $\foB^2_i$ to be the lower line, or vice versa.
In Case III (a), we partition by taking the two lower lines and
the upper line, and in Case III (c) we take the two upper
lines and the lower line. In each of these cases, the argument is
then identical to the argument given in Case I (a).

Case III (b) is slightly more delicate. Here we take $\foB^1_i$ be
the left-most line and $\foB^2_i$ to be the other two lines,
or vice versa. In this case, it is simplest to calculate the contributions
of the two sets of broken lines directly. Choose $n_1,n_2\in N$ primitive,
annihilating the tangent spaces to $\fod_1$, $\fod_2$ respectively,
and so that $\langle n_1,r(m_{\fod_2})\rangle >0$ and $\langle n_2,r(m_{\fod_1})
\rangle <0$.
Then $n_3=(w_1n_1+w_2n_2)/w_{\out}\in N$, annihilates the tangent space
to $\fod_3$, and $\langle n_3,r(m_{\fod_1})\rangle <0$, as is easily checked.

Label the broken lines
from left to right $\beta_1,\beta_2,\beta_3$.
Now the monomials attached to the $j+2$-nd segment of $\beta_1$ and $\beta_2$
and the $j+1$-st segment of $\beta_3$ are
\begin{align*}
\beta_1:\quad & \langle n_1,r(m)+r(m_{\fod_2})\rangle \langle n_2,r(m)\rangle
cc_1c_2w_1w_2 z^{m+m_{\fod_1}+m_{\fod_2}}\\
\beta_2:\quad & \langle n_2,r(m)+r(m_{\fod_1})\rangle \langle n_1,r(m)\rangle
cc_1c_2w_1w_2 z^{m+m_{\fod_1}+m_{\fod_2}}\\
\beta_3:\quad & \langle n_3,r(m)\rangle |r(m_{\fod_1})\wedge r(m_{\fod_2})|
cc_1c_2w_{\out} z^{m+m_{\fod_1}+m_{\fod_2}}
\end{align*}
But $|r(m_{\fod_1})\wedge r(m_{\fod_2})|=w_1\langle n_1,r(m_{\fod_2})\rangle=-w_2\langle n_2,
r(m_{\fod_1})\rangle$.
So the total contribution from $\beta_2$ and $\beta_3$ is, leaving off
the common term $cc_1c_2z^{m+m_{\fod_1}+m_{\fod_2}}$,
\begin{align*}
&\langle n_2,r(m)+r(m_{\fod_1})\rangle\langle n_1,r(m)\rangle w_1w_2+
\langle n_{3},r(m)\rangle |r(m_{\fod_1})\wedge r(m_{\fod_2})|w_{\out}\\
= {} & \big(\langle n_2,r(m)\rangle\langle n_1,r(m)\rangle w_1w_2
+\langle n_2,r(m_{\fod_1})\rangle\langle n_1,r(m)\rangle w_1w_2\big)\\
&+(-\big \langle n_1,r(m)\rangle\langle n_2,r(m_{\fod_1})\rangle w_1w_2
+\langle n_2,r(m)\rangle\langle n_1,r(m_{\fod_2})\rangle w_1w_2\big)\\
= {} & \big(\langle n_2,r(m)\rangle\langle n_1,r(m)\rangle
+\langle n_2,r(m)\rangle\langle n_1,r(m_{\fod_2})\rangle\big)w_1w_2\\
= {} & \langle n_1,r(m)+r(m_{\fod_2})\rangle\langle n_2,r(m)\rangle w_1w_2.
\end{align*}
Thus the contribution from $\beta_1$ agrees with the sum of the contributions
from $\beta_2$ and $\beta_3$.

There is one degenerate situation we still have to deal with, namely
that one of $r(m)$ or $-r(m+m_{\fod_1}+m_{\fod_2})$
are parallel to one of the rays in $\foD_P$. These can be viewed
as a degenerate version of one of the cases considered above,
and requires another
case-by-case analysis, which we shall leave to the reader.
In all cases, one of the broken lines disappears.
Let us illustrate what happens in one case, in which, in Case I (a),
the incoming segment of the broken lines is parallel to $\fod_2$. In
this case, we only have two possible broken lines, as depicted as in
Figure \ref{parallelcase}. There is no broken line which starts
in $\sigma_5$ and bends first at $\fod_2$ and then at $\fod_1$.

This does not present any difficulty for our analysis, however.
As before, pick two paths $\gamma_1$ and $\gamma_2$ going from
$\sigma_1$ to $\sigma_3$, with $\gamma_1$ on the left. Note that
from the explicit analysis given in Case I (a), the broken line
on the left, $\beta_1$, contributes the $z^{m+m_{\fod_1}+m_{\fod_2}}$ term
to $\theta_{\gamma_1,\foD}(cz^m)$, while the broken line
on the right provides the only contribution to the $z^{m+m_{\fod_1}+m_{\fod_2}}$
term in $\theta_{\gamma_2,\foD}(cz^m)$, precisely because
$\theta_{\gamma_2,\fod_2}(z^m)=0$.

The same happens in all parallel cases, and we omit the details.
\qed

\begin{figure}
\input{parallelcase.pstex_t}
\caption{}
\label{parallelcase}
\end{figure}

\bigskip

We now turn our attention to studying the dependence of
$W_k(Q)$ on the points $P_1,\ldots,P_k$. 
For this discussion, it will be useful to explain that the automorphisms
arising from scattering diagrams are elements of an interesting
group, variants of which first appeared in \cite{ks} and then was used
in \cite{GS} and \cite{GPS}. We summarize this point of view.

We denote the \emph{module of log derivations} of $\CC[T_{\Sigma}]
\otimes_{\CC} R_k\lfor y_0\rfor$ to be the module
\[
\Theta(\CC[T_{\Sigma}]\otimes_{\CC} R_k\lfor y_0\rfor):=
\Hom_{\ZZ}(M,\CC[T_{\Sigma}]\otimes_{\CC} R_k\lfor y_0\rfor)=
(\CC[T_{\Sigma}]\otimes_{\CC} R_k\lfor y_0\rfor))\otimes_{\ZZ} N.
\]
An element $f\otimes n$ is written as $f\partial_n$, and acts
as a derivation on $\CC[T_{\Sigma}]\otimes_{\CC} R_k\lfor y_0\rfor$
over $\CC[K_{\Sigma}]\otimes_{\CC}R_k\lfor y_0\rfor$ via
\[
f\partial_n(z^m)=f\langle n,r(m)\rangle z^m.
\]
Given $\xi\in \fom_{R_k}\Theta(\CC[T_{\Sigma}]\otimes_{\CC} R_k
\lfor y_0\rfor)$, where $\fom_{R_k}=(u_1,\ldots,u_k)$ is the maximal
ideal of $R_k$, we define
\[
\exp(\xi)\in\Aut(\CC[T_{\Sigma}]\otimes_{\CC} R_k\lfor y_0\rfor)
\]
by
\[
\exp(\xi)(a)=a+\sum_{i=1}^{\infty} {\xi^i(a)\over i!}.
\]
This is a finite sum given the assumption on $\xi$.

Now let
\[
\fov_{\Sigma,k}=\bigoplus_{m\in T_{\Sigma}\atop r(m)\not=0}
z^m(\fom_{R_k}\otimes r(m)^{\perp})\subseteq \Theta(\CC[T_{\Sigma}]
\otimes_{\CC} R_k\lfor y_0\rfor).
\]
Set 
\[
\VV_{\Sigma,k}=\{\exp(\xi)\,|\,\xi\in\fov_{\Sigma,k}\}.
\]
Note the Lie bracket on $\fov_{\Sigma,k}$ is given by
\[
[z^m\partial_n,z^{m'}\partial_{n'}]=z^{m+m'}\big(
\langle n,r(m')\rangle \partial_{n'}-\langle n',r(m)\rangle \partial_n
\big).
\]
Then $\fov_{\Sigma,k}$ is closed under Lie bracket, and
hence $\VV_{\Sigma,k}$ is a group, with multiplication
given by the Baker-Campbell-Hausdorff formula.
We note that for $m\in T_{\Sigma}$ with $r(m)\not=0$, $n\in N$ with
$\langle n,r(m)\rangle=0$, $I\subseteq\{1,\ldots,k\}$ non-empty, $c\in\CC$,
\[
\exp(cu_Iz^m\partial_n)(z^{m'})
=z^{m'}(1+cu_I\langle n,r(m')\rangle z^m),
\]
so for any scattering diagram $\foD$, we have
$\theta_{\gamma,\foD}\in\VV_{\Sigma,k}$.
Furthermore, $\VV_{\Sigma,k}$ is generated by automorphisms of the form
$\exp(cu_Iz^m\partial_n)$.

The original version of this group introduced in \cite{ks} was defined
as a group of Hamiltonian symplectomorphisms, and it is convenient to
use this identification here. In particular, the holomorphic symplectic form
$\Omega$ was given by a choice of generator of $\bigwedge^2 M$,
i.e., an identification
$\bigwedge^2 M\cong\ZZ$. This gives an isomorphism
$M\cong N$ well-defined up to sign; we view $m\in M$
as an element $X_m$ of $N$ via the linear map $m'\mapsto m\wedge m'
\in\bigwedge^2 M\cong\ZZ$. 

Suppose $\fod$ is a ray
or line with $f_{\fod}=1+w(m)cu_Iz^m$ where $w(m)$ is the index
of $r(m)$. Then the automorphism $\theta$ obtained by crossing $\fod$ is
$\exp(\pm cu_I z^m X_{r(m)})$. 
In fact $z^mX_{r(m)}$ is the Hamiltonian vector field\footnote{Given
a function $f$ on $M\otimes \CC^{\times}$, the Hamiltonian
vector field associated to it is the vector field $X_f$ such that
$\iota(X_f)\Omega=df$.}
associated
to the function $f=-z^m$. 
So for $f\in \fom_{R_k}\big(\CC[T_{\Sigma}]\otimes_{\CC} 
R_k\lfor y_0\rfor\big)$, write 
\[
X_f
\in\fom_{R_k}\Theta(\CC[T_{\Sigma}]\otimes_{\CC} R_k\lfor y_0\rfor)
\]
for the Hamiltonian vector field induced by $f$.
This is convenient for writing the following
easily checked standard lemma: 

\begin{lemma}
\label{adjointlemma}
If $f\in \fom_{R_k}(\CC[T_{\Sigma}]\otimes_{\CC} R_k\lfor y_0\rfor)$
and $\theta\in\VV_{\Sigma,k}$, then 
\[
\theta\circ X_f\circ\theta^{-1}=X_{\theta(f)}.
\]
\end{lemma}

%

We will need a three-dimensional version of scattering diagrams.

\begin{definition}
Let $L\subseteq \RR$ be a closed interval. Let $\pi_1$ and $\pi_2$
be the projections of $M_{\RR}\times L$ onto $M_{\RR}$ and $L$ respectively.
A \emph{scattering
diagram} in $M_{\RR}\times L$ is a set $\foD$ consisting of
pairs $(\fod,f_{\fod})$ such that
\begin{itemize}
\item $\fod\subseteq M_{\RR}\times L$ is a convex polyhedral
subset of dimension two such that $\pi_2(\fod)$ is one-dimensional. 
Furthermore there is a one-dimensional
subset $\fob\subseteq M_{\RR}\times L$ and an element
$m_0\in T_{\Sigma}$ with $r(m_0)\not=0$ such that
\[
\fod=\fob-\RR_{\ge 0}(r(m_0),0).
\] 
\item $f_{\fod}\in \CC[z^{m_0}]\otimes_{\CC} R_k
\subseteq \CC[T_{\Sigma}]\otimes_{\CC} R_k\lfor y_0\rfor$.
\item $f_{\fod}\equiv 1\mod (u_1,\ldots,u_k)z^{m_0}$.
\end{itemize}
We define
\[
\Sing(\foD)=\bigcup_{\fod\in\foD} \partial\fod
\cup \bigcup_{\fod_1,\fod_2\in\foD\atop
\dim\fod_1\cap\fod_2=1} \fod_1\cap\fod_2.
\]
This is a one-dimensional subset of $M_{\RR}\times L$.
Let $\Interstices(\foD)$ be the finite set of points where
$\Sing(\foD)$ is not a manifold. In keeping with the language of
\cite{GS}, we will denote by $\Joints(\foD)$ the set
of closures of the connected components of $\Sing(\foD)
\setminus\Interstices(\foD)$, calling elements of 
$\Joints(\foD)$ and $\Interstices(\foD)$ \emph{joints}
and \emph{interstices} respectively. We call a joint 
\emph{horizontal} if its image under $\pi_2$ is a point;
otherwise we call a joint \emph{vertical}.

For a path $\gamma$ in $(M_{\RR}\times L)\setminus \Sing(\foD)$,
one can define an element
\[
\theta_{\gamma,\foD}\in\VV_{\Sigma,k}
\]
exactly as in the case of a scattering diagram in $M_{\RR}$.
Indeed, we just need to define the automorphism $\theta_{\gamma,\fod}$
when $\gamma$ crosses $(\fod,f_{\fod})$ at time $t_i$.
Assuming $\gamma$ crosses $\fod$ transversally, then $\pi_{1*}(\gamma'(t_i))$
is not parallel to $r(m_0)$. So we can choose $n_0\in N$
primitive with $\langle n_0,r(m_0)\rangle=0$ and $\langle n_0,
\pi_{1*}(\gamma'(t_i))\rangle<0$. We define as usual
\[
\theta_{\gamma,\fod}(z^m)=z^mf_{\fod}^{\langle n_0,r(m)\rangle}.
\]
Again, it is easy to check that $\theta_{\gamma,\foD}$ only
depends on the homotopy type of the path $\gamma$ inside
$(M_{\RR}\times L)\setminus \Sing(\foD)$.

A \emph{broken line} in $M_{\RR}\times L$ is a map
$\beta:(-\infty,0]\rightarrow M_{\RR}\times L$, along with data
$t_0<\cdots<t_n$ and monomials $c_iz^{m_i^{\beta}}$, such that
\begin{enumerate}
\item
$\pi_2\circ\beta$ is constant, say with image $P\in L$. 
\item $\pi_1\circ\beta$ is a broken line
in the sense of Definition \ref{brokenlinedef} with respect
to the scattering diagram $\foD_P$ in $M_{\RR}$ given, after identifying
$M_{\RR}\times \{P\}$ with $M_{\RR}$, by
\[
\foD_P:=\big\{\big(\fod\cap (M_{\RR}\times \{P\}),f_{\fod}\big)\,\big|\,
\hbox{$(\fod,f_{\fod})\in \foD$ such that $\fod\cap (M_{\RR}\times \{P\})
\not=\emptyset$}\big\}.
\]
\end{enumerate}
\end{definition}

We now reach the last goal of our section: we wish to understand
how $W_k(Q)$ varies as the points $P_1,\ldots,P_k$ are varied.

\begin{theorem}
\label{Pindeptheorem}
Let $W$ and $W'$ be $W_k(Q)$ for two different choices of general
points $P_1,\ldots,P_k$ and $P_1',\ldots,P_k'$. Then 
\[
W'=\theta(W)
\]
for some $\theta\in\VV_{\Sigma,k}$.
\end{theorem}

\proof We shall show this result by induction on $k$, noting that the
base case $k=1$ is obvious, as moving $P_1$ and keeping $Q$
fixed is the same thing as moving $Q$ and keeping $P_1$ fixed.

It is clearly enough to show this result in the case that only $P_1$
changes. So consider
a choice of general points $P_1,\ldots,P_k$ and $P_1'$.
Consider the line segment $L$ joining $P_1$ and $P_1'$. For all but 
a finite number of points $P\in L$, we can assume $P,P_2,\ldots,P_k$
will be sufficiently general so that $\Trees(\Sigma,P,P_2,\ldots,P_k)$
is finite, and all elements of this set are trivalent. This gives rise
to a family of scattering diagrams $\foD(\Sigma,P,P_2,\ldots,P_k)$,
$P\in L$. We can put these scattering diagrams together into
a scattering diagram 
$\tilde\foD=\tilde\foD(\Sigma,L,P_2,\ldots,P_k)$ in $M_{\RR}
\times L$.
$\tilde\foD$ is determined by the requirement
that for $P\in L$ general,
\[
\foD(\Sigma,P,P_2,\ldots,P_k)=\big\{\big(\tilde\fod\cap (M_{\RR}\times \{P\}),
f_{\tilde\fod}\big)\,\big|\, \hbox{$(\tilde\fod, f_{\tilde\fod})\in
\tilde\foD$ such that $\tilde\fod\cap (M_{\RR}\times \{P\})\not=\emptyset$}
\big\}.
\]

To keep track of the dependence of $W_k(Q)$ on the point $P\in L$, we
write $W_k(Q;P)$. We wish now to show that if $\gamma$ is a general path in
$M_{\RR}\times L$ joining $(Q,P_1)$ to $(Q,P_1')$, then 
\begin{equation}
\label{Wktransport}
W_k(Q;P_1')=\theta_{\gamma,\tilde\foD}(W_k(Q;P_1)).
\end{equation}
It is enough to show 
\begin{enumerate}
\item $W_k(Q;P)$ is constant for $(Q,P)\in M_{\RR}\times L$
varying within a connected component of $(M_{\RR}\times L)
\setminus \Supp(\tilde\foD)$.
\item For two such connected components separated by a wall
$(\tilde\fod,f_{\tilde\fod})$ and points $(Q,P)$, $(Q',P')$ on
either side of the wall, we have \eqref{Wktransport} for $\gamma$
a short path joining $(Q,P)$ with $(Q',P')$.
\end{enumerate}

Once we show (1), Theorem \ref{wallcrossingtheorem} already shows (2):
as there are no walls in $\tilde\foD$ projecting to points in $L$,
we can always choose points $(Q,P)$, $(Q',P')$ on opposite sides
of a wall with $P=P'$, and then we are in the case already shown in Theorem
\ref{wallcrossingtheorem}. So we only need to show (1).

To show (1), we use the same technique we used for the variation of
$Q$, deforming broken lines. 
Take $(Q,P)$ and $(Q',P')$ general within a connected
component of $(M_{\RR}\times L)\setminus \Supp(\tilde\foD)$ and move from
$(Q,P)$ to $(Q',P')$ via a general path $\gamma$. Consider broken
lines in $M_{\RR}\times L$ 
with endpoint $\gamma(t)$. As $t$ varies, we can continuously
deform a broken line with endpoint $\gamma(t)$ unless the broken
line converges to one passing through a singular point of
$\tilde\foD$. However, since such a family of broken lines
traces out a two-dimensional subset of $M_{\RR}\times L$,
by choosing $\gamma$ sufficiently general, we can be sure
that none of these broken lines converge to broken lines
passing through interstices,
as interstices are codimension three. 
However, they can pass through
joints, and this requires some care.

The first observation
is that we have already analyzed in the proof of Theorem 
\ref{wallcrossingtheorem} what happens if a broken line passes through a
vertical joint. Indeed, we can just as well assume that
$\gamma$ has been chosen so that at a time $t_0$ when 
a broken line passes through a vertical joint, $\pi_2(\gamma(t))$
remains constant for $t$ in a neighbourhood of $t_0$. Then
we are in precisely the situation analyzed in Theorem 
\ref{wallcrossingtheorem}.

So we only need to see what happens if a broken line
passes through a horizontal joint. Note that horizontal
joints occur when two or more
parallel rays in a scattering diagram come together as the
point $P$ varies; this can typically lead to values of $P$ with families
of Maslov index zero disks or the existence of Maslov index $-2$ disks.

In fact, it is enough to show that if $\foj$
is a horizontal joint and $\gamma_{\foj}$ is a small loop in
$M_{\RR}\times L$ around the joint,
then $\theta_{\gamma_{\foj},\tilde\foD}=\id$. Indeed, if $\foj$ projects
to $P\in L$, $\foj$ is contained in some polygons $\tilde\fod_1,
\ldots,\tilde\fod_n\in\tilde\foD$, and necessarily for $P'\in L$
near $P$, $\tilde\fod_i\cap (M_{\RR}\times\{P'\})$ is either a ray
parallel to $\foj$ or is empty. Thus
as $P'\in L$ moves from one side of $P$ to
the other, some parallel rays $\fod_1,\ldots,\fod_p$ in 
$\foD(\Sigma,P',P_2,\ldots,P_k)$
come together to yield the joint and then turn into parallel rays $\fod'_1,
\ldots,\fod'_{p'}$ on the other side of $P$. Let $\foD_1,\foD_2$
be the scattering diagrams in $M_{\RR}$ given by $\foD(\Sigma,P',P_2,
\ldots,P_k)$ for $P'$ very close to $P$, but on opposite sides of $P$.
Let $\gamma$ be a path which is a short line segment crossing $\foj$,
so that we can write
\begin{align*}
\theta_{\gamma,\foD_1}{} =&\theta_{\gamma,\fod_1}\circ\ldots
\theta_{\gamma,\fod_p},\\
\theta_{\gamma,\foD_2}{} =&\theta_{\gamma,\fod'_1}\circ\ldots
\theta_{\gamma,\fod'_{p'}}.
\end{align*}
(Note the ordering is immaterial as all these automorphisms commute).
But
$\theta_{\gamma_{\foj},\tilde\foD}=\theta_{\gamma,\foD_1}^{-1}\circ
\theta_{\gamma,\foD_2}$, so if $\theta_{\gamma_{\foj},\tilde\foD}=
\id$, we have $\theta_{\gamma,\foD_1}=\theta_{\gamma,\foD_2}$.
This means that, by Definition \ref{brokenlinedef}, (4), broken lines
will behave in the same way on either side of $P$ near the joint $\foj$.
Note that the actual set of broken lines on either side may be different,
because we are not claiming that the set $\{\theta_{\gamma,\fod_1},
\ldots,\theta_{\gamma,\fod_p}\}$ coincides with $\{\theta_{\gamma,\fod_1'},
\ldots,\theta_{\gamma,\fod'_{p'}}\}$, but rather the total contribution from
bends along the two sets of broken lines remains the same.

To show $\theta_{\gamma_{\foj},\tilde\foD}=\id$ for each horizontal 
joint, we use
a method introduced in \cite{GS}, Proposition 3.17.
For $I\subseteq \{1,\ldots, k\}$, define
\[
\Ideal(I)
:=\langle u_i | 
i\not\in I\rangle\subseteq \CC[T_{\Sigma}]\otimes_{\CC} R_k\lfor y_0\rfor.
\]
We will proceed by induction, showing 

\medskip

\emph{Claim}. For $k'\ge 0$ and $\#I=k'$,
we have
\[
\theta_{\gamma_{\foj},\tilde\foD}\equiv\id\mod 
\Ideal(I)
\]
for every horizontal joint $\foj$.
\medskip

The base case with $k'=0$ is trivial, because all automorphisms
are trivial modulo the ideal $(u_1,\ldots,u_k)$. 
So assume the claim for all $k''<k'$.
Fix a set $I$ with $\#I=k'$.
Fix an orientation on $M_{\RR}\times L$, so that if any
joint $\foj$ is given an orientation, this determines the orientation
of a loop $\gamma_{\foj}$ around $\foj$, using, say, the right-hand rule.
We wish to study $\theta_{\gamma_{\foj},\tilde\foD}$ for $\foj$
horizontal.

Note that as 
$\theta_{\gamma_{\foj},\tilde\foD}$ for $\foj$ horizontal
only involves a composition of
automorphisms associated to parallel rays, we can in fact write
\[
\theta_{\gamma_{\foj},\tilde\foD}(z^{m'})=
f_{\foj}^{\langle n_{\foj},r(m')\rangle}z^{m'}
\]
for some $n_{\foj}\in N$ primitive and zero on the tangent space
to $\foj$. Also,
\[
f_{\foj}\in\CC[\{m\in T_{\Sigma}|\hbox{$r(m)$ is tangent to $\foj$}\}]
\otimes_{\CC} R_k\lfor y_0\rfor.
\]
Note that $f_{\foj}$ depends on the choice of sign of $n_{\foj}$. 
Assume we have chosen these consistently, in the sense that if any two joints
$\foj$ and $\foj'$ have the same tangent space, then $n_{\foj}=n_{\foj'}$.

We need to show $f_{\foj}\equiv 1\bmod \Ideal(I)$. 
Fix some $m\in T_{\Sigma}$.
For each horizontal joint $\foj$,
let the term in $f_{\foj}\bmod \Ideal(I)$ involving $z^m$ be $c_{m,\foj}z^m$.
Thus $c_{m,\foj}=\bar c_{m,\foj} \prod_{i\in I} u_i$
for some $\bar c_{m,\foj}\in\CC$ since $f_{\foj}\equiv 1\mod
\Ideal(I')$ for any $I'\subsetneq I$. 
Note that $f_{\foj}$ is
a product of polynomials of the form $1+cz^{m'}$ with $r(m')\not=0$.
It then follows
that if $\bar c_{m,\foj}\not=0$, then $r(m)\not=0$. So we will assume
$r(m)\not=0$.
We will also include here the case that $\foj$ is a vertical joint, by setting
$\bar c_{m,\foj}=0$ for vertical joints. 
Note that $\bar c_{m,\foj}$ depends on the orientation
on $\foj$. A change of orientation of $\foj$ changes the direction of
$\gamma_{\foj}$, replacing $f_{\foj}$ with $f_{\foj}^{-1}$.
This changes the sign of
$\bar c_{m,\foj}$. As a result, we can view 
\[
\foj\mapsto \bar c_{m,\foj}
\]
as a 1-chain for the one-dimensional simplicial complex $\Sing(\tilde\foD)$. 
Here the choice of orientation on $\foj$ implicit.

\medskip

\emph{Subclaim}. $\foj\mapsto\bar c_{m,\foj}$ is a 1-cycle.

\proof We need to check the 1-cycle condition at each interstice
of $\tilde\foD$, so let $(Q,P)\in \Interstices(\tilde\foD)$. Consider
a small two-sphere $S$ in $M_{\RR}\times L$ with center $(Q,P)$. Then
suppose that $x_1,\ldots,x_s\in S$ are distinct points such that
\[
\{x_1,\ldots,x_s\}=\bigcup_{\foj\in \Joints(\tilde\foD)} \foj\cap S.
\]
Choose a base-point $y\in S$, $y\not\in \Supp(\tilde\foD)$.
We can choose small counterclockwise loops $\gamma_1,\ldots,
\gamma_s$ in $S$ around $x_1,\ldots,x_s$ and paths $\beta_i$ joining
$y$ with the base-point of $\gamma_i$ in such a way that
\[
\beta_1\gamma_1\beta_1^{-1}\ldots\beta_s\gamma_s\beta_s^{-1}=1
\]
in $\pi_1(S\setminus\{x_1,\ldots,x_s\},y)$. Because
$\theta_{\gamma,\tilde\foD}$ only depends on the homotopy type of
$\gamma$ in $(M_{\RR}\times L)\setminus\Sing(\tilde\foD)$, we obtain
the equality
\begin{equation}
\label{standardrel}
\theta_{\beta_s}^{-1}\circ\theta_{\gamma_s}\circ\theta_{\beta_s}
\circ\cdots\circ\theta_{\beta_1}^{-1}\circ\theta_{\gamma_1}\circ
\theta_{\beta_1}=\id.
\end{equation}
Here, we have dropped the $\tilde\foD$'s in the subscripts.

We now distinguish between two cases.

\medskip

\emph{Case 1}. The interstice $(Q,P)$ does not satisfy $Q\in \{P,P_2,
\ldots,P_k\}$. Then by Proposition \ref{gammaid}, $\theta_{\gamma_i}=
\id$ for each $\gamma_i$ which is a loop around a vertical joint
containing $(Q,P)$. 
On the other hand, modulo $\Ideal(I)$, for $\gamma_i$
around a horizontal joint $\foj_i$, by the induction hypothesis,
$f_{\foj_i}$ is of the form $1+(\cdots)\prod_{i\in I} u_i$. 
One then checks 
$\theta_{\gamma_i}$ necessarily commutes, modulo $\Ideal(I)$, with 
any element of $\VV_{\Sigma,k}$. This can easily be seen as in Example
\ref{basicscatterexample}, using the fact that $u_j\prod_{i\in I} u_i
\equiv 0\mod I$ for any $j$. Thus in particular, $\theta_{\gamma_i}$
commutes with
$\theta_{\beta_i}$. Thus \eqref{standardrel} becomes
\[
\prod \theta_{\gamma_i}\equiv\id\mod\Ideal(I)
\]
where the product is over all $\gamma_i$ around horizontal joints. Applying
this identity to a monomial $z^{m'}$, we obtain
\[
\prod f_{\foj_i}^{\langle n_{\foj_i},r(m')\rangle} z^{m'}=
z^{m'}\mod\Ideal(I),
\]
which, after expansion, gives the identity
\begin{equation}
\label{standardrel2}
\sum \langle n_{\foj_i},r(m')\rangle \bar c_{m,\foj_i}=0\mod\Ideal(I)
\end{equation}
for any $m'\in T_{\Sigma}$.
Now a monomial $z^m$ can only appear in $f_{\foj_i}$
if $r(m)$ is in fact tangent to $\foj_i$, so the only horizontal
joints containing
$(Q,P)$ with $\bar c_{m,\foj}\not=0$ are the joints contained in the
affine line $(Q,P)+\RR (r(m),0)$. Let $s$ be the number of joints
contained in this line containing $(Q,P)$. Then either $s=0,1$ or $2$.
If $s=0$, there is nothing to prove. If $s=1$, with $\foj_i$ the
only such joint,
it follows from \eqref{standardrel2} that $\bar c_{m,\foj_i}=0$.
If $s=2$, let
$\foj_{i_1},\foj_{i_2}$ be the two such joints. Then
\eqref{standardrel2} implies that
$\bar c_{m,\foj_{i_1}}=\bar c_{m,\foj_{i_2}}$,
assuming $\foj_{i_1}$ and $\foj_{i_2}$ are oriented in the same direction.
This shows that the 1-cycle condition holds at $(Q,P)$.

\medskip

\emph{Case 2}. The interstice $(Q,P)$ satisfies $Q\in \{P,P_2,\ldots,P_k\}$,
say $Q=P_i$. We'll write $i=1$ if $Q=P$. The argument is almost the same,
but now there are two vertical joints, say $\foj_1$ and $\foj_2$ with endpoint
$(Q,P)$, with $\foj_1,\foj_2\subseteq \{P_i\}\times L$ if $i>1$ and
$\foj_1,\foj_2\subseteq \{(P',P')|P'\in L\}\subseteq M_{\RR}\times L$
if $i=1$.
Without loss of generality
we can take the base-point $y$ near $x_1$ and assume $\beta_1$ is
a constant path so that $\theta_{\beta_1}=\id$. 
The argument will be the same as in Case 1 once we show
\begin{equation}
\label{newrel}
\theta_{\beta_2}^{-1}\circ\theta_{\gamma_2}\circ\theta_{\beta_2}
\circ\theta_{\gamma_1}=\id.
\end{equation}

To do so, consider the scattering diagram $\foD(\Sigma,P',P_2,\ldots,
P_k)$ for $P'\in L$, $P'$ near $P$ but $P'\not=P$. 
By Remark \ref{Piremark}, the rays emanating
from $P_i$ ($P'$ if $i=1$) in $\foD(\Sigma,P',P_2,\ldots,P_k)$ are
in one-to-one correspondence with the terms in $W_{k-1}(P_i;P')-y_0$,
where $W_{k-1}(P_i;P')$ denotes $W_{k-1}(P_i)$ computed using the marked points
$P',\ldots,{P_{i-1}},P_{i+1},\ldots P_k$ (or
$P_2,\ldots,P_k$ if $i=1$). In particular, given a term
$cz^m$ in $W_{k-1}(P_i;P')-y_0$, the corresponding ray carries the function
$1+u_icw(m)z^m$, where $w(m)$ is the index of $r(m)$. 
Note that if $\gamma$ is a simple loop around $P_i$, then the contribution
to $\theta_{\gamma,\foD(\Sigma,P',P_2,\ldots,P_k)}$ from such
a ray is $\exp(\pm X_{u_icz^m})$. Here the sign only depends
on the orientation of $\gamma$ and the chosen identification of $\bigwedge^2
M$ with $\ZZ$. 
All automorphisms attached to the rays
emanating from $P_i$ commute by Example \ref{basicscatterexample} because
$u_i^2=0$, so
\[
\theta_{\gamma,\foD(\Sigma,P',P_2,\ldots,P_k)}
=\prod \exp(\pm X_{u_icz^m})
=\exp(\pm X_{u_i(W_{k-1}(P_i;P')-y_0)}).
\]
Here the product is over all terms $cz^m$ appearing in $W_{k-1}(P_i;P')-y_0$.
Furthermore, if $P'\in \pi(\foj_1)\setminus\{P\}$ 
and $P''\in \pi(\foj_2)\setminus\{P\}$,
then by \eqref{Wktransport} applied inductively
to $k-1$ points if $i\not=1$, and by Theorem
\ref{wallcrossingtheorem}, if $i=1$,
\[
u_iW_{k-1}(P_i;P'')=u_i\theta_{\beta_2}(W_{k-1}(P_i;P')).
\]
It then follows from Lemma \ref{adjointlemma} that
\[
\theta_{\gamma_2}=\big( \theta_{\beta_2}\circ\theta_{\gamma_1}\circ
\theta_{\beta_2}^{-1}\big)^{-1},
\]
the last inverse on the right since $\gamma_1$ and $\gamma_2$ are homotopic
to loops in $M_{\RR}\times \{P'\}$ and $M_{\RR}\times \{P''\}$
respectively with \emph{opposite orientations}. This shows 
\eqref{newrel}. We can then finish as in Case 1. 

This completes the proof of the subclaim.
\qed

\medskip

To complete the proof of the claim, hence the theorem, 
we now note that the cycle $\sigma$
given by
$\foj\mapsto \bar c_{m,\foj}$ is in fact zero. Indeed, picking a given
joint with $\bar c_{m,\foj}\not=0$, the fact that $\sigma$ is a cycle
implies that the line containing $\foj$ can be written as a union
of joints $\foj'$ with orientation compatible with that on $\foj$,
with $\bar c_{m,\foj'}=\bar c_{m,\foj}$. However, there must be one
joint $\foj'$ contained in this line which is unbounded in the direction
$r(m)$. But none of the polyhedra
of $\tilde\foD$ containing $\foj'$ can involve a monomial of the
form $z^{m}$, since a ray carrying a monomial $z^m$ is unbounded
only in the direction $-r(m)$. Thus $0=\bar c_{m,\foj'}=\bar c_{m,\foj}$
as desired.
\qed

\begin{lemma}
\label{twistedderhamlemma}
Let $\theta\in \VV_{\Sigma,k}$, and suppose $f$ is in the ideal
generated by $u_1,\ldots,u_k$ in $\CC[T_{\Sigma}]\otimes_{\CC}
R_k\lfor y_0\rfor$. 
Then for any cycle $\Xi\in H_2((\check\shX_{\Sigma,k})_{\kappa},
\Re(qW_0(Q))\ll 0;\CC)$,
\[
\int_{\Xi} e^{q(W_0(Q)+f)}\Omega=\int_{\Xi} e^{q\theta(W_0(Q)+f)}\Omega.
\]
\end{lemma}

\proof 
We continue to use the fixed identification $\bigwedge^2 M\cong
\ZZ$ given by $\Omega$.
It is enough to show the lemma for $\theta=\exp(cz^{m_0}X_{r(m_0)})$
with $m_0\in T_{\Sigma}$, $r(m_0)\not=0$ and $c^2=0$, as such elements
generate $\VV_{\Sigma,k}$.
Note that if $W_0(Q)+f=\sum_m c_mz^m$, then
\[
\theta(W_0(Q)+f)=\sum_m c_m(z^m+\langle X_{r(m_0)}, r(m)\rangle c z^{m_0+m})
\]
and
\[
e^{q\theta(W_0(Q)+f)}=
e^{q(W_0(Q)+f)}\big(1+\sum_m 
qcc_m \langle X_{r(m_0)}, r(m)\rangle z^{m_0+m}\big).
\]
Furthermore, $d(z^m\dlog(z^{m_0}))=-\langle X_{r(m_0)}, r(m)\rangle
z^m\Omega$.
Thus
\begin{eqnarray*}
(e^{q\theta(W_0(Q)+f)}-e^{q(W_0(Q)+f)})\Omega&=&
e^{q(W_0(Q)+f)}\bigg(q\sum_m cc_m \langle X_{r(m_0)}, r(m)\rangle
z^{m_0+m}\bigg)\Omega\\
&=&-d(cz^{m_0}e^{q(W_0(Q)+f)}\dlog (z^{m_0})).
\end{eqnarray*}
The result then follows from Stoke's theorem and the fact that
$e^{q(W_0(Q)+f)}$ goes to zero rapidly on the unbounded part of $\Xi$.
\qed

\bigskip

\emph{Proof of Lemma \ref{integralindependentlemma}.}
This now follows immediately from Theorems \ref{wallcrossingtheorem},
\ref{Pindeptheorem}, and Lemma \ref{twistedderhamlemma}. \qed

\section{Evaluation of the period integrals}

Our main goal in this section is the computation of the integrals
\[
\int_{\Xi} e^{qW_k(Q)} \Omega
\]
in the case of $\PP^2$. In doing so, we will prove Theorems
\ref{PQinvariant} and \ref{maintheorem}. We continue with
the notation $\Sigma,T_{\Sigma}, t_i,\rho_i$
of Example \ref{firstP2example} as well as $m_i=r(t_i)$.

\begin{lemma}
\label{basicintegral}
There is a (multi-valued) basis $\Xi_0,\Xi_1,\Xi_2$ of sections
of $\shR$ such that,
with $x_0x_1x_2=1$,
\begin{align*}
\sum_{i=0}^2& \alpha^i\int_{\Xi_i} e^{q(x_0+x_1+x_2)}
\Omega=\\
&q^{3\alpha}\left(\sum_{d=0}^{\infty} {q^{3d}\over (d!)^3}
-3\alpha\sum_{d=1}^{\infty} {q^{3d}\over (d!)^3}\sum_{k=1}^d {1\over k}
+{9\over 2}\alpha^2\sum_{d=1}^{\infty}
{q^{3d}\over (d!)^3}
\left(\left(\sum_{k=1}^d{1\over k}\right)^2+{1\over 3}
\sum_{k=1}^d {1\over k^2}\right)\right)\\
=:&q^{3\alpha}\sum_{d=0}^{\infty} q^{3d}\big(B_0(d)+\alpha B_1(d)
+\alpha^2 B_2(d)\big),
\end{align*}
where the last equality defines the numbers $B_0(d),B_1(d),B_2(d)$.
\end{lemma}

\proof This was shown in \cite{Bar}, Proposition 3.1. 
In particular, each integral
$f_i=\int_{\Xi_i} e^{q(x_0+x_1+x_2)}\Omega$ can be shown to 
satisfy the differential equation
\[
(qd/dq)^3f_i=27q^3f_i,
\]
which can then be solved using a recurrence relation to obtain the
above solutions.
\qed

\medskip

We can use this to compute the integrals we are interested in by writing
\[
\int_{\Xi_i} e^{qW_k(Q)}\Omega=
e^{qy_0}\int_{\Xi_i}e^{q(x_0+x_1+x_2)}e^{q(W_k(Q)-W_0(Q))}
\Omega.
\]
The factor $e^{q(W_k(Q)-W_0(Q))}$ can then be expanded in a
Taylor series, noting that in any term, each monomial in $W_k(Q)-
W_0(Q)$ can appear at most once, because it has a coefficient
of square zero; thus this expansion is quite easy and is finite. Thus
we only need to calculate, with $x_0x_1x_2=\kappa$,
\[
\sum_{i=0}^2\alpha^i\int_{\Xi_i} e^{q(x_0+x_1+x_2)}x_0^{n_0}x_1^{n_1}
x_2^{n_2}\Omega.
\]

\begin{lemma} 
\label{unbasicintegral}
With $\Xi_0,\Xi_1,\Xi_2$ as in Lemma 
\ref{basicintegral}, but with $x_0x_1x_2=\kappa$,
\[
\sum_{i=0}^2\alpha^i\int_{\Xi_i} e^{q(x_0+x_1+x_2)}x_0^{n_0}x_1^{n_1}
x_2^{n_2}\Omega
=q^{3\alpha}\kappa^{\alpha}\sum_{i=0}^2\psi_i(n_0,n_1,n_2)\alpha^i,
\]
where
\[
\psi_i(n_0,n_1,n_2)=\sum_{d=0}^{\infty} D_i(d,n_0,n_1,n_2)q^{3d-n_0-n_1-n_2}
\kappa^d
\]
with $D_i$ given as follows. First,
\[
D_0(d,n_0,n_1,n_2)=
\begin{cases}
{1\over (d-n_0)!(d-n_1)!(d-n_2)!}&\hbox{ if $d\ge n_0,n_1,n_2$}\\
0&\hbox{otherwise.}
\end{cases}
\]
Second, if $d\ge n_0,n_1,n_2$, then 
\[
D_1(d,n_0,n_1,n_2)=-{\sum_{k=1}^{d-n_0}{1\over k}+\sum_{k=1}^{d-n_1}{1\over k}
+\sum_{k=1}^{d-n_2}{1\over k}\over (d-n_0)!(d-n_1)!(d-n_2)!}
\]
while if $n_0,n_1\le d <n_2$, then 
\[
D_1(d,n_0,n_1,n_2)={(-1)^{n_2-d-1}(n_2-d-1)!\over
(d-n_0)!(d-n_1)!},
\]
with similar expressions if instead $d<n_0$ or $d<n_1$. If $d$
is smaller than two of $n_0,n_1,n_2$, then 
\[
D_1(d,n_0,n_1,n_2)=0.
\]
Third, if $d\ge n_0,n_1,n_2$, then
\[
D_2(d,n_0,n_1,n_2)={\bigg(\sum_{l=0}^2\sum_{k=1}^{d-n_l} {1\over k}\bigg)^2
+\sum_{l=0}^2\sum_{k=1}^{d-n_l} {1\over k^2}
\over
2 (d-n_0)!(d-n_1)!(d-n_2)!}
\]
while if $n_0,n_1\le d< n_2$,
\[
D_2(d,n_0,n_1,n_2)=
{(-1)^{d-n_2}(n_2-d-1)!\over
(d-n_0)!(d-n_1)!} \bigg( \sum_{k=1}^{d-n_0}{1\over k}
+\sum_{k=1}^{d-n_1} {1\over k}+\sum_{k=1}^{n_2-d-1} {1\over k}\bigg),
\]
with similar expressions if instead $d<n_0$ or $d<n_1$.
If $n_0\le d <n_1,n_2$, then
\[
D_2(d,n_0,n_1,n_2)=
{(-1)^{n_1+n_2}(n_1-d-1)!(n_2-d-1)!\over (d-n_0)!},
\]
with similar expressions if instead $n_1\le d < n_0,n_2$ or $n_2\le d<n_0,n_1$.
Finally, if $d<n_0,n_1,n_2$, then
\[
D_2(d,n_0,n_1,n_2)=0.
\]
\end{lemma}

\proof
Consider the integral 
\[
I_i(a_0,a_1,a_2)=\int_{\Xi_i}e^{a_0x_0+a_1x_1+a_2x_2}\Omega,
\]
with $a_0,a_1,a_2\in \CC^{\times}$ and $x_0x_1x_2=1$. 
Then
\[
{\partial^{n_0+n_1+n_2}
\over \partial a_0^{n_0}\partial a_1^{n_1}\partial a_2^{n_2}}
I_i=\int_{\Xi_i}e^{a_0x_0+a_1x_1+a_2x_2}x_0^{n_0}x_1^{n_1}x_2^{n_2}
\Omega.
\]
Evaluate this at $a_0=a_1=a_2=q\kappa^{1/3}$ and make the
change of variables $x_i\mapsto x_i/\kappa^{1/3}$ in the
integral. Note that as $\Omega={dx_1\wedge dx_2\over x_1x_2}$, such
a change of variables does not affect $\Omega$. Then using
$x_0x_1x_2=\kappa$ we obtain
\[
{\partial^{n_0+n_1+n_2}\over \partial a_0^{n_0}\partial a_1^{n_1}\partial a_2^{n_2}}
I_i\bigg|_{a_i=q\kappa^{1/3}}
=\int_{\Xi_i}e^{q(x_0+x_1+x_2)}\kappa^{-(n_0+n_1+n_2)/3}
x_0^{n_0}x_1^{n_1}x_2^{n_2}\Omega.
\]
On the other hand, $I_i$
can be calculated by making the substitution 
\begin{eqnarray*}
x_0&\mapsto& (a_1a_2/a_0^2)^{1/3}x_0\\
x_1&\mapsto& (a_0a_2/a_1^2)^{1/3}x_1\\
x_2&\mapsto& (a_0a_1/a_2^2)^{1/3}x_2
\end{eqnarray*}
in $I_i$ which gives 
\[
I_i(a_0,a_1,a_2)=
\int_{\Xi_i}e^{(a_0a_1a_2)^{1/3}(x_0+x_1+x_2)}\Omega.
\]
Thus we can compute $\sum_{i=0}^2 \alpha^i I_i(a_0,a_1,a_2)$
by substituting in $q=(a_0a_1a_2)^{1/3}$ in the formula of Lemma
\ref{basicintegral}. To differentiate the resulting expression, note that 
under this substitution, $q^{3\alpha+3d}$ becomes $(a_0a_1a_2)^{\alpha+d}$
and
\begin{align*}
&{\partial^{n_0+n_1+n_2}
\over \partial a_0^{n_0}\partial a_1^{n_1}\partial a_2^{n_2}}
(a_0a_1a_2)^{\alpha+d}\bigg|_{a_i=q\kappa^{1/3}}=\\
&\quad q^{3\alpha+3d-n_0-n_1-n_2}\kappa^{\alpha+d-(n_0+n_1+n_2)/3}
\prod_{k=1}^{n_0}(\alpha+d-k+1)
\prod_{k=1}^{n_1}(\alpha+d-k+1)
\prod_{k=1}^{n_2}(\alpha+d-k+1)\\
&=
q^{3\alpha+3d-n_0-n_1-n_2}\kappa^{\alpha+d-(n_0+n_1+n_2)/3}
(C_0(d,n_0,n_1,n_2)+\alpha C_1(d,n_0,n_1,n_2)+\alpha^2 C_2(d,n_0,n_1,n_2))),
\end{align*}
where the last equality defines $C_0,C_1$ and $C_2$. One then sees that
\[
\psi_i(n_0,n_1,n_2)=\sum_{d=0}^{\infty} \sum_{k=0}^i
B_k(d)C_{i-k}(d,n_0,n_1,n_2)q^{3d-n_0-n_1-n_2}
\kappa^d
\]
with the $B_i$'s defined in Lemma \ref{basicintegral}. Furthermore,
computing the $C_i$'s, we see
\[
C_0(d,n_0,n_1,n_2)=
\begin{cases}
{(d!)^3\over (d-n_0)!(d-n_1)!(d-n_2)!}&\hbox{ if $d\ge n_0,n_1,n_2$}\\
0&\hbox{otherwise.}
\end{cases}
\]
If $d\ge n_0,n_1,n_2$, then
\[
C_1(d,n_0,n_1,n_2)=
{(d!)^3\over (d-n_0)!(d-n_1)!(d-n_2)!}\bigg(
\sum_{k=d-n_0+1}^d {1\over k}
+\sum_{k=d-n_1+1}^d {1\over k}
+\sum_{k=d-n_2+1}^d {1\over k}\bigg),
\]
while if \emph{one} of $n_0,n_1,n_2$ is larger
than $d$, we have
\[
C_1(d,n_0,n_1,n_2)=
\prod_{k=d-n_0+1\atop k\not=0}^d k
\prod_{k=d-n_1+1\atop k\not=0}^d k
\prod_{k=d-n_2+1\atop k\not=0}^d k.
\]
Otherwise
\[
C_1(d,n_0,n_1,n_2)=0.
\]
If $d\ge n_0,n_1,n_2$ then 
\begin{align*}
C_2&(d,n_0,n_1,n_2)=\\
&{(d!)^3\over 2(d-n_0)!(d-n_1)!(d-n_2)!}\bigg(
\bigg(\sum_{k=d-n_0+1}^d {1\over k}
+\sum_{k=d-n_1+1}^d{1\over k}
+\sum_{k=d-n_2+1}^d {1\over k}\bigg)^2\\
&\quad\quad-
\bigg(\sum_{k=d-n_0+1}^d {1\over k^2}
+\sum_{k=d-n_1+1}^d{1\over k^2}
+\sum_{k=d-n_2+1}^d {1\over k^2}\bigg)\bigg).
\end{align*}
If $n_1,n_2\le d< n_0$, then
\begin{align*}
C_2&(d,n_0,n_1,n_2)=\\
&\bigg(\prod_{k=d-n_0+1\atop k\not=0}^d k\bigg) {(d!)^2\over (d-n_1)!(d-n_2)!}
\bigg(
\sum_{k=d-n_0+1\atop k\not=0}^d {1\over k}
+\sum_{k=d-n_1+1}^d {1\over k}+\sum_{k=d-n_2+1}^d {1\over k}
\bigg).
\end{align*}
We have similar expressions if $d<n_1$ or $d<n_2$. 
If two of $n_0,n_1$ and $n_2$ are larger than $d$,
then
\[
C_2(d,n_0,n_1,n_2)=\prod_{k=d-n_0+1\atop k\not=0}^d k
\prod_{k=d-n_1+1\atop k\not=0}^d k
\prod_{k=d-n_2+1\atop k\not=0}^d k.
\]
Finally, if $n_0,n_1,n_2>d$, then
\[
C_2(d,n_0,n_1,n_2)=0.
\]
A laborious calculation now gives the forms given in the Lemma for the 
coefficients $D_i$.
\qed

\begin{remark}
The coefficient $D_i(d,n_0,n_1,n_2)$ can be written in a more uniform
way, which is convenient for computation on a computer algebra system.
Obviously
\[
D_0(d,n_0,n_1,n_2)={1\over \Gamma(d-n_0+1)\Gamma(d-n_1+1)\Gamma(d-n_2+1)}.
\]
For $i=1$, we obtain
\[
D_1(d,n_0,n_1,n_2)=-{3\gamma+\Psi_0(d-n_0+1)+\Psi_0(d-n_1+1)+\Psi_0(d-n_2+1)
\over \Gamma(d-n_0+1)\Gamma(d-n_1+1)\Gamma(d-n_2+1)}.
\]
Here $\gamma$ is Euler's constant and 
$\Psi_0$ is the digamma function, 
$\Psi_0=\Gamma'/\Gamma$, with $\Psi_0(d)=-\gamma+\sum_{k=1}^{d-1} 1/k$
for $d\ge 1$. If $\max(n_0,n_1,n_2)>d$, this is interpreted as a limit,
using $\lim_{x\rightarrow -d} \Psi_0(x)/\Gamma(x)=(-1)^{d+1}d!$, for $d>0$
an integer. 

Finally, $D_2(d,n_0,n_1,n_2)$ is
\[
{\pi^2+18\gamma^2+12\gamma\bigg(\sum_{k=0}^2 \Psi_0(d-n_k+1)\bigg)
+2\bigg(\sum_{k=0}^2 \Psi_0(d-n_k+1)\bigg)^2-2
\sum_{k=0}^2 \Psi_1(d-n_k+1)
\over 4\Gamma(d-n_0+1)\Gamma(d-n_1+1)\Gamma(d-n_2+1)}.
\]
Here $\Psi_1=\Psi_0'$, and $\Psi_1(d)=\pi^2/6-\sum_{k=1}^{d-1} 1/k^2$
for $d>0$ an integer. Again, with the appropriate limit interpretation,
this covers all cases.
\end{remark}

\begin{definition}
For $m\in T_{\Sigma}$, $m=\sum_{i=0}^2 n_it_i$ with $n_i\ge 0$ for all $i$,
define
\begin{eqnarray*}
\psi_i(m)&:=&\psi_i(n_0,n_1,n_2)\\
D_i(d,m)&:=&D_i(d,n_0,n_1,n_2)
\end{eqnarray*}
and
\[
|m|:=n_0+n_1+n_2.
\]
\end{definition}

We will now start on our proof of Theorems 
\ref{PQinvariant} and \ref{maintheorem}, beginning
with the following definition.

\begin{definition} Fix $P_1,\ldots,P_k$ general. For $Q$
general, let $S_k$ (or $S_k(Q)$ if the dependence on $Q$
needs to be emphasized) be a finite set of triples $(c,\nu,m)$
with $c\in R_k$ a monomial such that 
\begin{equation}
\label{expexpansion}
e^{q(W_k(Q)-W_0(Q))}
=\sum_{(c,\nu,m)\in S_k} cq^{\nu}z^m,
\end{equation}
with each term $cq^{\nu}z^m$ of the form $\prod_{i=1}^{\nu}
\Mono(h_i)$ for $h_1,\ldots,h_{\nu}$ distinct Maslov index two
tropical disks with boundary $Q$.

Let
\[
L_i^d=L_i^d(Q):=
\sum_{(c,\nu,m)\in S_k} cq^{3d+\nu-|m|} D_i(d,m).
\]
\end{definition}

We can now clarify what needs to be proved. The following
lemma reduces the two theorems to three equalities.

\begin{lemma}
\label{threeforms}
Let $Q$ be chosen generally, and let $L$ be the tropical
line with vertex $Q$. The three equalities
\begin{eqnarray}
\label{eq56-}
L^d_0&=&\delta_{0,d}+\sum_{\nu\ge 0}
\sum_{\scriptstyle I\subseteq\{1,\ldots,k\}\atop 
{\scriptstyle
I=\{i_1,\ldots,i_{3d-\nu-2}\}\atop \scriptstyle
i_1<\cdots<i_{3d-\nu-2}}}
\langle P_{i_1},\ldots,P_{i_{3d-\nu-2}},\psi^{\nu}L\rangle_d^{\trop}
u_Iq^{\nu+2}\\
\label{eq56}
L^d_1&=&\sum_{\nu\ge 0}
\sum_{\scriptstyle I\subseteq\{1,\ldots,k\}\atop 
{\scriptstyle
I=\{i_1,\ldots,i_{3d-\nu-2}\}\atop \scriptstyle
i_1<\cdots<i_{3d-\nu-1}}}
\langle P_{i_1},\ldots,P_{i_{3d-\nu-1}},\psi^{\nu}L\rangle_d^{\trop}
u_Iq^{\nu+1}\\
\label{eq57}
L^d_2&=&y_2q^{-1}\delta_{0,d}+
\sum_{\nu\ge 0}\sum_{\scriptstyle I\subseteq\{1,\ldots,k\}\atop
{\scriptstyle
I=\{i_1,\ldots,i_{3d-\nu}\}\atop \scriptstyle
i_1<\cdots<i_{3d-\nu}}}
\langle P_{i_1},\ldots,P_{i_{3d-\nu}},\psi^{\nu}M_{\RR}\rangle_d^{\trop}
u_Iq^{\nu}.
\end{eqnarray}
imply Theorems
\ref{PQinvariant} and \ref{maintheorem}.
\end{lemma}

\proof
Let us be precise about what needs to be shown to prove Theorems
\ref{PQinvariant} and \ref{maintheorem}.
If we write, for $0\le i
\le 2$,
\[
K_i^{\trop}=\sum_{d\ge 1}\sum_{\nu\ge 0} \langle T_2^{3d+i-2-\nu},
\psi^{\nu}T_{2-i}\rangle^{\trop}_d q^{\nu+2}\kappa^d{y_2^{3d+i-2-\nu}
\over (3d+i-2-\nu)!},
\]
then 
\begin{eqnarray*}
J_0^{\trop}&=&e^{qy_0}(1+K_0^{\trop})\\
J_1^{\trop}&=&e^{qy_0}(qy_1(1+K_0^{\trop})+K_1^{\trop})\\
J_2^{\trop}&=&e^{qy_0}\big({q^2y_1^2\over 2}(1+K_0^{\trop})
+qy_1K_1^{\trop}+qy_2+K_2^{\trop}\big)
\end{eqnarray*}
We wish to compare these expressions with the expressions
obtained via period integrals over $\Xi_0$, $\Xi_1$ and $\Xi_2$.
Take for the $\Xi_i$ the cycles given by Lemma \ref{basicintegral}.
Consider the $\varphi_i$'s defined using these cycles
in \eqref{phidefeq}. 
Expanding the integral in \eqref{phidefeq}
by using Lemma \ref{unbasicintegral}
and $\kappa^{\alpha}=e^{y_1\alpha}=1+y_1\alpha+y_1^2\alpha^2/2$,
the left-hand side of \eqref{phidefeq} is
\[
\sum_{(c,\nu,m)\in S_k}
ce^{qy_0}q^{3\alpha+\nu}\kappa^{\alpha}
\sum_{i=0}^2 \psi_i(m)\alpha^i\\
=
q^{3\alpha}e^{qy_0}\sum_{(c,\nu,m)\in S_k} 
cq^{\nu}\sum_{i=0}^2\sum_{k=0}^i {y_1^k\over k!}\psi_{i-k}(m)\alpha^i.
\]
Comparing this with the right-hand side of \eqref{phidefeq},
we get
\begin{eqnarray*}
\varphi_0&=&e^{qy_0}\sum_{(c,\nu,m)\in S_k} cq^{\nu}
\psi_0(m)\\
\varphi_1&=&e^{qy_0}\sum_{(c,\nu,m)\in S_k} cq^{\nu+1}
(y_1\psi_0(m)+\psi_1(m))\\
\varphi_2&=&e^{qy_0}\sum_{(c,\nu,m)\in S_k} cq^{\nu+2}
\big({y_1^2\over 2}\psi_0(m)+y_1\psi_1(m)
+\psi_2(m)\big).
\end{eqnarray*}
Thus to show $\varphi_i=J_i^{\trop}$, we need to show the following
three equalities:
\begin{eqnarray}
\sum_{(c,\nu,m)\in S_k} cq^{\nu}\psi_0(m)
&=&1+K_0^{\trop}
\label{eq51}\\
\sum_{(c,\nu,m)\in S_k} cq^{\nu}\psi_1(m)
&=&q^{-1}K_1^{\trop}
\label{eq52}\\
\sum_{(c,\nu,m)\in S_k} cq^{\nu}\psi_2(m)
&=&q^{-2}(qy_2+K_2^{\trop}).
\label{eq53}
\end{eqnarray}
Then using the expansion for
$\psi_i$ in Lemma \ref{unbasicintegral}, \eqref{eq51},
\eqref{eq52} and \eqref{eq53}
are equivalent, if we compare the coefficients of $\kappa^d$
on both sides, to:
\begin{eqnarray}
\label{eq54-}
L^d_0&=&\delta_{0,d}+
\sum_{\nu\ge 0}
\langle T_2^{3d-\nu-2},\psi^{\nu}T_2\rangle_d^{\trop}
{y_2^{3d-\nu-2}\over (3d-\nu-2)!}q^{\nu+2},\\
\label{eq54}
L^d_1&=&
\sum_{\nu\ge 0}
\langle T_2^{3d-\nu-1},\psi^{\nu}T_1\rangle_d^{\trop}
{y_2^{3d-\nu-1}\over (3d-\nu-1)!}q^{\nu+1},\\
\label{eq55}
L^d_2&=&y_2q^{-1}\delta_{0,d}+
\sum_{\nu\ge 0}
\langle T_2^{3d-\nu},\psi^{\nu}T_0\rangle_d^{\trop}
{y_2^{3d-\nu}\over (3d-\nu')!}q^{\nu}.
\end{eqnarray}
Now suppose we have shown \eqref{eq56-}, \eqref{eq56} and
\eqref{eq57}. The left-hand sides of these equations
come from the period integrals, and hence are independent
of the locations of $Q$ and $P_1,\ldots,P_k$ by Lemma 
\ref{integralindependentlemma}. So the right-hand side is also
independent of the locations of $Q$ and $P_1,\ldots,P_k$.
So in particular, once we show \eqref{eq56-}, \eqref{eq56} and \eqref{eq57}, 
we find that the invariants $\langle T_2^{3d+i-2-\nu},\psi^{\nu}
T_{2-i}\rangle_d^{\trop}$ are well-defined, showing Theorem 
\ref{PQinvariant}, and also showing \eqref{eq54-},
\eqref{eq54} and \eqref{eq55},
hence $\varphi_i=J_i^{\trop}$. In particular, $\varphi_{i,1}=y_i$
for $0\le i\le 2$. This gives Theorem \ref{maintheorem}. \qed

\medskip

We have in fact already taken care of \eqref{eq56-}:

\begin{lemma}
\label{zerolemma}
\eqref{eq56-} holds.
\end{lemma}

\proof In fact Proposition \ref{J0case} shows the equivalent statement
that $\varphi_0=J_0^{\trop}$. Note however that the proof of
Proposition \ref{J0case} was carried out using a specific choice
of $\Xi_0$, which a priori may not be the same $\Xi_0$ given by
Lemma \ref{basicintegral}. However, one checks easily that the
values for the integrals $\int_{\Xi_0}e^{q(x_0+x_1+x_2)}x_0^{n_0}
x_1^{n_1}x_2^{n_2}\Omega$ given in Lemma \ref{unbasicintegral}
agree with the integrals over the $\Xi_0$ used in the proof of
Proposition \ref{J0case}. Thus the argument still works.
\qed

\medskip

We will now refine the expressions $L^d_i$ which we need
to compute.

\begin{definition}
For each cone $\sigma\in\Sigma$, $\sigma$ is the image
under $r$ of a proper face $\tilde\sigma$ of the cone 
$K\subseteq T_{\Sigma}\otimes\RR$ generated by $t_0,t_1,t_2$
(i.e., the first octant). For $d\ge 0$, denote by $K_d
\subseteq K$ the cube
\[
K_d=\big\{\sum_{i=0}^2 n_i t_i\,\big|\, 0\le n_i \le d\big\}
\]
and for $\sigma\in\Sigma$, define
\[
\tilde\sigma_d:=
(\tilde\sigma+K_d)\setminus \bigcup_{\scriptstyle \tau\subsetneq \sigma
\atop \scriptstyle \tau\in\Sigma} (\tilde\tau+K_d).
\]
Here $+$ denotes Minkowski sum. 
\end{definition}

\begin{example}
We have the following examples of $\tilde\sigma_d$.
Let $m=\sum_i n_it_i \in K$.
\begin{itemize}
\item
If $\sigma=\{0\}$, then $m\in \tilde\sigma_d$ if and only if
$d\ge \max\{n_0,n_1,n_2\}$.
\item If $\sigma=\rho_0$, then
$m\in \tilde\sigma_d$ if and only if $n_1,n_2\le d < n_0$. 
\item If $\sigma=\rho_1+\rho_2$, then $m\in\tilde\sigma_d$ if and only
if $n_0\le d < n_1,n_2$.
\end{itemize}
\end{example}

\begin{definition}
For $\sigma\in\Sigma$, define
\[
L^d_{i,\sigma}=L^d_{i,\sigma}(Q)
:=\sum_{\scriptstyle (c,\nu,m)\in S_k\atop \scriptstyle 
m\in\tilde\sigma_d}
cq^{3d+\nu-|m|}D_i(d,m).
\]
\end{definition}

\begin{lemma}
\label{stupidLdsigmafacts}
\begin{enumerate}
\item
$L_i^d=\sum_{\sigma\in\Sigma} L^d_{i,\sigma}$.
\item
\[
L^d_{i,\{0\}}(Q)=\sum_{\nu\ge i}
\sum_{\scriptstyle I\subseteq\{1,\ldots,k\}\atop
{\scriptstyle I=\{i_1,\ldots,i_{3d-2+i-\nu}\}\atop\scriptstyle 
i_1<\cdots<i_{3d-2+i-\nu}}}
\langle P_{i_1},\ldots,P_{i_{3d-2+i-\nu}},\psi^{\nu}S
\rangle_{d,\{0\}}^{\trop} u_Iq^{\nu+2-i}
\]
where $S=Q,L$ the tropical line with vertex $Q$,
or $M_{\RR}$ in the cases $i=0,1$ and $2$.
Here, the meaning of the notation on the right-hand side with
subscript $\{0\}\in\Sigma$ is defined in \eqref{tropdesczonedef}.
\end{enumerate}
\end{lemma}

\proof
(1) just follows from Lemma \ref{unbasicintegral}, which tells
us that $D_i(d,m)=0$ if $m\not\in\bigcup_{\sigma\in\Sigma}
\tilde\sigma_d$.

(2) This is essentially the same argument as made in the proof of Proposition
\ref{J0case}. Let $(c,\nu,m)\in S_k$ with $m=\sum_{i=0}^2 n_it_i$. Then
$(c,\nu,m)$ contributes to $L^d_{i,\{0\}}$ only if $n_0,n_1,n_2\le d$.
Write 
\[
cq^{\nu}z^m=q^{\nu}\prod_{i=1}^{\nu}\Mult(h_i)z^{\Delta(h_i)}u_{I(h_i)}
\]
for $h_i:\Gamma_i'\rightarrow M_{\RR}$, $1\le i\le\nu$, 
Maslov index two disks with boundary $Q$. Let $\Gamma$ be the graph
obtained by identifying the outgoing vertices $V_{\out,i}$
of $\Gamma_1',\ldots,
\Gamma_{\nu}'$ to get a single vertex $V_{\out}$
and then adding $(d-n_0)+(d-n_1)+(d-n_2)+1$ additional unbounded edges with
vertex $V_{\out}$. We define $h:\Gamma\rightarrow M_{\RR}$ to be
$h_i$ on each subgraph $\Gamma_i'\subseteq\Gamma$. Furthermore, for $0\le i
\le 2$, $h$ maps $d-n_i$ of the new unbounded edges to the ray
$Q+\RR_{\ge 0}m_i$. Finally, the last unbounded ray is labelled with
an $x$ and is contracted by $h$. Just as in the argument of Proposition
\ref{J0case}, $h$ is now a balanced tropical curve. 

The contribution of this term $cq^{\nu}z^m$ to $L^d_{i,\{0\}}$
is then 
\begin{eqnarray*}
&&q^{3d+\nu-n_0-n_1-n_2}D_i(d,n_0,n_1,n_2)\prod_{i=1}^{\nu}\Mult(h_i)u_{I(h_i)}
\\
&=&q^{3d+\nu-n_0-n_1-n_2}u_{I(h)}\Mult^i_x(h)\prod_{V\in\Gamma^{[0]}
\atop V\not\in E_x}\Mult_V(h),
\end{eqnarray*}
comparing the definitions of $\Mult^i_x(h)$ and $D_i(d,n_0,n_1,n_2)$.
Note that the valency $\Val(V_{\out})$ of the vertex $V_{\out}$ in $h$ is
$\nu+3d-(n_0+n_1+n_2)+1$. 
Suppose that $I(h)=\{i_1,\ldots,i_{3d-2+i-\nu'}\}$ for some $\nu'$.
Noting that $h$ is obtained by gluing $\Val(V_{\out})-1$ Maslov
index two disks, we see that
\begin{eqnarray*}
\Val(V_{\out})-1 &=& \sum_{i=1}^{\nu} \big(|\Delta(h_i)|-
\#I(h_i)\big)+(d-n_1)+(d-n_2)+(d-n_3)\\
&=&3d-(3d-2+i-\nu')=\nu'+2-i.
\end{eqnarray*}
Then
the curve $h$ contributes precisely the correct
contribution, as given by Definition \ref{tropicaldescinv}, (1) (a),
(2) (b), or (3) (d), to
\[
\langle P_{i_1},\ldots,P_{i_{3d-2+i-\nu'}},\psi^{\nu'}S\rangle^{\trop}_{d,\{0\}}
u_{I(h)}q^{\nu'+2-i}.
\]

Conversely, given any curve $h$ contributing to the above quantity,
it follows from Lemma \ref{curvechop}, (3), that $h$ will
arise in the above manner from some term $(c,\nu,m)\in S_k$.
\qed

\medskip

Next, we need to understand the asymptotic behaviour of
$L^d_{i,\sigma}(Q)$.

\begin{lemma}
\label{asympbehav}
Let $\omega\in\Sigma$, and
let $v\in\omega$
be non-zero (hence ruling out $\omega=\{0\}$).
Then
\begin{equation}
\label{omegavanish}
\lim_{s\rightarrow\infty}L_{i,\omega}^d(Q+sv)
=0.
\end{equation}
\end{lemma}

\proof 
We first note that with $\omega\not=\{0\}$,
\begin{equation}
\label{omegatildefact}
\hbox{if $m\in\tilde\omega_d$, then 
$r(m)\in \bigcup_{\sigma\supseteq\omega\atop\sigma\in\Sigma} \Int(\sigma)$.}
\end{equation}

Next, for sufficiently large $s$, $Q+sv$ lies in an unbounded
connected component $\shC$ of $M_{\RR}\setminus \Supp(\foD)$,
where $\foD=\foD(\Sigma,P_1,
\ldots,P_k)$. By taking $s$ sufficiently large, we can assume $\shC$ is
the last
component entered as $s\rightarrow\infty$. To show \eqref{omegavanish},
it will be enough to show that if $Q+sv\in \shC$, there exists a 
convex cone $K'\subseteq M_{\RR}$ with $K'\cap \bigcup_{\sigma\supseteq
\omega\atop \sigma\in\Sigma}\Int(\sigma)=\emptyset$
such that $W_k(Q+sv)-W_0(Q+sv)$ 
only contains monomials $z^m$ with $r(m)\in K'$. It then follows
that all monomials $z^m$ in $\exp(q(W_k(Q+sv)-W_0(Q+sv)))$ 
satisfy $r(m)\in K'$, and hence by \eqref{omegatildefact},
$m\not\in \tilde\omega_d$.
This implies \eqref{omegavanish}. 

So we study monomials $z^m$ appearing in $W_k(Q+sv)-W_0(Q+sv)$ 
and construct a cone $K'$ with the desired properties.
We will make use of
the asymptotic cone to the closure $\overline{\shC}$ of $\shC$, 
$\Asym(\overline{\shC})$, which is defined to
be the Hausdorff limit $\lim_{\epsilon\rightarrow 0} 
\epsilon\overline{\shC}$.
Note that the connected components of $M_{\RR}\setminus
\foD(\Sigma,P_1)$ are $P_1-\Int(\sigma)$ where $\sigma$ runs over the maximal
cones of $\Sigma$. Since $\Supp(\foD(\Sigma,P_1))\subseteq
\Supp\foD$, one sees that $\Asym(\overline{\shC})$ is contained in some cone
$-\sigma$ with $\sigma\in\Sigma$ maximal and $(-\sigma)\cap\omega\not=
\{0\}$. Note also that $\Asym(\overline{\shC})$ can be a ray
if the unbounded edges of $\overline{\shC}$ are parallel.
Let $\fod_1$, $\fod_2$ denote the two unbounded edges of $\overline{\shC}$.

Now for general $s$, a term $cz^m$ in $W_k(Q+sv)$ corresponds to a broken line
$\beta$ with given data $-\infty=t_0<\cdots <t_p=0$,
$m_i^{\beta}\in T_{\Sigma}$ as in Definition 
\ref{brokenlinedef}, and $m=m_{p}^{\beta}$. 
If $-r(m)\not\in \RR_{>0}v$, then
for $s$ sufficiently large, with $+$ denoting Minkowski sum,
\[
Q+sv\not\in \RR_{\ge 0}(-r(m))+(\partial\overline{\shC}\setminus
(\fod_1\cup\fod_2)).
\]
Indeed, $\partial\overline{\shC}\setminus (\fod_1\cup\fod_2)$
is bounded, so the asymptotic cone of the right-hand side is
$\RR_{\ge 0}(-r(m))$, which does not contain $v$ by assumption.
Thus, taking a sufficiently large $s$, we note $\beta$ cannot last
enter $\shC$ via $\partial\overline{\shC}\setminus (\fod_1\cup
\fod_2)$ since the last line segment of $\beta$ is in the direction
$-r(m)$.
So for sufficiently large $s$,
$\beta$ must enter $\overline{\shC}$ by crossing
one of $\fod_1$ or $\fod_2$.
In what follows, we will not need to
study the case $-r(m)\in \RR_{>0} v$ as the cone $K'$ we construct
will always contain $-v$.

We can now assume that for large $s$, $\beta$
enters $\shC=\shC_n$ from another unbounded connected
component $\shC_{n-1}$ of $M_{\RR}\setminus \Supp(\foD)$.
Necessarily, the $m_i^{\beta}$ attached to $\beta$
while $\beta$ passes through $\shC_{n-1}$ satisfies $-r(m_i^{\beta})
\not\in\Asym(\overline{\shC_{n-1}})$. Indeed,
otherwise $\beta$ could not hit an unbounded
edge of $\overline{\shC_{n-1}}$. Again, for large enough $s$, one sees similarly
that $\beta$ must enter $\shC_{n-1}$ through the other unbounded
edge of $\overline{\shC_{n-1}}$, and we can then continue this process
inductively, with $\beta$ passing only through unbounded edges
via a sequence of unbounded components $\shC_0,\ldots,\shC_n$.
When $\beta$ bends, it then always bends outward, as depicted in 
Figure \ref{bendingoutward}. From this we make the following two
observations:
\begin{itemize}
\item[(C1)] If the edges corresponding to $\fod_1$ and $\fod_2$
of $\Asym(\overline{\shC})$ are generated by $v_1$, $v_2$ respectively (possibly
$v_1=v_2$) and
$\beta$ enters $\shC$ by crossing $\fod_i$, then
$-r(m)$ lies in a half-plane with boundary $\RR v_i$
containing $\Asym(\overline{\shC})$; otherwise, $\beta$ cannot reach the
interior of $\overline{\shC}$.
\item[(C2)] For any $j$, $1\le j\le p$,
$-r(m)$ lies in the half-plane with boundary
$\RR r(m_j^{\beta})$ containing $v_i$ corresponding to the edge $\fod_i$
that $\beta$ crosses to enter $\shC$. This follows
from the behaviour described above about how $\beta$ bends.
\end{itemize}
\begin{figure}
\input{bendingoutward.pstex_t}
\caption{}
\label{bendingoutward}
\end{figure}

Without loss of generality, let us assume for the ease
of drawing pictures that $\omega=\rho_2$ or $\rho_1+\rho_2$
and $\Asym(\overline{\shC})\subseteq -(\rho_0+\rho_1)$. See Figure
\ref{conelocations}. Note that as depicted there, we must have
$v_2\in\rho_1+\rho_2$.

\begin{figure}
\input{conelocations.pstex_t}
\caption{}
\label{conelocations}
\end{figure}

We analyze the possibilities for $\beta$:
we have three cases, based on whether the initial direction
of $\beta$ is $-m_1$, $-m_2$, or $-m_0$.

\emph{Case 1.}
$r(m^{\beta}_0)=m_1$. Then $\beta$ must enter $\shC$ via
$\fod_2$. By (C1), $-r(m)$ lies in the half-plane with
boundary $\RR v_2$ containing $\Asym(\overline\shC)$, and by
(C2), $-r(m)$ lies in the half-plane with boundary $\RR m_1$
containing $\Asym(\overline\shC)$. Thus $-r(m)\in (-\RR_{\ge 0} m_1
+\RR_{\ge 0}v_2)$.

\emph{Case 2}. $r(m^{\beta}_0)=m_2$. Then either $\RR_{\ge 0}m_2
\subseteq \Asym(\overline{\shC})$ or 
$\Asym(\overline{\shC})\subseteq \rho_1+\rho_2$ since
$v\in\Asym(\overline{\shC})$. In the first
case, $\beta$ has no opportunity to bend, so corresponds to
the monomial $x_2$, which doesn't appear in $W_k(Q+sv)-W_0(Q+sv)$.
In the second case, $\beta$ bends at time $t_1$ as it crosses
a ray $\fod\in\foD$ with
$f_{\fod}=1+c_{\fod}z^{m_{\fod}}$ with $-r(m_{\fod})\in
\Int(\rho_1+\rho_2)$. Now $r(m_1^{\beta})=m_2+r(m_{\fod})$,
so it follows that $-r(m_1^{\beta})\in \rho_1+\rho_2$. (Here we
use integrality of $m_\fod$ and $m_2=(0,1)$.) Thus by (C1) and (C2),
$-r(m) \in (\RR_{\ge 0} m_1+\RR_{\ge 0} v_1)$.

\emph{Case 3}. $r(m^{\beta}_0)=m_0$. In this case $\beta$ must enter
$\shC$ through the edge $\fod_1$ since $\Asym(\overline{\shC})
\subseteq -(\rho_0+\rho_1)$.
Then one sees from
(C1) and (C2) that
$-r(m)\in (\RR_{\ge 0}(-m_0)+\RR_{\ge 0} v_1)$.

\medskip

We now see that if $\RR_{\ge 0}m_2\subseteq\Asym(\shC)$,
(which always happens if $v$ is proportional to $m_2$,
in particular when $\omega=\rho_2$), then
of these three cases, only cases 1 and 3 can occur,
and in fact $r(m),-v
\in \rho_0+\rho_1$.
Thus $K'=\rho_0+\rho_1$ is the desired cone,
proving the claim in this case. 

If $\RR_{\ge 0}m_2 \not\subseteq \Asym(\shC)$, then $v$ is not proportional
to $m_2$ and $\omega=\rho_1+\rho_2$. In this case, the above
three cases show that $-r(m)$ is always contained in the
upper half-plane. Thus $K'$ the lower half-plane is the desired
cone, proving the claim in this case.
\qed

\medskip

The next step is to explain how $L^d_{i,\sigma}(Q)$ depends on
$Q$ via a wall-crossing formula. While of course $L^d_i$ is independent
of $Q$, the way the terms in $L^d_i$ are redistributed among
the expressions $L^d_{i,\sigma}(Q)$ is key to the calculations.

\begin{definition}
\label{Ldarrowdef}
Let $\foD=\foD(\Sigma,P_1,\ldots,P_k)$.
Let $\shC_1$, $\shC_2$ be two connected components
of $M_{\RR}\setminus\Supp(\foD)$ with $\dim\overline{\shC}_1\cap
\overline{\shC}_2=1$. Let $Q_i\in\shC_i$ be general points, and let
$\gamma$ be a path from $Q_1$ to $Q_2$, passing through $\Supp(\foD)$
only at one time $t_0$, with $\gamma(t_0)\not\in
\Sing(\foD)$. Let $\fod\in\foD$ be a ray with $\gamma(t_0)\in\fod$,
and let $n_{\fod}\in N$ be a primitive vector which is orthogonal
to $\fod$ and satisfies $\langle n_{\fod},\gamma'(t_0)\rangle <0$.
Writing $f_{\fod}=1+c_{\fod}z^{m_{\fod}}$, note that
\[
\theta_{\gamma,\fod}(z^m)=z^m+c_{\fod}\langle n_{\fod},r(m)\rangle
z^{m+m_{\fod}}.
\]
Now take a pair $\omega\subsetneq\tau$ with $\omega,\tau\in\Sigma$
and $\dim\tau=\dim\omega+1$. Note there is a unique index $j\in
\{0,1,2\}$ such that $m_j\not\in\omega$ but $m_j\in\tau$;
call this index $j(\omega,\tau)$.
Then define
\[
L^d_{i,\fod,\gamma,\omega\rightarrow\tau}
:=
\sum_{(c,\nu,m)} \langle n_{\fod},m_{j(\omega,\tau)}\rangle c_{\fod}c
D_i(d,m+m_{\fod}+t_{j(\omega,\tau)})q^{\nu+3d-|m+m_{\fod}|},
\]
where the sum is over all $(c,\nu,m)\in S_k(Q_1)$ such that $m+m_{\fod}
\in \tilde\omega_d$ but $m+m_{\fod}+t_{j(\omega,\tau)}\in \tilde\tau_d$.
If $(c,\nu,m)\in S_k(Q_1)$ satisfies this condition, 
then we say \emph{the term $cq^{\nu}z^m$
contributes to $L^d_{i,\fod,\gamma,\omega\rightarrow\tau}$}.

Define
\[
L^d_{i,\gamma,\omega\rightarrow\tau}
:=\sum_{\fod} L^d_{i,\fod,\gamma,\omega\rightarrow\tau},
\]
where the sum is over all $\fod\in\foD$ with $\gamma(t_0)\in\fod$.

For an arbitrary path $\gamma$ in $M_{\RR}\setminus \Sing(\foD)$
with $\gamma(0)=Q$, $\gamma(1)=Q'$, choose a partition of
$[0,1]$, $0=t_0< t_1<\cdots t_n=1$ such that $\gamma|_{[t_{j-1},t_j]}$
is a path of the sort considered above, connecting endpoints
in adjacent connected components. Then define
\[
L^d_{i,\gamma,\omega\rightarrow\tau}
:=\sum_{j=1}^n L^d_{i,\gamma|_{[t_{j-1},t_j]},\omega\rightarrow\tau}.
\]
\qed
\end{definition}

\begin{lemma} Let $P_1,\ldots,P_k$ be general.
Let $\gamma$ be a path in $M_{\RR}\setminus\Sing(\foD)$ with
$\gamma(0)=Q$, $\gamma(1)=Q'$.
Then for $\dim\rho=1$, $\rho\in\Sigma$,
\begin{equation}
\label{eqrhochange}
L^d_{i,\rho}(Q')-L^d_{i,\rho}(Q)
=L^d_{i,\gamma,\{0\}\rightarrow\rho}
-\sum_{\scriptstyle \sigma\in\Sigma\atop
\scriptstyle \rho\subsetneq\sigma}
L^d_{i,\gamma,\rho\rightarrow\sigma}
\end{equation}
while for $\dim\sigma=2$, $\sigma\in\Sigma$,
\begin{equation}
\label{eqsigmachange}
L^d_{i,\sigma}(Q')-L^d_{i,\sigma}(Q)
=
\sum_{\scriptstyle \rho\in\Sigma\atop{\scriptstyle \dim\rho=1\atop\scriptstyle 
\rho\subsetneq\sigma}}
L^d_{i,\gamma,\rho\rightarrow\sigma}.
\end{equation}
\end{lemma}

\proof It is enough to show this for $\gamma$ a short path
connecting $Q$ and $Q'$ in two adjacent components $\shC_1$ and
$\shC_2$ of $M_{\RR}\setminus\Supp(\foD)$ as in Definition
\ref{Ldarrowdef}. Suppose that at time $t_0$,
$\gamma(t_0)\in \fod_1\cap\cdots \cap\fod_s$ for
rays $\fod_1,\ldots,\fod_s\in\foD$. Of course, $\dim\fod_i\cap\fod_j=1$.
We can then write, for $n_{\fod}=n_{\fod_i}$ for any $i$,
\begin{align*}
\theta_{\gamma,\foD}(z^m)={} &
z^m\prod_{i=1}^s f_{\fod_i}^{\langle n_{\fod},r(m)\rangle} \\
= {} & z^m\prod_{i=1}^s(1+c_{\fod_i}\langle n_{\fod},r(m)\rangle
z^{m_{\fod_i}})\\
={}&
z^m+\sum_{i=1}^s c_{\fod_i}\langle n_{\fod},r(m)\rangle
z^{m+m_{\fod_i}}.
\end{align*}
Here the last equality follows from $c_{\fod_i}c_{\fod_j}=0$
for $i\not=j$. This is the case by the assumption that
$P_1,\ldots,P_k$ are general. Indeed, if $c_{\fod_i}c_{\fod_j}\not=0$,
then the Maslov index zero trees $h_i$ and $h_j$
corresponding to $\fod_i$ and $\fod_j$ would have $I(h_i)\cap I(h_j)
=\emptyset$. However, a generic perturbation of the marked points with
indices in $I(h_i)$ would deform $\fod_i$ without deforming
$\fod_j$, so that $\fod_i\cap\fod_j=\emptyset$. 

Now
\[
W_k(Q')=\theta_{\gamma,\foD}(W_k(Q))
\]
by Theorem \ref{wallcrossingtheorem}. 
Using the expansion \eqref{expexpansion} and $W_0(Q)=y_0+\sum_{j=0}^2
z^{t_j}$,
\begin{eqnarray*}
&&\exp\big(q(W_k(Q')-W_0(Q'))\big)\\
&=&
\exp\big(q(\theta_{\gamma,\foD}(W_k(Q))-W_0(Q))\big)\\
&=&
\theta_{\gamma,\foD}\big(\exp(q(W_k(Q)-W_0(Q)))\big)
\cdot \exp\big(q(\theta_{\gamma,\foD}(W_0(Q))
-W_0(Q))\big)\\
&=&\theta_{\gamma,\foD}\bigg(\sum_{(c,m,\nu)\in S_k(Q)} cq^{\nu}z^m\bigg)
\big(1+q\sum_{\ell=1}^s\sum_{j=0}^2 c_{\fod_\ell}\langle n_{\fod},m_j
\rangle z^{m_{\fod_\ell}+t_j}\big)\\
&=&\exp\big(q(W_k(Q)-W_0(Q))\big)\\
&&+\sum_{(c,\nu,m)\in S_k(Q)}
\sum_{\ell=1}^s
\bigg(c_{\fod_\ell}c q^{\nu}\big(\langle n_{\fod},r(m)\rangle z^{m+m_{\fod_\ell}}
+q\sum_{j=0}^2 \langle n_{\fod},m_j\rangle z^{m+m_{\fod_\ell}+t_j}\big)\bigg).
\end{eqnarray*}

We interpret this as follows. For each $(c,\nu,m)\in S_k(Q)$ and
each $\ell$, look at the four terms
\[
c_{\fod_\ell}c q^{\nu}\big(\langle n_{\fod},r(m)\rangle z^{m+m_{\fod_{\ell}}}
+q\sum_{j=0}^2 \langle n_{\fod},m_j\rangle z^{m+m_{\fod_\ell}+t_j}\big).
\]
These four terms contribute the expression
\[
c_{\fod_\ell}cq^{3d+\nu-|m+m_{\fod_\ell}|} \big(\langle n_{\fod},r(m)
\rangle
D_i(d,m+m_{\fod_\ell})+\sum_{j=0}^2\langle n_{\fod},m_j\rangle
D_i(d,m+m_{\fod_\ell}+t_j)\big)
\]
to $L^d_i(Q')$. One can check that in fact this total contribution is
zero, either by direct but tedious checking from the formulas
for $D_i$, or by applying Lemma \ref{twistedderhamlemma} with
$f=cq^{\nu-1}z^m$ and $\theta=\theta_{\gamma,\foD}$.

Now if $m+m_{\fod_\ell}$ and $m+m_{\fod_\ell}+t_j$, $0\le j\le 2$,
all lie in the \emph{same} $\tilde\omega_d$, then these terms
produce no total contribution to $L^d_{i,\tau}(Q')$ for any $\tau\in\Sigma$,
including $\tau=\omega$. On the other hand, these four terms can
contribute to different $L^d_{i,\omega}(Q')$'s
if $m+m_{\fod_{\ell}}$ and $m+m_{\fod_{\ell}}+t_j$, $j=0,1,2$, don't 
all lie in $\tilde\omega_d$ for the same $\omega\in\Sigma$.
This can happen only if $m+m_{\fod}\in\tilde\omega_d$ but $m+m_{\fod}+t_j
\in\tilde\tau_d$ for some $j$ with $\omega\subsetneq\tau\in\Sigma$
with $\dim\tau=\dim\omega+1$ and $m_j\in\tau$, $m_j\not\in\omega$.
In this case, $L^d_{i,\tau}(Q')-L^d_{i,\tau}(Q)$ has a contribution
of the form $cc_{\fod_{\ell}}\langle n_{\fod},m_j\rangle q^{3d+\nu-
|m+m_{\fod_{\ell}}|}
D_i(d,m+m_{\fod_{\ell}}+t_j)$. Thus
$L^d_{i,\omega}(Q')-L^d_{i,\omega}(Q)$ must have a contribution coming
from the same term, but with opposite sign.
This gives the Lemma.
\qed

\medskip

We can now use the asymptotic behaviour of the expressions
$L^d_{i,\omega}(Q)$ and the above wall-crossing formula to
rewrite the needed expressions:

\begin{lemma}
\label{basiclooplemma}
Let $\gamma_j$ be the straight line path joining
$Q$ with $Q+sm_j$ for $s\gg 0$.
Let $\gamma_{j,j+1}$ be the loop based at $Q$
which passes linearly from $Q$ to $Q+sm_j$, then takes
a large circular arc to $Q+sm_{j+1}$, and then proceeds linearly
from $Q+sm_{j+1}$ to $Q$. Here we take $j$ modulo $3$, and $\gamma_{j,j+1}$
is always a counterclockwise loop. Let $\sigma_{j,j+1}=\rho_j+\rho_{j+1}$,
a two-dimensional cone in $\Sigma$.
Then
\[
L^d_i(Q)-L^d_{i,\{0\}}(Q)=-\sum_{j=0}^2 L^d_{i,\gamma_j,\{0\}\rightarrow
\rho_j}
-\sum_{j=0}^2 L^d_{i,\gamma_{j,j+1},\rho_{j+1}\rightarrow \sigma_{j,j+1}}.
\]
\end{lemma}

\proof
By Lemma \ref{asympbehav},
$L^d_{i,\sigma}(Q+sm_j)=0$ for any $\sigma\in\Sigma$ with 
$\rho_j\subseteq\sigma$. Thus
by \eqref{eqrhochange} and \eqref{eqsigmachange}, we have
\begin{eqnarray*}
L^d_{i,\rho_j}(Q)&=&-L^d_{i,\gamma_j,\{0\}\rightarrow\rho_j}
+\sum_{\scriptstyle \sigma\in\Sigma\atop\scriptstyle 
\rho_j\subsetneq\sigma} L^d_{i,\gamma_j,
\rho_j\rightarrow\sigma},\\
L^d_{i,\sigma_{j,j+1}}(Q)&=&-\sum_{\scriptstyle \rho\in\Sigma\atop{\scriptstyle
\dim\rho=1
\atop\scriptstyle\rho\subsetneq\sigma_{j,j+1}}} L^d_{i,\gamma_j,\rho\rightarrow
\sigma_{j,j+1}}.
\end{eqnarray*}
Note we have broken symmetry for the second equation.

Adding together contributions from the $\rho_j$'s and $\sigma_{j_1,j_2}$'s,
we see from Lemma \ref{stupidLdsigmafacts}, (1), that
\begin{eqnarray*}
L^d_i(Q)-L^d_{i,\{0\}}(Q)&=&-\sum_{j=0}^2 L^d_{i,\gamma_j,\{0\}\rightarrow
\rho_j}\\
&&-(L^d_{i,\gamma_0,\rho_1\rightarrow\sigma_{0,1}}
-L^d_{i,\gamma_1,\rho_1\rightarrow\sigma_{0,1}})\\
&&-(L^d_{i,\gamma_1,\rho_2\rightarrow\sigma_{1,2}}
-L^d_{i,\gamma_2,\rho_2\rightarrow\sigma_{1,2}})\\
&&-(L^d_{i,\gamma_2,\rho_0\rightarrow\sigma_{2,0}}
-L^d_{i,\gamma_0,\rho_0\rightarrow\sigma_{2,0}}).
\end{eqnarray*}
Again by Lemma \ref{asympbehav}, it follows that the contribution to
$L^d_{i,\gamma_{j,j+1},\rho_{j+1}\rightarrow \sigma_{j,j+1}}$
from the large circular arc is zero. Hence
\[
L^d_i(Q)-L^d_{i,\{0\}}(Q)=-\sum_{j=0}^2 L^d_{i,\gamma_j,\{0\}\rightarrow
\rho_j}
-\sum_{j=0}^2 L^d_{i,\gamma_{j,j+1},\rho_{j+1}\rightarrow \sigma_{j,j+1}},
\]
the desired result. \qed

\medskip

We have already interpreted $L^d_{i,\{0\}}(Q)$ in Lemma
\ref{stupidLdsigmafacts}, (2),
so it remains to interpret the remaining terms on the right-hand side
of the above lemma.

\begin{lemma}
\label{Llemma}
\begin{equation*}
-L^d_{i,\gamma_j,\{0\}\rightarrow\rho_j}=\sum_{\nu\ge i-1}\sum_{
\scriptstyle I\subseteq \{1,
\ldots,k\}\atop {\scriptstyle 
I=\{i_1,\ldots,i_{3d-2+i-\nu}\}\atop \scriptstyle
i_1<\cdots<i_{3d-2+i-\nu}}}
\langle P_{i_1},\ldots,P_{i_{3d-2+i-\nu}},\psi^{\nu}S
\rangle^{\trop}_{d,\rho_j} u_Iq^{\nu+2-i}
\end{equation*}
for $S=Q,L$ or $M_{\RR}$ for $i=0,1$ and $i=2$ respectively. Here,
as usual, $L$ is a tropical line with vertex $Q$.
\end{lemma}

\proof
This is vacuous for $i=0$, so we assume $i\ge 1$.
Without loss of generality, consider $L^d_{i,\gamma_0,\{0\}\rightarrow
\rho_0}$. This quantity is a sum of contributions from each
point $P\in Q+(\rho_0\setminus \{0\})$ which is the intersection of
$Q+\rho_0$ with a ray $\fod\in\foD=\foD(\Sigma,P_1,\ldots,P_k)$.
Let us consider the contribution to $L^d_{i,\gamma_0,\{0\}\rightarrow\rho_0}$
from a small segment $\gamma$ of $\gamma_0$ which only crosses $\fod$.
Let $\gamma$ run from $Q_1$ to $Q_2$. Now $\fod$ corresponds to a
Maslov index zero tree passing through $P$, and by cutting it at
$P$, we obtain a Maslov index zero disk $h_1:\Gamma'_1\rightarrow M_{\RR}$
with boundary $P$. Then
\[
f_{\fod}=1+w_{\Gamma_1'}(E_{\out,1})\Mult(h_1)z^{\Delta(h_1)}u_{I(h_1)}.
\]
Furthermore a term $cz^mq^{\nu}$ in $\exp(q(W_k(Q_1)-W_0(Q_1)))$
arises
from $\nu$ distinct Maslov index two disks with boundary $Q_1$, say
$h_2,\ldots,h_{\nu+1}$ (each with at least one marked point), 
and the term contributed is
\[
q^{\nu}\prod_{i=2}^{\nu+1}\Mult(h_i)z^{\Delta(h_i)}u_{I(h_i)}.
\]
In order for this term to contribute to $L^d_{i,\gamma,\{0\}\rightarrow
\rho_0}$, $m+m_{\fod}=
\sum_{i=1}^{\nu+1}\Delta(h_i)$ must be of the form $dt_0+n_1t_1+n_2t_2$
with $n_1,n_2\le d$.
The disks $h_2,\ldots,h_{\nu+1}$ 
deform to disks with boundary at $P$, which we also call
$h_2,\ldots,h_{\nu+1}$. Write these disks
as $h_i:\Gamma_i'\rightarrow M_{\RR}$. Each $\Gamma_i'$, $1\le i\le\nu+1$,
has a vertex $V_{\out,i}$. 

Using this data, we can construct an actual tropical curve as follows.
Let $\Gamma$ be the graph obtained by identifying all the outgoing
vertices $V_{\out,i}$ in $\Gamma_1',\ldots,\Gamma_{\nu+1}'$, to obtain
a graph with a distinguished vertex $V_{\out}$, and then attaching
$(d-n_1)+(d-n_2)+1$ additional unbounded edges with vertex $V_{\out}$.
We then define $h:\Gamma\rightarrow M_{\RR}$ to agree with $h_i$
on $\Gamma_i'\subseteq\Gamma$. We have $h$ taking 
the first $d-n_1$ new unbounded edges to $P+\RR_{\ge 0}m_1$; the
second $d-n_2$ new unbounded edges to $P+\RR_{\ge 0}m_2$; and the
last unbounded edge is contracted, and marked with the label $x$.
Note $\Gamma$ has valency at $V_{\out}$ given by
$\Val(V_{\out})=\nu+1+(d-n_1)+(d-n_2)+1$. Thus
we obtain a parameterized curve $h:\Gamma\rightarrow M_{\RR}$ with $h(x)=P$. 
The
balancing condition needs to be checked at $V_{\out}$, but as in
\eqref{balanceeq1} and \eqref{balanceeq2}, the fact that 
\[
\sum_{i=1}^{\nu+1}\Delta(h_i)+(d-n_1)t_1+(d-n_2)t_2=d(t_0+t_1+t_2)
\]
shows the balancing condition indeed holds at $V_{\out}$.

The contribution of this term to $-L^d_{i,\gamma,\{0\}\rightarrow
\rho_0}$ is
\begin{equation}
\label{temporaryanswer}
-\langle n_{\fod},m_0\rangle u_{I(h)}w_{\Gamma_1'}(E_{\out,1})
\bigg(\prod_{i=1}^{\nu+1}\Mult(h_i)\bigg)D_i(d,d+1,n_1,n_2)
q^{\nu+3d-(d+n_1+n_2)}.
\end{equation}
Note $n_{\fod}$ is primitive, annihilates $r(m_{\fod})$, and must be
positive on $-m_0$. Furthermore, after choosing an isomorphism
$\bigwedge^2 M\cong\ZZ$, $w(E_{1,\out})n_{\fod}$ can be identified,
up to sign, with $X_{r(m_{\fod})}$. Thus setting $m(h_1)
=r(m_{\fod})$ as in Definition \ref{tropicaldescinv}, we see that
\[
-\langle n_{\fod},m_0\rangle w_{\Gamma'_1}(E_{\out,1})=|m(h_1)\wedge m_0|.
\]
Thus \eqref{temporaryanswer} coincides with 
\[
|m(h_1)\wedge m_0|u_{I(h)}D_i(d,d+1,n_1,n_2)
\bigg(\prod_{V\in\Gamma^{[0]}\atop V\not\in E_x} \Mult_V(h)\bigg)
q^{\Val(V_{\out})-2}.
\]
Now $D_i(d,d+1,n_1,n_2)=\Mult^{i-1}_x(h)$ as defined in
Definition \ref{tropicaldescinv} via direct comparison with the
definitions of the $D_i$'s. Furthermore, if $I(h)=\{i_1,\ldots,
i_{3d-2+i-\nu'}\}$ for some $\nu'$, we see that, as $h$ is obtained
by gluing one Maslov index zero disk to $\Val(V_{\out})-2$
Maslov index two disks, we have
\begin{eqnarray*}
\Val(V_{\out})-2&=&\sum_{i=1}^{\nu+1}\big(|\Delta(h_i)|-\#I(h_i)\big)+(d-n_1)+(d-n_2)\\
&=&3d-(3d-2+i-\nu')=\nu'+2-i.
\end{eqnarray*}
Thus, by Definition \ref{tropicaldescinv},
the term under consideration contributes to $-L^d_{i,\gamma_0,\{0\}\rightarrow
\rho_0}$ by exactly the same amount that the curve $h$ contributes
to
\[
\langle P_{i_1},\ldots,P_{i_{3d-2+i-\nu'}},\psi^{\nu'}S
\rangle^{\trop}_{d,\rho_0}u_{I(h)}q^{\nu'+2-i},
\]
as desired.

Conversely, given any curve $h$ contributing to 
$\langle P_{i_1},\ldots,P_{i_{3d-2+i-\nu'}},\psi^{\nu'}S
\rangle^{\trop}_{d,\rho_0}$ with $h(E_x)=P\in \rho_0\setminus \{0\}$,
the procedure of Lemma \ref{curvechop} shows $h$ must arise in
precisely the way described above. 
\qed

\medskip

\begin{lemma}
\label{MRlemma}
\[
-L^d_{i,\gamma_{j,j+1},\rho_{j+1}
\rightarrow\sigma_{j,j+1}}=
\sum_{\hbox{$i$ s.t.}\atop \scriptstyle P_i\in\sigma_{j,j+1}}
u_iq^{-1}+
\sum_{\nu\ge 0}\sum_{\scriptstyle
I\subseteq \{1,
\ldots,k\}\atop {
\scriptstyle
I=\{i_1,\ldots,i_{3d-2+i-\nu}\}\atop
\scriptstyle
i_1<\cdots<i_{3d-2+i-\nu}}}
\langle P_{i_1},\ldots,P_{i_{3d-2+i-\nu}},\psi^{\nu}S
\rangle^{\trop}_{d,\sigma_{j,j+1}} u_Iq^{\nu+2-i}
\]
for $S=Q$, $L$ or $M_{\RR}$ for $i=0,1$ and $i=2$ respectively.
\end{lemma}

\proof
First note this is vacuous for $i=0$ or $1$ as both sides are zero, so we
can assume $i=2$.
Second, if $\shC_1$ and $\shC_2$ are closures of
two connected components of
$M_{\RR}\setminus\Supp\foD$, where $\foD=\foD(\Sigma,P_1,\ldots,P_k)$, 
$\dim\shC_1\cap\shC_2=1$, and $\gamma$ is a short path from $\Int(\shC_1)$
into $\Int(\shC_2)$ just crossing $\Int(\shC_1\cap\shC_2)$ once,
then $L^d_{i,\gamma,\rho_{j+1}\rightarrow\sigma_{j,j+1}}$ is
independent of $\gamma$ and its endpoints. Furthermore, 
reversing the direction of $\gamma$ changes the sign of
$L^d_{i,\gamma,\rho_{j+1}\rightarrow\sigma_{j,j+1}}$. 
So a simple homological argument shows that
\[
L^d_{i,\gamma_{j,j+1},\rho_{j+1}\rightarrow\sigma_{j,j+1}}
=\sum_{P\in \Sing(\foD)\cap\sigma_{j,j+1}}
L^d_{i,\gamma_{P},\rho_{j+1}\rightarrow\sigma_{j,j+1}}
\]
where $\gamma_P$ is a small counterclockwise loop around the singular
point $P$. This localizes the calculation to the singular points of
$\foD$ in $\sigma_{j,j+1}$. Now such a singular point $P$ is either
in $\{P_1,\ldots,P_k\}$ or not; this will give us cases (3) (b) and
(3) (a) of Definition \ref{tropicaldescinv} respectively. 
To save on typing, we set
\[
L_{P,j}:=L^d_{2,\gamma_P,\rho_{j+1}\rightarrow\sigma_{j,j+1}}.
\]

\emph{Case 1}. $P\not\in\{P_1,\ldots,P_k\}$. Fix a base-point $Q'$ near
$P$. Consider a term $cq^{\nu}z^m$ in $\exp(q(W_k(Q')-W_0(Q')))$ of the form
\begin{equation}
\label{caseoneeq}
cq^{\nu}z^m=q^{\nu}\prod_{i=3}^{\nu+2} \big(\Mult(h_i)z^{\Delta(h_i)}
u_{I(h_i)}\big)
\end{equation}
where the $h_i$'s are Maslov index two disks with boundary $Q'$, but
none of the $h_i$'s come from broken lines which bend near
$P$. As a result, this term appears in $\exp(q(W_k(Q'')-W_0(Q'')))$ 
for all $Q''$ general in
a small open neighbourhood of $P$.

Suppose that such a term $cq^{\nu}z^m$ contributes to $L_{P,j}$
when $\gamma_P$ crosses a ray
$\fod\in\foD$ with $P\in\fod$, $P\not=\Init(\fod)$. But
$\gamma_P$ crosses this ray $\fod$ twice, in opposite directions,
so $cq^{\nu}z^m$ will contribute to $L_{P,j}$ twice, but
with opposite signs. Thus these contributions cancel, and don't
contribute to the total in $L_{P,j}$.

Thus we only need analyze contributions arising when $\gamma_P$
crosses a ray $\fod$ with $\Init(\fod)=P$ or contributions from
monomials as in \eqref{caseoneeq} where some of the $h_i$'s come from
broken lines which do bend near $P$.
As in the proof of Theorem \ref{wallcrossingtheorem}, 
we can in fact assume that there are precisely three rays,
$\fod_1$, $\fod_2$, $\fod_3$ passing through $P$, with $\Init(\fod_1),
\Init(\fod_2)\not=P$ and $\Init(\fod_3)=P$. Now $\fod_1,\fod_2$
correspond to Maslov index zero trees passing through $P$,
and by cutting them, we obtain Maslov index zero disks
$h_i:\Gamma_i'\rightarrow M_{\RR}$, $i=1,2$ with boundary $P$,
and for $i=1,2$,
\[
f_{\fod_i}=1+w_{\Gamma_i'}(E_{\out,i})\Mult(h_i)z^{\Delta(h_i)}u_{I(h_i)}.
\]

We now analyze how additional terms $cq^{\nu}z^m$ which can contribute
to $L_{P,j}$ may arise. In what follows, assume that $cq^{\nu}z^m$
is as in \eqref{caseoneeq}
in which none of the broken lines corresponding to $h_3,\ldots,h_{\nu+2}$
bend at $\fod_1$, $\fod_2$ or $\fod_3$.

Write 
\[
m+\Delta(h_1)+\Delta(h_2)=
\sum_{j=0}^2 n_jt_j.
\]

We have the following possibilities of additional contributions:
\begin{itemize}
\item[(I)] $cq^{\nu}z^m$ may contribute to $L_{P,j}$ when $\gamma_P$
crosses $\fod_3$. This contribution can only occur if $n_{j+2}\le n_j=d
<n_{j+1}$.
\item[(II)] 
After crossing $\fod_1$, new terms of the form (leaving
off the coefficients) $z^{m+\Delta(h_1)}$ and $z^{m+\Delta(h_1)+t_k}$,
$k=0,1,2$ may appear in $\exp(q(W_k-W_0))$. Thus when we cross
$\fod_2$, these new terms may contribute to $L_{P,j}$. 
Note that
$z^{m+\Delta(h_1)}$ only contributes when crossing $\fod_2$ if $n_{j+2}
\le n_j=d<n_{j+1}$. The term $z^{m+\Delta(h_1)+t_j}$ 
only contributes if 
$n_{j+2}\le d$, $n_j= d-1$ and $d<n_{j+1}$. The term
$z^{m+\Delta(h_1)+t_{j+1}}$ only contributes if $n_{j+2}\le n_j=d\le n_{j+1}$.
The term $z^{m+\Delta(h_1)+t_{j+2}}$ only contributes if
$n_{j+2}< n_j=d<n_{j+1}$.
\item[(III)]
After crossing $\fod_2$, new terms of the form (leaving
off the coefficients) $z^{m+\Delta(h_2)}$ and $z^{m+\Delta(h_2)+t_k}$,
$k=0,1,2$ may appear in $\exp(q(W_k-W_0))$. Thus when we cross
$\fod_1$, these new terms may contribute to $L_{P,j}$.
Note that
$z^{m+\Delta(h_2)}$ only contributes when crossing $\fod_2$ if $n_{j+2}
\le n_j=d<n_{j+1}$. The term $z^{m+\Delta(h_2)+t_j}$ 
only contributes if 
$n_{j+2}\le d$, $n_j=d-1$ and $d<n_{j+1}$. The term
$z^{m+\Delta(h_2)+t_{j+1}}$ only contributes if $n_{j+2}\le n_j=d\le n_{j+1}$.
The term $z^{m+\Delta(h_2)+t_{j+2}}$ only contributes if
$n_{j+2}< n_j=d<n_{j+1}$.
\end{itemize}
There are now three cases when these additional contributions to $L_{P,j}$
occur.

\emph{Case 1(a).} 
$n_{j+2}\le n_j=d<n_{j+1}$. In this case, (leaving off the
coefficients), $z^m$ gives a
contribution to $L_{P,j}$ of type (I)
when $\gamma_P$ crosses $\fod_3$, and $z^{m+\Delta(h_i)}$,
$z^{m+\Delta(h_i)+t_{j+1}}$,
or $z^{m+\Delta(h_i)+t_{j+2}}$ (if $n_{j+2}<d$) may give contributions
of type (II) or (III) when $\gamma_P$ crosses $\fod_1$ and $\fod_2$.
Now note that the total change
to $L^d_{i,\sigma_{j,j+1}}$ due to these monomials as we traverse the loop
$\gamma_P$ is 
the sum of the contributions of these monomials to
$L^d_{i,\gamma_P,\rho_j\rightarrow\sigma_{j,j+1}}$ and 
$L^d_{i,\gamma_P,\rho_{j+1}\rightarrow\sigma_{j,j+1}}$. However, the
total contribution to the change of $L^d_{i,\sigma_{j,j+1}}$
is necessarily zero, as $\gamma_P$ is a loop, and because $d<n_{j+1}$,
none of these monomials contribute to any change of 
$L^d_{i,\gamma_P,\rho_j\rightarrow\sigma_{j,j+1}}$. Hence the total
contribution of these monomials to 
$L^d_{i,\gamma_P,\rho_{j+1}\rightarrow\sigma_{j,j+1}}=L_{P,j}$
is also zero.

\emph{Case 1(b).}
$n_{j+2}\le d$, $n_j=d-1$, $d<n_{j+1}$. In this case only
the terms $z^{m+\Delta(h_i)+t_j}$, $i=1,2$, may contribute. However,
the same argument as in Case 1(a) shows that the total contribution from
these terms is zero.

\emph{Case 1(c).} 
$n_{j+2}\le n_j=n_{j+1}=d$. In this case, contributions to 
$L_{P,j}$ only arise from terms
of the form $z^{m+\Delta(h_i)+t_{j+1}}$. 
Choose $n_{\fod_i}$ so that at the first time $t_i$ when
$\gamma_P$ passes through $\fod_i$, $\langle n_{\fod_i},\gamma'_P(t_i)\rangle
<0$. By interchanging the labelling of $\fod_1$ and $\fod_2$
and choosing the base-point $Q'$ appropriately, we can assume firstly
that $\gamma_P$ passes initially through $\fod_1$ and then $\fod_2$,
and secondly that $\langle n_{\fod_i},m_{j+1}\rangle\ge 0$ for $i=1,2$.
Write $f_{\fod_i}=1+c_{\fod_i}z^{m_{\fod_i}}$ for $i=1,2$.

Then the term 
\[
\langle n_{\fod_1},m_{j+1}\rangle cc_{\fod_1}q^{\nu+1}
z^{m+\Delta(h_1)+t_{j+1}}
\]
appears in $\exp\big(q(W_k(\gamma_P(t))-W_0(\gamma_P(t)))\big)$ right
after $\gamma_P$ crosses
$\fod_1$ the first time (and disappears when we cross $\fod_1$ for the
second time), and hence when $\gamma_P$ crosses
$\fod_2$ for the first time, we obtain a contribution to
$L_{P,j}$ of
\[
\langle n_{\fod_2},m_j\rangle\langle
n_{\fod_1},m_{j+1}\rangle cc_{\fod_1}c_{\fod_2}
D_2(d,m+\Delta(h_1)+\Delta(h_2)+t_j+t_{j+1})
q^{\nu+3d-|m+\Delta(h_1)+\Delta(h_2)|}.
\]
On the other hand, the term $\langle n_{\fod_2},m_{j+1}\rangle
cc_{\fod_2} q^{\nu+1}z^{m+\Delta(h_2)+t_{j+1}}$ appears
after $\gamma_P$ crosses $\fod_2$ for the first time (and disappears
when we cross $\fod_2$ for the second time), and hence
when $\gamma_P$ crosses $\fod_1$ for the second time, we obtain
a contribution to $L_{P,j}$ of 
\[
\langle -n_{\fod_1},m_j\rangle\langle
n_{\fod_2},m_{j+1}\rangle cc_{\fod_1}c_{\fod_2}
D_2(d,m+\Delta(h_1)+\Delta(h_2)+t_j+t_{j+1})
q^{\nu+3d-|m+\Delta(h_1)+\Delta(h_2)|}.
\]
Note
\begin{eqnarray*}
\langle n_{\fod_2},m_j\rangle\langle
n_{\fod_1},m_{j+1}\rangle 
-\langle n_{\fod_1},m_j\rangle\langle
n_{\fod_2},m_{j+1}\rangle
&=&-|n_{\fod_1}\wedge n_{\fod_2}|\\
&=&-|m^{\prim}(h_1)\wedge m^{\prim}(h_2)|
\end{eqnarray*}
as $n_{\fod_1}, n_{\fod_2}$ form a positively oriented basis of
$N_{\RR}$, and $m_j,m_{j+1}$ form a positively oriented basis of $M_{\RR}$. 

Now the Maslov index two disks $h_3,\ldots,h_{\nu+2}$ deform to
disks with boundary $P$, which we also call $h_3,\ldots,h_{\nu+2}$.
We can then glue together the disks $h_1,\ldots,h_{\nu+2}$
along with $d-n_{j+2}$ copies of the Maslov index two disks
with no marked points in the direction $m_{j+2}$. These are glued
at their respective outgoing vertices, yielding a vertex $V_{\out}$,
and we add one additional unbounded edge $E_x$ with the label $x$,
also attached to the vertex $V_{\out}$. This yields a graph
$\Gamma$, whose valency at $V_{\out}$ is
$\Val(V_{\out})=\nu+3+d-n_{j+2}$. Thus we obtain a parameterized
curve $h:\Gamma\rightarrow M_{\RR}$ with $h(x)=P$. Again one
easily checks the balancing condition at $V_{\out}$.

Thus the total contribution arising in the ways
analyzed from $cq^{\nu}z^m$ to $-L_{P,j}$ is
\begin{eqnarray*}
&&|m^{\prim}(h_1)\wedge m^{\prim}(h_2)|w_{\Gamma_1'}(E_{\out,1})
w_{\Gamma_2'}(E_{\out,2})
\Mult^0_x(h)\bigg(\prod_{V\in\Gamma^{[0]}
\atop V\not\in E_x} \Mult_V(h)\bigg)q^{\nu+3d-(2d+n_{j+2})}u_{I(h)}\\
&=&|m(h_1)\wedge m(h_2)|
\Mult^0_x(h)\bigg(\prod_{V\in\Gamma^{[0]}
\atop V\not\in E_x} \Mult_V(h)\bigg)q^{\Val(V_{\out})-3}u_{I(h)}.
\end{eqnarray*}
One sees that if $I(h)=\{i_1,\ldots,i_{3d-\nu'}\}$ for some $\nu'$, then
since $h$ is obtained by gluing two Maslov index zero disks
with $\Val(V_{\out})-3$ Maslov index two disks, we have
\begin{eqnarray*}
\Val(V_{\out})-3&=&\sum_{i=1}^{\nu+2} (|\Delta(h_i)|-\#I(h_i))+d-n_{j+2}\\
&=&3d-(3d-\nu')=\nu'.
\end{eqnarray*}
Thus we see that the coefficient of the contributions analyzed above
from $cq^{\nu}z^m$ to $-L_{P,j}$
is precisely the contribution of $h$ to
\begin{equation}
\label{step6eq1}
\langle P_{i_1},\ldots,P_{i_{3d-\nu'}},
\psi^{\nu'}M_{\RR}\rangle^{\trop}_{d,\sigma_{j,j+1}}u_{I(h)}q^{\nu'}
\end{equation}
as desired.

Conversely, given an $h$ contributing to \eqref{step6eq1} with $h(x)=P$,
one can cut it at $P$, using Lemma \ref{curvechop}, 
decomposing it into tropical disks. Then we see $h$ arises precisely as above.
Thus we see that $-L_{P,j}$ is the contribution
to \eqref{step6eq1} from maps with $h(E_x)=P$.

\emph{Case 2}. $P=P_i$ for some $i$.
Again, choose a basepoint $Q'$ near $P_i$. 
By Remark \ref{Piremark}, there is a one-to-one correspondence
between rays in $\foD$ containing $P_i$ and Maslov index two disks
with boundary $P_i$ not having $P_i$ as a marked point. With
$Q'$ sufficiently near $P_i$, these Maslov index two disks deform
to ones with boundary at $Q'$, so 
the Maslov index two disks with boundary $P_i$ not having
$P_i$ as a marked point are in
one-to-one correspondence with the Maslov index two disks with boundary $Q'$
not having $P_i$ as a marked point.

If we are interested in
terms in $\exp(q(W_k(Q')-W_0(Q')))$ 
which may contribute to $L_{P_i,j}$,
we only need to look at those terms in $\exp(q(W_k(Q')-W_0(Q')))$ 
which
do not have $u_i$  as a factor, as any term that does will
not produce any new terms as we cross a ray through $P_i$.
So consider a term $cq^{\nu}z^m$ of the form
\begin{equation}
\label{secondcaseterm}
cq^{\nu}z^m=q^{\nu}\prod_{p=1}^{\nu} \Mult(h_p)z^{\Delta(h_p)} u_{I(h_p)},
\end{equation}
where each of these disks 
$h_p$ with boundary $Q'$ does not pass through $P_i$, and hence
corresponds to a disk with boundary
$P_i$, which we also write as $h_p:\Gamma_p'\rightarrow M_{\RR}$.
By extending these disks to trees and marking $P_i$, we obtain Maslov
index zero trees, corresponding to rays $\fod_p$
in $\foD$ with initial point $P_i$. 
In addition, we have rays $\foc_p\in \foD$, $p=0,1,2$,
with initial point $P_i$, corresponding to the three Maslov index
two disks with boundary $Q'$ with no marked points. These do not appear
in $W_k(Q')-W_0(Q')$, so are distinct from the $\fod_p$'s.

In what follows, we 
write $m=\sum_{p=1}^{\nu} \Delta(h_p)=\sum_{j=0}^2 n_jt_j$,
and take $n_{\fod_p}$ and $n_{\foc_p}$ to have their sign chosen
so that they are negative on $\gamma'_P$ when $\gamma_P$ crosses
the corresponding ray. Note that as $\gamma_P$ is counterclockwise,
if we use the identification
$\bigwedge^2 M\cong\ZZ$ given by the standard orientation,
i.e., $m_1\wedge m_2\mapsto 1$, then $w_{\Gamma'_p}(E_{\out,j})n_{\fod_p}
=X_{r(\Delta(h_p))}$. So
\begin{equation}
\label{balancingequation}
\sum_{p=1}^{\nu} w_{\Gamma'_p}(E_{\out,j})n_{\fod_p}
=X_{r(m)}.
\end{equation}
On the other hand,
$\langle n_{\foc_{j+1}},m_j\rangle=-1$ and
$\langle n_{\foc_{j+2}},m_j\rangle=1$.

We can now view this term $cq^{\nu}z^m$ as giving rise to contributions
to $L_{P,j}$ in the following four ways: 
\begin{itemize}
\item[(I)] $\gamma_P$ crosses $\fod_l$ for some $1\le l\le \nu$.
Then the term 
\[
q^{\nu-1}\prod_{p=1\atop p\not=l}^{\nu}
\Mult(h_p)z^{\Delta(h_p)} u_{I(h_p)}
\]
contributes to $L_{P,j}$ if
$n_{j+2}\le n_j=d<n_{j+1}$, in which case the contribution is
\[
\langle n_{\fod_l},m_j\rangle \bigg(\prod_{p=1}^{\nu}
\Mult(h_p)u_{I(h_p)}\bigg) u_iw_{\Gamma'_l}(E_{\out,l})
D_2(d,m+t_j)q^{\nu+3d-|m|-1}.
\]
Note such a contribution requires $\nu>0$.
\item[(II)] $\gamma_P$ crosses $\foc_j$. If $cq^{\nu}z^m$ 
contributes to $L_{P,j}$ when $\gamma_P$ crosses $\foc_j$, 
its contribution would involve a factor
of $\langle n_{\foc_j},m_j\rangle=0$, hence there is no contribution.
\item[(III)] $\gamma_P$ crosses $\foc_{j+1}$. We get a contribution
from $cq^{\nu}z^m$
if $n_{j+2}\le n_j=d\le n_{j+1}$, in which case the contribution is
\[
\langle n_{\foc_{j+1}},m_j\rangle \bigg(\prod_{p=1}^{\nu}
\Mult(h_p)u_{I(h_p)}\bigg) u_i
D_2(d,m+t_j+t_{j+1})q^{\nu+3d-|m|-1}.
\]
\item[(IV)] $\gamma_P$ crosses $\foc_{j+2}$. We get a contribution
from $cq^{\nu}z^m$
if $n_{j+2}<n_j=d<n_{j+1}$, in which case we get
\[
\langle n_{\foc_{j+2}},m_j\rangle \bigg(\prod_{p=1}^{\nu}
\Mult(h_p)u_{I(h_p)}\bigg) u_i
D_2(d,m+t_j+t_{j+2})q^{\nu+3d-|m|-1}.
\]
\end{itemize}
We now consider three cases. 

\emph{Case 2(a).} $n_{j+2}\le n_j=d<n_{j+1}$, $\nu>0$. In this
case ignoring the common factors 
\[
q^{\nu+3d-|m|-1}u_i\prod_{p=1}^{\nu}
\Mult(h_p)u_{I(h_p)},
\]
the total contribution is, using 
\eqref{balancingequation} and Lemma \ref{unbasicintegral},
\begin{eqnarray*}
&&\langle\sum_{l=1}^{\nu}w_{\Gamma_l'}(E_{\out,l})n_{\fod_l},
m_j\rangle D_2(d,m+t_j)
-D_2(d,m+t_j+t_{j+1})\\
&&+\begin{cases}D_2(d,m+t_j+t_{j+2})&n_{j+2}<d\\
0&n_{j+2}=d\end{cases}\\
&=&r(m)\wedge m_j (-1)^{n_j+n_{j+1}+1}
{(n_j-d)!(n_{j+1}-d-1)!\over (d-n_{j+2})!}\\
&&-(-1)^{n_j+n_{j+1}+2}{(n_j-d)!(n_{j+1}-d)!\over (d-n_{j+2})!}\\
&&+(-1)^{n_j+n_{j+1}+1}{(n_j-d)!(n_{j+1}-d-1)!\over (d-n_{j+2})!}
(d-n_{j+2})\\
&=&\big((n_{j+2}-n_{j+1})+(n_{j+1}-d)+(d-n_{j+2})\big)
 (-1)^{n_j+n_{j+1}+1}
{(n_j-d)!(n_{j+1}-d-1)!\over (d-n_{j+2})!}\\
&=&0.
\end{eqnarray*}
So there is no contribution to $L_{P,j}$ from this case.

\emph{Case 2(b).} $n_{j+2}\le n_j=d=n_{j+1}$, $\nu>0$.
In this case we only get
a contribution from (III). In this case, we can glue together
the disks $h_1,\ldots,h_{\nu}$ along with $d-n_{j+2}$ copies
of the Maslov index two disk with no marked points in the direction
$m_{j+2}$. These are glued at their respective outgoing vertices,
yielding a vertex $V_{\out}$, and we add two additional marked unbounded
edges $E_x$ and $E_{p_l}$ for some $l$ attached to $V_{\out}$.
This yields a graph $\Gamma$, whose valency at $V_{\out}$
is $\Val(V_{\out})=\nu+(d-n_{j+2})+2=\nu+3d-|m|+2$. Thus we obtain
a parameterized curve $h:\Gamma\rightarrow M_{\RR}$ with $h(V_{\out})=h(x)=
h(p_l)=P_i$.
The contribution
to $-L_{P,j}$ is then easily seen by inspection to be
\begin{equation}
\label{step6eq2}
\Mult^0_x(h)\left(\prod_{V\in\Gamma^{[0]}\atop V\not\in E_x}
\Mult_V(h)\right) u_{I(h)}q^{\Val(V_{\out})-3}.
\end{equation}
Suppose that $I(h)=\{i_1,\ldots,i_{3d-\nu'}\}$ for some
$\nu'$, recalling $i\in I(h)$ since we added the marked edge
$E_{p_l}$ mapping to $P_i$. Since $h$ is obtained by
gluing $\Val(V_{\out})-2$ Maslov index two disks, we have
\begin{eqnarray*}
\Val(V_{\out})-2&=&\sum_{i=1}^{\nu}(|\Delta(h_i)|-\#I(h_i))+d-n_{j+2}\\
&=&3d-(3d-\nu'-1)=\nu'+1.
\end{eqnarray*}
Thus we see that \eqref{step6eq2} is precisely
the contribution of $h$ to
\[
\langle P_{i_1},\ldots,P_{i_{3d-\nu'}},\psi^{\nu'} M_{\RR}
\rangle^{\trop}_{d,\sigma_{j,j+1}} u_{I(h)}q^{\nu'}
\]
from Definition \ref{tropicaldescinv} (3) (b). As in the other
cases we have considered, conversely any such curve $h$ will
give rise to the correct monomial $cq^{\nu}z^m$ by cutting
the curve at $P$.

\emph{Case 2(c).} $\nu=0$. There is only one element $(c,\nu,m)\in S_k$
with $\nu=0$, namely $(1,0,0)$ corresponding to the constant monomial $1$.
So $n_0=n_1=n_2=0$ and we have no contribution unless $d=0$. Again,
this contribution to $L_{P,j}$ only arises from (III), and is
\[
\langle n_{\foc_{j+1}},m_j\rangle u_i D_2(0,t_j+t_{j+1})q^{-1}=-u_iq^{-1}.
\]
This gives the remaining claimed terms in $-L_{P,j}$.
\qed

\medskip

We finally have:

\medskip

\emph{Proof of Theorems \ref{PQinvariant} and \ref{maintheorem}.}
By Lemma \ref{threeforms}, it is enough to prove \eqref{eq56-},
\eqref{eq56} and \eqref{eq57}. We have \eqref{eq56-} by Lemma \ref{zerolemma}.
The remaining two equations follow from Lemmas \ref{stupidLdsigmafacts},
(2), \ref{basiclooplemma}, \ref{Llemma}, and \ref{MRlemma}.
\qed

\section{Generalizations}

I would like to comment briefly on possible generalizations of the results
of this paper. One expects the version of mirror symmetry considered
in this paper to hold more broadly for any non-singular projective
toric variety. It should be no particular problem generalizing
these results to $\PP^n$; one simply includes collections of
linear subspaces of all codimensions $\ge 2$ to obtain the
necessary perturbations. This produces higher dimensional 
scattering diagrams, and the arguments of \S 5 should go through,
although the formulas will get progressively more complicated.

More of an issue is generalizing to other toric varieties. As mentioned
in Remark \ref{correctLGremark}, a basic
problem with tropical geometry as currently understood is that
it can't account for curves with irreducible components mapping
into the toric boundary. Even for surfaces, the only cases 
in which tropical geometry correctly computes quantum cohomology
is $\PP^2$ and $\PP^1\times\PP^1$; in particular, the arguments
of this paper should work equally well for $\PP^1\times\PP^1$.
To go further, one needs to improve the tropical understanding
of these issues, which we leave for others. But the expectation
is that once one correctly counts possible disks, then the methods
of this paper should still work.

It is also possible to combine the ideas in this paper of using
tropical geometry to count Maslov index two disks with the
techniques of \cite{GS}. If one starts with an integral affine
manifold with singularities which is non-compact, then one
obtains something like a scattering diagram governing a degeneration
of non-compact varieties. This data is called a \emph{structure}
in \cite{GS}. One expects the structure to encode Maslov index zero
disks in the mirror. Then one can again describe a Landau-Ginzburg
potential in terms of broken lines. This will be taken up elsewhere.



\begin{thebibliography}{cccccccc}
\bibitem{Auroux} D.~Auroux: \emph{Mirror symmetry and $T$-duality in 
       the complement of an anticanonical divisor,}  J. G\"okova Geom.\ 
       Topol.\ GGT  {\bf 1}  (2007), 51--91. 
\bibitem{Bar} S.~Barannikov: \emph{Semi-infinite Hodge
        structures and mirror symmetry for projective spaces,} preprint, 2000.
\bibitem{ChanLeung} K.~Chan, N.-C.~Leung: \emph{Mirror symmetry for 
         toric Fano manifolds via SYZ transformations}, preprint, 2008. 
\bibitem{ChoOh} C.-H.~Cho, Y.-G.~Oh: \emph{Floer cohomology and disc 
         instantons of Lagrangian torus fibers in Fano toric manifolds,}
         Asian J.\ Math.\ {\bf 10} (2006), 773--814. 
\bibitem{CoxKatz} D.~Cox, S.~Katz:
        \emph{Mirror symmetry and algebraic geometry,}
        Mathematical Surveys and Monographs, 68. American Mathematical Society,
        Providence, RI, 1999.
\bibitem{DSI}  A.~Doaui, C.~Sabbah: \emph{Gauss-Manin systems,
        Brieskorn lattices and Frobenius structures. I,} 
        Proceedings of the International Conference in Honor of Fr\' ed\' eric 
        Pham (Nice, 2002). Ann. Inst. Fourier (Grenoble)  {\bf 53}  (2003),  
        1055--1116.
\bibitem{DSII}  A.~Doaui, C.~Sabbah: \emph{Gauss-Manin systems,
        Brieskorn lattices and Frobenius structures. II,} 
        Frobenius manifolds,  1--18, Aspects Math., E36, Vieweg, Wiesbaden, 
        2004.
\bibitem{Dub} B.~Dubrovin: \emph{Geometry of $2$D topological field theories,}
        Integrable systems and quantum groups (Montecatini Terme, 1993),  
        120--348, Lecture Notes in Math., 1620, Springer, Berlin, 1996.
\bibitem{FO3II} K.~Fukaya, Y.-G.~Oh, H.~Ohta, K.~Ono:
         \emph{Lagrangian Floer theory on compact toric manifolds II: Bulk
         deformations,} preprint, 2008.
\bibitem{Fulton} W.~Fulton, \emph{Introduction to toric varieties},
        Annals of Mathematics Studies, {\bf 131}, The William H. Roever 
		Lectures in Geometry. Princeton University Press, Princeton, 
		NJ, 1993. xii+157 pp.
\bibitem{GathMark} A.~Gathmann, H.~Markwig: \emph{Kontsevich's
        formula and the WDVV equations in tropical geometry,} Adv.\ Math.\
        {\bf 217} (2008) 537--560.
\bibitem{Giv} A.~Givental: \emph{Homological geometry and mirror symmetry,}
        Proceedings of the International Congress of Mathematicians, Vol. 1, 2 
        (Z\" urich, 1994),  472--480, Birkh\" auser, Basel, 1995.
\bibitem{Gr4} M.~Gross:
        \emph{The Strominger-Yau-Zaslow conjecture: From torus fibrations
        to degenerations}, to appear in Procedings of Symposia of
        Pure Mathematics, Seattle 2005. 
\bibitem{Announce} M.~Gross, B.~Siebert:
        \emph{Affine manifolds,  log structures, and mirror symmetry},
        Turkish J.\ Math.\ \textbf{27} (2003), 33--60.
\bibitem{logmirror} M.\ Gross, B.\ Siebert:
        \emph{Mirror symmetry via logarithmic degeneration data I},
        J.\ Differential Geom.~\textbf{72}
        (2006), 169--338.
\bibitem{part II} M.\ Gross, B.\ Siebert:
        \emph{Mirror symmetry via logarithmic degeneration data II},
        preprint 2007.
\bibitem{GS} M.\ Gross, B.\ Siebert: \emph{From real affine geometry
        to complex geometry}, preprint, 2007.
\bibitem{GPS} M.\ Gross, R.\ Pandharipande, B.\ Siebert:
        \emph{The tropical vertex}, preprint, 2009.
\bibitem{Hien} M.~Hien: \emph{Periods for flat algebraic connections,}
         preprint, 2008.
\bibitem{Iritani} H.~Iritani: \emph{Quantum $D$-modules and generalized
         mirror transformations,} Topology {\bf 47} (2008), 225--276.
\bibitem{iritani} H.~Iritani: \emph{An integral structure in quantum 
        cohomology and mirror symmetry for toric orbifolds,} preprint, 2009.
\bibitem{KaKoPa} L.\ Katzarkov, M.\ Kontsevich, T.\ Pantev: \emph{Hodge
         theoretic aspects of mirror symmetry,} preprint, 2008.
\bibitem{ks} M.\ Kontsevich, Y.\ Soibelman:
	\emph{Affine structures and non-Archimedean analytic spaces}, in:
	\textsl{The unity of mathematics} (P.~Etingof, V.~Retakh,
	I.M.~Singer, eds.),  321--385, Progr.\ Math.~244,
	Birkh\"auser~2006.
\bibitem{KSDT} M.\ Kontsevich, Y.\ Soibelman: \emph{Stability
        structures, motivic Donaldson-Thomas invariants and cluster
        transformations}, preprint, 2008.
\bibitem{MarkRau} H.~Markwig, J.~Rau: \emph{Tropical
        descendent Gromov-Witten invariants}, preprint, 2008.
\bibitem{mikhalkin} G.\ Mikhalkin:
        \emph{Enumerative tropical  algebraic geometry in $\RR^2$},
        J.\ Amer.\ Math.\ Soc.\ \textbf{18}  (2005), 313--377.
\bibitem{Nishinou} T.\ Nishinou: \emph{Disc counting on toric
        varieties via tropical curves,} preprint (2006).
\bibitem{nisi} T.\ Nishinou, B.\ Siebert:
        \emph{Toric degenerations of toric varieties and tropical
        curves},
        Duke Math.\ J.~\textbf{135} (2006), 1--51.
\bibitem{Sabbah} C.~Sabbah: \emph{D\' eformations isomonodromiques et 
        vari\' et\' es de Frobenius}, Savoirs Actuels (Les Ulis). 
        Math\' ematiques (Les Ulis), EDP Sciences, Les Ulis; CNRS \' Editions,          Paris, 2002. xvi+289 pp.
\bibitem{Sturmfels} J.~Richter-Gebert, B.~Sturmfels, T.~Theobald:
        \emph{First steps in tropical geometry}, in
        Idempotent mathematics and mathematical physics,  289--317, Contemp. 
        Math., 377, Amer. Math. Soc., Providence, RI, 2005. 
\bibitem{Saito} K.~Saito: \emph{Period mapping associated to a primitive
        form,} Publ.\ RIMS, Kyoto Univ.\ {\bf 19}, (1983), 1231-1264.
\bibitem{SYZ} A.\ Strominger, S.-T.\ Yau, and E.~Zaslow:
        \emph{Mirror Symmetry is $T$-duality},
        Nucl.\ Phys.~\textbf{B479} (1996), 243--259.
\end{thebibliography}
\end{document}